\documentclass[11pt]{amsart}
\usepackage{amssymb}
\usepackage{amsmath}
\usepackage{fancyhdr}
\usepackage[british]{babel}
\usepackage{geometry}
\usepackage{enumitem}
\usepackage{graphicx}
\usepackage[font={footnotesize}]{caption}
\usepackage[usenames,dvipsnames]{xcolor}
\usepackage{pstricks,tikz}
\usetikzlibrary{positioning,chains,fit,shapes,calc,patterns}
\usepackage[h]{esvect}
\usepackage[
		bookmarksopen=true,
		bookmarksopenlevel=1,
		colorlinks=true,
		linkcolor=darkblue, 
        linktoc=page,
		citecolor=darkblue,
]{hyperref}

\title[]{On the decomposition threshold of a given graph}

\date{\today}
\author[S.~Glock, D.~K\"uhn, A.~Lo, R.~Montgomery and D.~Osthus]{Stefan Glock, Daniela K\"uhn, Allan Lo, Richard Montgomery and Deryk Osthus}
\thanks{The research leading to these results was partially supported by the EPSRC, grant no. EP/M009408/1 (D.~K\"uhn and D.~Osthus), 
by the Royal Society and the Wolfson Foundation (D.~K\"uhn) as well as by the European Research Council
under the European Union's Seventh Framework Programme (FP/2007--2013) / ERC Grant
Agreements no. 258345 (D.~K\"uhn, R.~Montgomery) and 306349 (S.~Glock and D.~Osthus).}

\newtheorem{firstthm}{Proposition}[section]
\newtheorem{theorem}[firstthm]{Theorem}
\newtheorem{prop}[firstthm]{Proposition}
\newtheorem{lemma}[firstthm]{Lemma}
\newtheorem{cor}[firstthm]{Corollary}

\newtheorem{defin}[firstthm]{Definition}
\newtheorem{conj}[firstthm]{Conjecture}

\newtheorem{fact}[firstthm]{Fact}

\numberwithin{equation}{section}
\definecolor{darkblue}{rgb}{0,0,0.5}
\def\noproof{{\unskip\nobreak\hfill\penalty50\hskip2em\hbox{}\nobreak\hfill%
       $\square$\parfillskip=0pt\finalhyphendemerits=0\par}\goodbreak}
\def\endproof{\noproof\bigskip}
\newdimen\margin   
\def\textno#1&#2\par{
   \margin=\hsize
   \advance\margin by -4\parindent
          \setbox1=\hbox{\sl#1}
   \ifdim\wd1 < \margin
      $$\box1\eqno#2$$
   \else
      \bigbreak
      \hbox to \hsize{\indent$\vcenter{\advance\hsize by -3\parindent
      \it\noindent#1}\hfil#2$}
      \bigbreak
   \fi}
\def\proof{\removelastskip\penalty55\medskip\noindent{\bf Proof. }}
\def\lateproof#1{\removelastskip\penalty55\medskip\noindent{\bf Proof of #1. }}

\begin{document}

\def\COMMENT#1{}
\def\TASK#1{}

\def\eps{{\varepsilon}}
\newcommand{\ex}{\mathbb{E}}
\newcommand{\pr}{\mathbb{P}}
\newcommand{\cB}{\mathcal{B}}
\newcommand{\cE}{\mathcal{E}}
\newcommand{\cS}{\mathcal{S}}
\newcommand{\cF}{\mathcal{F}}
\newcommand{\cH}{\mathcal{H}}
\newcommand{\cC}{\mathcal{C}}
\newcommand{\cM}{\mathcal{M}}
\newcommand{\bN}{\mathbb{N}}
\newcommand{\bR}{\mathbb{R}}
\def\O{\mathcal{O}}
\newcommand{\cP}{\mathcal{P}}
\newcommand{\cQ}{\mathcal{Q}}
\newcommand{\cR}{\mathcal{R}}
\newcommand{\cK}{\mathcal{K}}
\newcommand{\cD}{\mathcal{D}}
\newcommand{\cI}{\mathcal{I}}
\newcommand{\cV}{\mathcal{V}}
\newcommand{\cT}{\mathcal{T}}
\newcommand{\1}{{\bf 1}_{n\not\equiv \delta}}
\newcommand{\eul}{{\rm e}}
\newcommand{\Erd}{Erd\H{o}s}
\newcommand{\cupdot}{\mathbin{\mathaccent\cdot\cup}}

\newcommand{\defn}{\emph}

\newcommand\restrict[1]{\raisebox{-.5ex}{$|$}_{#1}}

\newcommand{\prob}[1]{\mathrm{\mathbb{P}}(#1)}
\newcommand{\expn}{\mathrm{\mathbb{E}}}
\def\gnp{G_{n,p}}
\def\G{\mathcal{G}}
\def\lflr{\left\lfloor}
\def\rflr{\right\rfloor}
\def\lcl{\left\lceil}
\def\rcl{\right\rceil}

\newcommand{\brackets}[1]{\left(#1\right)}
\def\sm{\setminus}
\newcommand{\Set}[1]{\{#1\}}
\newcommand{\set}[2]{\{#1\,:\;#2\}}
\def\In{\subseteq}

\begin{abstract}  \noindent
We study the $F$-decomposition threshold $\delta_F$ for a given graph $F$.
Here an $F$-decomposition of a graph $G$ is a collection of edge-disjoint copies of $F$ in $G$ which together cover every edge of $G$. (Such an $F$-decomposition can only exist if $G$ is $F$-divisible, i.e.~if $e(F)\mid e(G)$ and each vertex degree of $G$ can be expressed as a linear combination of the vertex degrees of $F$.)

The $F$-decomposition threshold $\delta_F$ is the smallest value ensuring that an $F$-divisible graph $G$ on $n$ vertices with $\delta(G)\ge(\delta_F+o(1))n$ has an $F$-decomposition.

Our main results imply the following for a given graph $F$, where $\delta_F^\ast$ is the fractional version of $\delta_F$ and $\chi:=\chi(F)$:
\begin{enumerate}[label=(\roman*)]
\item $\delta_F\le \max\Set{\delta_F^\ast,1-1/(\chi+1)}$;
\item if $\chi\ge 5$, then $\delta_F\in\Set{\delta_F^{\ast},1-1/\chi,1-1/(\chi+1)}$;
\item we determine $\delta_F$ if $F$ is bipartite.
\end{enumerate}
In particular, (i) implies that $\delta_{K_r}=\delta^\ast_{K_r}$. Our proof involves further developments of the recent `iterative' absorbing approach.
\end{abstract}

\maketitle

\section{Introduction} \label{sec:intro}

Let $F$ be a fixed graph. A fundamental theorem of Wilson~\cite{W} states that for all sufficiently large $n$, the complete graph $K_n$ has an $F$-decomposition (subject to the divisibility conditions discussed below). Here, an \defn{$F$-decomposition} of a graph $G$ is a collection of edge-disjoint copies of $F$ in $G$ which together cover every edge of $G$. The case when $F$ is a triangle is known as Kirkman's theorem \cite{Ki}.

The problem of determining whether an arbitrary graph $G$ has an $F$-decomposition is much more difficult (in fact, the corresponding decision problem is NP-complete (see~\cite{DT})). Recently there has been some significant progress in extending Wilson's theorem to dense graphs, i.e.~graphs of large minimum degree --- the current paper will build on this.

A clearly necessary condition for the existence of an $F$-decomposition is that $e(F)\mid e(G)$. If this is satisfied, then we say that $G$ is \defn{$F$-edge-divisible}. Moreover, for $r\in \bN$, we call $G$ \defn{$r$-divisible}, if $r \mid d_G(x)$ for all $x\in V(G)$. We say that $G$ is \defn{$F$-degree-divisible} if it is $gcd(F)$-divisible, where $gcd(F):=gcd\set{d_F(v)}{v\in V(F)}$ (this is another trivially necessary condition for the existence of an $F$-decomposition). If a graph $G$ is both $F$-edge-divisible and $F$-degree-divisible, then we simply say that $G$ is \defn{$F$-divisible}. 

For a fixed graph $F$, let $\delta_F$ be the minimum of the set of all non-negative real numbers $\delta$ with the following property: for all $\mu > 0$ there exists an $n_0\in\bN$ such that whenever $G$ is an $F$-divisible graph on $n\geq n_0$ vertices with $\delta(G)\geq(\delta+\mu)n$, then $G$ has an $F$-decomposition. Clearly, the minimum exists. Note that isolated vertices in $F$ are irrelevant and the $F$-decomposition problem is trivial if $F$ has only one edge. Thus in all our statements concerning a given graph $F$, we will assume that $F$ has no isolated vertices and $e(F)\ge 2$ without mentioning this explicitly.

The purpose of this paper is to investigate the above $F$-decomposition threshold $\delta_F$. In particular, we determine $\delta_F$ for all bipartite graphs, improve existing bounds for general $F$ and prove a `discretisation' result for the possible values of $\delta_F$.\COMMENT{We also note that it is not restrictive to require $\delta$ being non-negative in the above definition. $F=K_2$ is the only graph for which every ($F$-divisible) graph has an $F$-decomposition. Since $F=K_2$ is a trivial case for our problem, we will always assume that $F$ has at least two edges. Then, no $\delta<0$ would work in the above definition.
If $\chi(F)\ge 3$, it is clear, so assume $\chi(F)=2$. If $F$ is connected, we have $\delta_F\ge 1/2$. So assume $F$ has more than one component. An extremal example on $n$ vertices could then look as follows: Take $F$ and identify two vertices from different components. Add a lot of isolated vertices to reach $n$ vertices.}

\subsection{Bounding the decomposition threshold for arbitrary graphs}
Our first main result (Theorem~\ref{thm:main}) bounds the decomposition threshold $\delta_F$ in terms of the approximate decomposition threshold $\delta_F^{0+}$, the fractional decomposition threshold $\delta_F^\ast$, and the threshold $\delta_F^e$ for covering a given edge. We now introduce these formally.

Let $F$ be a fixed graph.
For $\eta \ge 0$, an \defn{$\eta$-approximate $F$-decomposition} of an $n$-vertex graph $G$ is a collection of edge-disjoint copies of $F$ contained in $G$ which together cover all but at most $\eta n^2$ edges of $G$. Let $\delta^{\eta}_F$ be the smallest $\delta\ge 0$ such that for all $\mu > 0$ there exists an $n_0\in\bN$ such that whenever $G$ is a graph on $n\geq n_0$ vertices with $\delta(G)\geq(\delta+\mu)n$, then $G$ has an $\eta$-approximate $F$-decomposition.
Clearly, $\delta_F^{\eta'} \ge \delta_F^{\eta}$ whenever $\eta'\le \eta$. We let $\delta_F^{0+}:=\sup_{\eta >0}\delta_F^{\eta}$.

Let $G^F$ be the set of copies of~$F$ in $G$.
A \defn{fractional $F$-decomposition of $G$} is a function $\omega \colon G^F \to [0,1]$ such that, for each $e \in E(G)$,
\begin{equation} \label{packing}
\sum_{F' \in G^F \colon e \in E(F')} \omega(F') = 1.
\end{equation}
Note that every $F$-decomposition is a fractional $F$-decomposition where $\omega(F) \in \{0, 1\}$.

Let $\delta^{\ast}_F$ be the smallest $\delta\ge 0$ such that for all $\mu > 0$ there exists an $n_0\in\bN$ such that whenever $G$ is an $F$-divisible graph on $n\geq n_0$ vertices with $\delta(G)\geq(\delta+\mu)n$, then $G$ has a fractional $F$-decomposition. Usually the definition considers all graphs $G$ (and not only those which are $F$-divisible) but it is convenient for us to make this additional restriction here as $\delta_F^\ast$ is exactly the relevant parameter when investigating $\delta_F$ (in particular, we trivially have $\delta_F^\ast\le \delta_F$). Haxell and R\"odl~\cite{HR} used Szemer\'edi's regularity lemma to show that a fractional $F$-decomposition of a graph $G$ can be turned into an approximate $F$-decomposition of $G$ (see Theorem~\ref{thm:Haxell-Rodl}). This implies $\delta_F^{0+}\le \delta_F^\ast$ (see Corollary~\ref{cor:threshold relations}).

Let $\delta^e_F$ be the smallest $\delta$ such that for all $\mu > 0$ there exists an $n_0\in\bN$ such that whenever $G$ is a graph on $n\geq n_0$ vertices with $\delta(G)\geq(\delta+\mu)n$, and $e'$ is an edge in $G$, then $G$ contains a copy of $F$ which contains~$e'$.

Our first result bounds $\delta_F$ in terms of the approximate decomposition threshold $\delta_F^{0+}$ and the chromatic number of $F$. Parts (ii) and (iii) give much more precise information if $\chi\ge 5$. We obtain a `discretisation result' in terms of the parameters introduced above. We do not believe that this result extends to $\chi=3,4$ (see Section~\ref{sec:conclusion} for a further discussion). On the other hand, we do have $\delta_F\in \Set{0,1/2,2/3}$ if $\chi(F)=2$ (see Section~\ref{subsec:bip}).
We also believe that none of the terms in the discretisation statement can be omitted.

\begin{theorem} \label{thm:main}
Let $F$ be a graph with $\chi:=\chi(F)$.
\begin{enumerate} [label=(\roman*)]
\item Then $\delta_F\le \max\Set{\delta_F^{0+},1-1/(\chi+1)}$.
\item If $\chi\ge 5$, then $\delta_F\in\Set{\max\Set{\delta_F^{0+},\delta_F^{e}},1-1/\chi,1-1/(\chi+1)}$.
\item If $\chi\ge 5$, then $\delta_F\in\Set{\delta_F^\ast,1-1/\chi,1-1/(\chi+1)}$.
\end{enumerate}
\end{theorem}

Theorem~\ref{thm:main}(i) improves a bound of $\delta_F\le \max\Set{\delta_F^{0+},1-1/3r}$ proved in \cite{BKLO} for $r$-regular graphs $F$. Also, the cases where $F=K_3$ or $C_4$ of (i) were already proved in \cite{BKLO}.

Since it is known that $\delta_{K_r}^{0+} \ge 1-1/(r+1)$ (see e.g.~\cite{Y05}),
Theorem~\ref{thm:main} implies that the decomposition threshold for cliques equals its fractional relaxation.

\begin{cor} \label{cor:int=rel}
For all $r\ge 3$, $\delta_{K_r}=\delta_{K_r}^\ast=\delta_{K_r}^{0+}$.
\end{cor}

\subsection{Explicit bounds} \label{subsec:explicit bounds}

Theorem~\ref{thm:main} involves several `auxiliary thresholds' and parameters that play a role in the construction of an $F$-decomposition. Bounds on these of course lead to better `explicit' bounds on $\delta_F$ which we now discuss.

The central conjecture in the area is due to Nash-Williams~\cite{NW} (for the triangle case) and Gustavsson~\cite{G} (for the general case).

\begin{conj}[Gustavsson~\cite{G}, Nash-Williams~\cite{NW}]\label{conj:Gustavsson}
For every $r\ge 3$, there exists an $n_0=n_0(r)$ such that every $K_r$-divisible graph $G$ on $n\ge n_0$ vertices with $\delta(G)\ge (1-1/(r+1))n$ has a $K_r$-decomposition.
\end{conj}

For general $F$, the following conjecture provides a natural upper bound for $\delta_F$ which would be best possible for the case of cliques.
It is not clear to us what a formula for general~$F$ might look like.

\begin{conj}\label{conj:general upper bound}
For all graphs $F$, $\delta_F\le 1-1/(\chi(F)+1)$.
\end{conj}

Note that by Theorem~\ref{thm:main} in order to prove Conjecture~\ref{conj:general upper bound} it suffices to show $\delta_F^{0+}\le 1-1/(\chi(F)+1)$.
This in turn implies that Conjecture~\ref{conj:general upper bound} is actually a special case of Conjecture~\ref{conj:Gustavsson}. Indeed, it follows from a result of Yuster~\cite{Y12} that for every graph $F$, $\delta_F^{0+}\le \delta_{K_{\chi(F)}}^{0+}$, and thus $\delta_F^{0+}\le \delta_{K_{\chi(F)}}^{\ast}\le \delta_{K_{\chi(F)}}$.

In view of this, bounds on $\delta_{K_r}^\ast$ are of considerable interest. The following result gives the best bound for general $r$ (see~\cite{BKLMO}) and triangles (see~\cite{D}).

\begin{theorem}[\cite{BKLMO}, \cite{D}]\label{thm:fractional} $\ $
\begin{itemize}
\item[(i)] For every $r\ge 3$, we have $\delta_{K_r}^\ast\le 1-10^{-4}r^{-3/2}$.
\item[(ii)] $\delta_{K_3}^\ast\le 9/10$.
\end{itemize}
\end{theorem}

This improved earlier bounds by Yuster~\cite{Y05} and Dukes~\cite{Du}. Together with the results in~\cite{BKLO},
part (ii) implies $\delta_{K_3}\le 9/10$. More generally, combining Theorem~\ref{thm:fractional} and Theorem~\ref{thm:main}(i) with the fact that $\delta_F^{0+}\le \delta_{K_{\chi(F)}}^{0+}\le \delta_{K_{\chi(F)}}^{\ast}$, one obtains the following explicit upper bound on the decomposition threshold.

\begin{cor}\label{cor:explit bounds}$\ $
\begin{enumerate}[label=(\roman*)]
\item For every graph $F$, $\delta_F\le 1-10^{-4}\chi(F)^{-3/2}$.
\item If $\chi(F)=3$, then $\delta_F\le 9/10$.
\end{enumerate}
\end{cor}

Here, (i) improves a bound of $1-1/\max\Set{10^{4}\chi(F)^{3/2},6e(F)}$ obtained by combining the results of \cite{BKLMO} and \cite{BKLO} (see~\cite{BKLMO}). It also improves earlier bounds by Gustavsson~\cite{G} and Yuster~\cite{Y05,Y14}. A bound of $1-\eps$ also follows from the
results of Keevash~\cite{Ke}.

In the $r$-partite setting an analogue of Corollary~\ref{cor:int=rel} was proved in~\cite{BKLOT}, an analogue of Theorem~\ref{thm:fractional}(i) (with weaker bounds) in~\cite{M} and an analogue of Theorem~\ref{thm:fractional}(ii) (again with weaker bounds) in~\cite{Du15}.
These bounds can be combined to give results on the completion of (mutually orthogonal) partially filled in Latin squares.
Moreover, it turns out that if $\delta_F>\delta^\ast_F$ (in the non-partite setting), then there exist extremal graphs that are extremely close to large complete partite graphs, which adds further relevance to results on the $r$-partite setting (see Section~\ref{sec:conclusion}).

\subsection{Decompositions into bipartite graphs}\label{subsec:bip}

Let $F$ be a bipartite graph. Yuster~\cite{Y02} showed that $\delta_F=1/2$ if $F$ is connected and contains a vertex of degree one. Moreover, Barber, K\"uhn, Lo and Osthus~\cite{BKLO} showed that $\delta_{C_4}=2/3$ and $\delta_{C_\ell}=1/2$ for all even $\ell\ge 6$ (which improved a bound of $\delta_{C_4}\le 31/32$ by Bryant and Cavenagh~\cite{BC}). Here we generalise these results to arbitrary bipartite graphs.

Note that if $F$ is bipartite, it is easy to see that $\delta_F^{0+}=0$, as one can trivially obtain an approximate decomposition via repeated applications of the \Erd-Stone theorem. This allows us to determine $\delta_F$ for any bipartite graph $F$. To state the result, we need the following definitions. A set $X\In V(F)$ is called \defn{$C_4$-supporting in $F$} if there exist distinct $a,b\in X$ and $c,d\in V(F)\sm X$ such that $ac,bd,cd\in E(F)$. We define 
\begin{align*}
\tau(F)&:=gcd\set{e(F[X])}{X\In V(F)\text{ is not }C_4\text{-supporting in }F},\\
\tilde{\tau}(F)&:=gcd\set{e(C)}{C\text{ is a component of }F}.
\end{align*} So for example $\tau(F)=1$ if there exists an edge in $F$ that is not contained in any cycle of length $4$, and $\tilde{\tau}(F)>1$ if $F$ is connected (and $e(F)\ge 2$).
The definition of $\tau$ can be motivated by considering the following graph $G$: Let $A,B,C$ be sets of size $n/3$ with $G[A],G[C]$ complete, $B$ independent and $G[A,B]$ and $G[B,C]$ complete bipartite. Note that $\delta(G)\sim 2n/3$. It turns out that the extremal examples which
we construct showing $\delta_F\ge 2/3$ for certain bipartite graphs $F$ are all similar to $G$ (see Proposition~\ref{prop:exex bip 1}). Moreover, $\tau(F)=1$ if for any large $c$ there is a set of copies of $F$ in $G$ whose number of edges in $G[A]$ add up to $c$.

We note that $\tau(F)\mid gcd(F)$ and $gcd(F)\mid \tilde{\tau}(F)$ (see Fact~\ref{fact:bip par rel}). The following theorem determines $\delta_F$ for every bipartite graph $F$. 

\begin{theorem} \label{thm:bipartite char}
Let $F$ be a bipartite graph. Then
$$\delta_F = \begin{cases}
2/3 &\mbox{if } \tau(F)>1; \\
0 &\mbox{if } \tilde{\tau}(F)=1 \mbox{ and } F \mbox{ has a bridge}; \\
1/2 &\mbox{otherwise. }
\end{cases}$$
\end{theorem}

We will prove Theorem~\ref{thm:bipartite char} 
in Section~\ref{sec:bip}.
The next corollary translates Theorem~\ref{thm:bipartite char} into explicit results for important classes of bipartite graphs.

\begin{cor} \label{cor:bip char}
The following hold.
\begin{enumerate}[label=(\roman*)]
\item Let $s,t\in \bN$ with $s+t>2$. Then $\delta_{K_{s,t}}=1/2$ if $s$ and $t$ are coprime and $\delta_{K_{s,t}}=2/3$ otherwise.
\item If $gcd(F)=1$ and $F$ is connected, then $\delta_F=1/2$.
\item If $F$ is connected and has an edge that is not contained in any cycle of length $4$, then $\delta_F=1/2$.
\end{enumerate}
\end{cor}
(For (ii) and (iii) recall that we always assume $e(F)\ge 2$.) Note that $\tau(K_{s,t})=gcd(s,t)$. Then (i)--(iii) follow from the definitions of $\tau$ and $\tilde{\tau}$.

\subsection{Near-optimal decompositions}\label{subsec:intro near optimal}

Along the way to proving Theorem~\ref{thm:main} we obtain the following bound guaranteeing a `near-optimal' decomposition.
For this, let $\delta^{vx}_F$ be the smallest $\delta\ge 0$ such that for all $\mu > 0$ there exists an $n_0\in\bN$ such that whenever $G$ is a graph on $n\geq n_0$ vertices with $\delta(G)\geq(\delta+\mu)n$, and $x$ is a vertex of $G$ with $gcd(F)\mid d_G(x)$, then $G$ contains a collection $\cF$ of edge-disjoint copies of $F$ such that $\set{xy}{y\in N_G(x)} \In \bigcup \cF$.
Loosely speaking, $\delta^{vx}_F$ is the threshold that allows us to cover all edges at one vertex. For example, if $F$ is a triangle, then $\delta_F^{vx}$ is essentially the threshold that $N_G(x)$ contains a perfect matching whenever $d_G(x)$ is even. Note that $\delta_F^{vx}\ge \delta_F^e$.

The following theorem roughly says that if we do not require to cover all edges of $G$ with edge-disjoint copies of $F$, but accept a bounded number of uncovered edges, then the minimum degree required can be less than if we need to cover all edges.

\begin{theorem} \label{thm:almost cover}
For any graph $F$ and $\mu >0$ there exists a constant $C=C(F,\mu)$ such that whenever $G$ is an $F$-degree-divisible graph on $n$ vertices satisfying $$\delta(G)\geq (\max\Set{\delta_F^{0+},\delta_F^{vx}}+\mu)n$$ then $G$ contains a collection of edge-disjoint copies of $F$ covering all but at most $C$ edges. 
\end{theorem}

Here, $\delta_F^{vx}\le 1-1/\chi(F)$ (see Corollary~\ref{cor:deltafa}). For many bipartite graphs $F$, e.g.~trees and complete balanced bipartite graphs, our results imply that $\max\Set{\delta_F^{0+},\delta_F^{vx}}<\delta_F$. It seems plausible to believe that there also exist graphs $F$ with $\chi(F)\ge 3$ such that $\max\Set{\delta_F^{0+},\delta_F^{vx}}<\delta_F$. However, the current bounds on $\delta_F^{0+}$ do not suffice to verify this. The proof of Theorem~\ref{thm:almost cover} can be found in Section~\ref{sec:near optimal}.

\subsection{Recent developments}

During the refereeing process of this paper, there have been quite a number of developments in the area of $F$-decompositions, which we briefly mention here.
Mostly relevant to the explicit bounds discussed in Section~\ref{subsec:explicit bounds}, we remark that the fractional decomposition threshold of cliques has been improved in~\cite{montgomery:ta}, showing that $\delta_{K_r}\le 1-1/(100r)$. Together with Theorem~\ref{thm:main}, this immediately improves the bound in Corollary~\ref{cor:explit bounds}(i) to $\delta_F\le 1-1/(100\chi(F))$ for every graph~$F$.

Our proof method, iterative absorption, has been further developed in the hypergraph setting to give a new proof~\cite{GKLO:16} of the existence conjecture for combinatorial designs (first proved by Keevash~\cite{Ke}).
The results in~\cite{GKLO:16} also give effective bounds in the minimum degree setting.
They were further extended to decompositions into arbitrary uniform hypergraphs~$F$
\cite{GKLO:17}, thus generalizing Wilson's theorem to hypergraphs. An alternative proof of the existence of such hypergraph $F$-decompositions was subsequently obtained in~\cite{keevash:18b}.

Building on some of our tools, Taylor~\cite{taylor:19} determined the exact minimum degree threshold for $C_\ell$-decompositions when $\ell =4$
(namely $2n/3-1$)
or when $\ell \ge 8$ is even (namely $n/2$). It would be interesting to obtain further 
such exact results.

For a short and self-contained exposition of the iterative absorption method in the setting of triangle decompositions, see~\cite{BGKLMO:18}.

\section{Overview of the proofs and organisation of the paper}  \label{sec:sketch}
One key ingredient in the proofs of Theorems~\ref{thm:main}, \ref{thm:bipartite char} and~\ref{thm:almost cover} is an iterative absorption method.
Very roughly, this means that we build our decomposition of a given graph $G$ in many iterations, 
where in each iteration we add copies of $F$ to our current partial decomposition. 
In the current proof, we can carry out this iteration until  we have a `near-optimal decomposition' which covers 
all but a bounded number of edges of $G$. 
Let $H$ be the graph consisting of the leftover (i.e.~uncovered) edges. 
This leftover graph $H$ can then be absorbed into a graph $A$ which we set aside at the beginning 
(i.e.~$H \cup A$ has an $F$-decomposition). Altogether this yields an $F$-decomposition of the original graph $G$.

More precisely, to obtain the near-optimal decomposition we proceed as follows.
At the beginning of the proof we will fix a suitable nested sequence of vertex sets
$V(G)=U_0 \supseteq U_1 \supseteq \dots \supseteq U_\ell$, which will be called a vortex in $G$.
After the $i$th iteration we can ensure that the uncovered edges all lie in $U_i$, which is much smaller than $U_{i-1}$.
We can also preserve the relative density of the leftover graph $G_i$ after the $i$th iteration, 
i.e.~$\delta(G_i[U_i])/|U_i| \sim \delta$, where $\delta:=\delta(G)/|G|$.
We will show that this can be achieved provided that $\delta \ge \delta_F^{0+}$ and $\delta \ge \delta_F^{vx}$.
These iterative steps are carried out in Section~\ref{sec:near optimal}. 
In particular, we obtain Theorem~\ref{thm:almost cover} as a byproduct of this iterative absorption argument.

We now turn to the absorption step itself. The final set $U_\ell$ in the iteration will have bounded size.
This immediately implies that the final leftover graph $H\In G[U_\ell]$ will also have a bounded number of edges.
In particular, there are only a bounded number $H_1,\dots,H_s$ of possibilities for $H$.
The graph $A$ will be constructed as the edge-disjoint union of absorbers $A_1, \dots, A_s$, where each $A_i$ is tailored towards $H_i$.
More precisely, the crucial property is that both $A_i$ and $A_i \cup H_i$ have an $F$-decomposition for each $i \in [s]$.
With this property, it is clear that $A$ has the required absorbing property, i.e.~$A \cup H_i$ has an $F$-decomposition for any of the 
permissible leftovers $H_i$.
 
The absorbers will be constructed in several steps: 
rather than constructing $A_i$ directly, we will obtain it as the `concatenation' 
(equivalent to the edge-disjoint union) of several `transformers' $T$.
The role of $T$ is  to transform $H_i$ into a suitable different graph $H_i'$
(more precisely, both $H_i' \cup T$ and $T \cup H_i$
have an $F$-decomposition).
We can then concatenate several such transformers to transform $H_i$ into a disjoint union of copies of $F$,
which trivially has an $F$-decomposition.

This reduces the absorption problem to that of constructing transformers.
Surprisingly, the main hurdle for the latter is the ability to construct a transformer $T$ which simply 
moves $H_i$ to a different position, i.e.~transforms $H_i$ into an isomorphic copy $H_i'$ of $H_i$ in $G$,
with a different vertex set. Once this is achieved, we can obtain more general transformers by simple modifications.

Yet again, we do not construct these transformers directly, but construct them from building blocks called `switchers'.
These switchers are transformers with more limited capabilities. 
The most important switchers are $C_6$-switchers and $K_{2,r}$-switchers.
A $C_6$-switcher $S$ transforms the perfect matching $E^+:=\{u_1u_2,u_3u_4,u_5u_6\}$
into its `complement' $E^-:=\{u_2u_3,u_4u_5,u_6u_1\}$ along a $6$-cycle.
(The formal requirement is that both $S \cup E^+$ and $S \cup E^-$ have an $F$-decomposition.)
A $K_{2,r}$-switcher transforms a star with $r$ leaves centred at $x$ into a star with the same leaves centred at $x'$.
Surprisingly, it turns out that these building blocks suffice to build the desired transformers (see Lemma~\ref{lem:switch2transform}). 

Occasionally, we build the above switchers from even more `basic' ones. For example, in Section~\ref{sec:switchers} we will build a $C_6$-switcher by combining $C_4$-switchers in a suitable way.
Apart from proving the existence of switchers, 
we also need to be able to find them in $G$.
This is where we may need the condition that $\delta(G) \ge (1-1/(\chi+1)+o(1))|G|$.
To achieve this, we will apply Szemer\'edi's regularity lemma to $G$ to obtain its reduced graph $R$.
We will then find a `compressed' version (i.e.~a suitable homomorphism) of the switcher in $R$.
This then translates to the existence of the desired switcher in $G$ via standard regularity techniques.
 
The switchers are also key to our discretisation results in Theorem~\ref{thm:main}(ii) and~(iii). We show that if $\delta_F<1-1/(\chi+1)$, then to find the relevant switchers (and hence, as described above, the relevant absorbers) we need the graph $G$ only to have minimum degree $(1-1/\chi+o(1))|G|$. Similarly, if $\delta_F<1-1/\chi$, the minimum degree we require is only $(1-1/(\chi-1)+o(1))|G|$. As discussed earlier we require the minimum degree to be at least $(\max\Set{\delta_F^{0+},\delta_F^{vx}}+o(1))|G|$ in order to iteratively cover all but a constant number of edges in $G$ (see Theorem~\ref{thm:almost cover}). This may not be sufficiently high to construct our absorbers, but this discretisation argument will allow us to conclude that if $\delta_F$ exceeds $\max\Set{\delta_F^{0+},\delta_F^{vx}}$ then it can take at most two other values, $1-1/(\chi+1)$ or $1-1/\chi$ (see Theorem~\ref{thm:almost main}).

Most of the above steps are common to the proof of Theorems~\ref{thm:main} and~\ref{thm:bipartite char},
i.e.~we can prove them in a unified way. The key additional difficulty in the bipartite case is proving the existence of a 
$C_6$-switcher for those $F$ with $\delta_F=1/2$.

The iterative absorption approach was initially introduced in \cite{KKO} and was further developed in the present context of
$F$-decompositions in \cite{BKLO}.
The present iteration procedure is much simpler than the one in \cite{BKLO}. 
The concept of transformers also originates in \cite{BKLO}, but the approach via switchers is a new feature which allows us to go
significantly beyond the results in \cite{BKLO}.

This paper is organised as follows. The following section contains the basic notation that we use. Sections~\ref{sec:vortex} and~\ref{sec:near optimal} deal with the `near-optimal decomposition'. More precisely, Section~\ref{sec:vortex} introduces the concept of vortices, and in Section~\ref{sec:near optimal} we perform the iteration based on these vortices, leading to the proof of Theorem~\ref{thm:almost cover}. 

Sections~\ref{sec:regularity}--\ref{sec:switchers} and~\ref{sec:absorbers} are devoted to absorbers, which are constructed in several steps as described above. Section~\ref{sec:regularity} recalls well-known results on $\eps$-regularity and introduces the setting of an $(\alpha,\eps,k)$-partition in which we intend to find absorbers. In Section~\ref{sec:general dec}, we combine the `near-optimal decomposition' with the concept of absorbing in order to prove a general decomposition theorem. In order to apply this theorem, we need to be able to find absorbers in a graph with an $(\alpha,\eps,k)$-partition. In Section~\ref{sec:compressions}, we introduce the concepts that allow us to achieve this. In Section~\ref{sec:transform} we will construct transformers out of switchers. Section~\ref{sec:switchers} deals with the construction of switchers. Finally, in Section~\ref{sec:absorbers} we use transformers to build absorbers, which we can then use to obtain upper bounds on $\delta_F$ using the general decomposition theorem. In addition, in Section~\ref{sec:lower bounds}, we prove some relations between the auxiliary thresholds which we need for our discretisation result. 

In Section~\ref{sec:bip}, we determine $\delta_F$ for all bipartite graphs $F$. Finally, in Section~\ref{sec:deltafa} we will investigate $\delta_F^{vx}$ and then combine all our results to prove Theorem~\ref{thm:main}.

\section{Notation and tools} \label{sec:notation}

For a graph $G$, we let $|G|$ denote the number of vertices of $G$, $e(G)$ the number of edges of $G$, and $\overline{G}$ the complement of $G$. For a vertex $v\in V(G)$, we write $N_G(v)$ for the neighbourhood of $v$ and $d_G(v)$ for its degree. More generally, for vertex subsets $S,V\In V(G)$, we let $N_G(S,V):=V\cap \bigcap_{v\in S}N_G(v)$ denote the set of vertices in $V$ that are adjacent to all vertices in $S$, where $N_G(\emptyset,V):=V$. We then define $d_G(S,V):=|N_G(S,V)|$ and write $d_G(x,V)$ if $S=\Set{x}$. 

A set of edge-disjoint copies of $F$ in a graph $G$ is called an \defn{$F$-collection}. Thus an $F$-decomposition of $G$ is an $F$-collection that covers every edge of $G$.
If $G$ and $H$ are two edge-disjoint graphs, we write $G\cupdot H$ for the union of $G$ and $H$.

Let $G$ be a graph and let $V_1,\dots,V_k$ be disjoint subsets of $V(G)$. We write $G[V_1]$ for the subgraph of $G$ induced by $V_1$. Moreover, if $k\ge 2$, then $G[V_1,\dots,V_k]$ denotes the $k$-partite subgraph of $G$ induced by $V_1,\dots,V_k$. If $V_1,\dots,V_k$ form a partition $\cP$ of $V(G)$, we write $G[\cP]$ instead of $G[V_1,\dots,V_k]$.

If $G$ is a graph and $H$ is a subgraph of $G$, then $G-H$ denotes the graph with vertex set $V(G)$ and edge set $E(G)\sm E(H)$. Moreover, if $X\In V(G)$, then $G-X:=G[V(G)\sm X]$.

For a graph $K$ and $t\in\bN$, we let $K(t)$ denote the graph obtained by replacing each vertex of $K$ with $t$ vertices and each edge of $K$ with a copy of $K_{t,t}$.
A homomorphism $\phi\colon H\to G$ from a graph $H$ to a graph $G$ is a map $\phi\colon V(H)\to V(G)$ such that $\phi(x)\phi(y)\in E(G)$ whenever $xy\in E(H)$. We let $\phi(H)$ denote the subgraph of $G$ with vertex set $\phi(V(H))$ and edge set $\set{\phi(x)\phi(y)}{xy\in E(H)}$.
We say that $\phi$ is \defn{edge-bijective} if $e(G)=e(\phi(H))=e(H)$.
We will sometimes identify a graph with its edge set if this enhances readability and does not affect the argument.
%

For $n\in\bN$ and distinct $i_1,\dots,i_k\in[n]$, we let $(i_1,\dots,i_k)$ denote the permutation $\pi\colon[n]\rightarrow[n]$ defined as $\pi(i_j):=i_{j+1}$ for all $j\in[k-1]$, $\pi(i_k):=i_1$ and $\pi(i):=i$ for all $i\in[n]\sm\Set{i_1,\dots,i_k}$.
We write $2^U$ for the power set of $U$.

We write $x\ll y$ to mean that for any $y\in (0,1]$ there exists an $x_0\in (0,1)$ such that for all $x\le x_0$ the subsequent statement holds. Hierarchies with more constants are defined in a similar way and are to be read from the right to the left. The expression $a= b \pm c$ means that $a\in [b-c,b+c]$.\COMMENT{Should we mention $[n]$ contains the natural numbers from $1$ to $n$ (and $[0]$ is the empty set)?}

Let $m,n,N\in \bN$ with $\max\Set{m,n}< N$. Recall that a random variable $X$ has hypergeometric distribution with parameters $N,n,m$ if $X:=|S\cap [m]|$, where $S$ is a random subset of $[N]$ of size $n$. 
We write $X\sim B(n,p)$ if $X$ has binomial distribution with parameters $n,p$. We will often use the following Chernoff-type bound.

\begin{lemma}[see {\cite[Remark 2.5 and Theorem 2.10]{JLR}}] \label{lem:chernoff}
Let $X\sim B(n,p)$ or let $X$ have a hypergeometric distribution with parameters $N,n,m$.
Then
\begin{align*}
\prob{|X - \expn(X)| \geq t} \leq 2\eul^{-2t^2/n}.
\end{align*}
\end{lemma}

\section{Vortices} \label{sec:vortex}

As explained earlier, our proof method involves an `iterative' absorption process, where in each iteration, we ensure that the leftover edges are all contained in some set $U_i$ where $U_i$ is much smaller than $U_{i-1}$. The underlying structure is that of a `vortex' which we introduce now.

\begin{defin}
Let $G$ be a graph on $n$ vertices and $W\In V(G)$. A $(\delta,\mu,m)$-vortex in $G$ surrounding $W$ is a sequence $U_0 \supseteq U_1 \supseteq \dots \supseteq U_\ell$ such that
\begin{enumerate}[label=(\roman*)]
\item[(V1)] $U_0=V(G)$;
\item[(V2)] $|U_i|=\lfloor \mu|U_{i-1}| \rfloor$ for all $i \in [\ell]$;
\item[(V3)] $|U_\ell|=m$ and $W\In U_\ell$;
\item[(V4)] $d_G(x,U_i) \ge \delta |U_i|$ for all $i \in [\ell]$, $x \in U_{i-1}$.
\end{enumerate}
\end{defin}

Often $W$ plays no role and we just refer to a $(\delta,\mu,m)$-vortex in $G$ in this case.

Our aim is to show that every large graph of high minimum degree contains a vortex such that the final set $U_\ell$ has constant size. We can also guarantee that a small given set $W$ is contained in this final set (i.e.~the vortex surrounds $W$). This will only be relevant in Section~\ref{sec:deltafa} (with $|W|=1$). In the main application, we will have $W=\emptyset$. The next proposition follows easily from Lemma~\ref{lem:chernoff}.

\begin{prop} \label{prop:random subset}
Let $\delta \in [0,1]$ and $1/n\ll \mu$. Suppose that $G$ is a graph on $n$ vertices with $\delta(G) \ge \delta n$ and $W\In V(G)$ with $|W|\le 1/\mu$. Then $V(G)$ contains a subset $U$ of size $\lfloor \mu n \rfloor$ such that $d_G(x,U) \ge (\delta - 2 n^{-1/3})|U|$ for every $x \in V(G)$ and $W\In U$.
\end{prop}\COMMENT{
\proof
Let $U'$ be a random subset of $V(G)$ of size $\lfloor \mu n \rfloor$. For every $x\in V(G)$, $$\prob{d_G(x,U')<\delta \lfloor \mu n \rfloor -\mu n^{2/3}} \le 2 \eul^{-2\mu^2 n^{4/3}/\lfloor \mu n\rfloor}.$$ Hence, for $n$ sufficiently large, we deduce that $d_G(x,U') \ge \delta |U'| -\mu n^{2/3}$ holds for all $x\in V(G)$. Let $U$ arise from $U'$ by adding all vertices from $W\sm U'$ and deleting $|W\sm U'|$ vertices from $U'\sm W$. Then, for every $x\in V(G)$, $d_G(x,U)\ge \delta |U| -\mu n^{2/3} -1/\mu \ge (\delta - 2 n^{-1/3})|U|$.
\endproof}

We now iterate the above proposition in order to obtain a vortex for a given graph $G$ such that the final set has constant size. 

\begin{lemma} \label{lem:get vortex}
Let $\delta\in[0,1]$ and $1/m'\ll \mu<1$. Suppose that $G$ is a graph on $n \ge m'$ vertices with $\delta(G) \ge \delta n$ and $W\In V(G)$ with $|W|\le 1/\mu$. Then $G$ has a $(\delta-\mu,\mu,m)$-vortex surrounding $W$ for some $\lfloor \mu m' \rfloor \le m \le m'$.
\end{lemma}

\proof
Recursively define  $n_0:=n$ and $n_i:=\lfloor \mu n_{i-1} \rfloor$. Observe that $\mu^i n \ge n_i \ge \mu^i n - 1/(1-\mu)$. Further, for $i\in \bN$, let $a_i:=n^{-1/3}\sum_{j\in [i]}\mu^{-(j-1)/3}$, with $a_0:=0$. Let $\ell:=1+\max\set{i\ge 0}{n_i\ge m'}$ and let $m:=n_\ell$. Note that $\lfloor \mu m' \rfloor \le m \le m'$. Now, suppose that for some $i\in[\ell]$, we have already found a $(\delta-3a_{i-1},\mu,n_{i-1})$-vortex $U_0,\dots,U_{i-1}$ in $G$ surrounding $W$. In particular, $\delta(G[U_{i-1}]) \ge (\delta-3a_{i-1})n_{i-1}$. By Proposition~\ref{prop:random subset}, there exists a subset $U_i$ of $U_{i-1}$ of size $n_i$ such that $d_G(x,U_i) \ge (\delta-3a_{i-1}-2n_{i-1}^{-1/3})n_i$ for all $x\in U_{i-1}$ and $W\In U_i$. Thus, $U_0,\dots,U_i$ is a $(\delta-3a_i,\mu,n_i)$-vortex in $G$ surrounding $W$. Finally, $U_0,\dots,U_\ell$ is a $(\delta-3a_\ell,\mu,m)$-vortex in $G$ surrounding $W$. Observe that $$a_\ell =   n^{-1/3}\frac{\mu^{-\ell/3}-1}{\mu^{-1/3}-1}  \le \frac{(\mu^{\ell-1}n)^{-1/3}}{1-\mu^{1/3}} \le \frac{m'^{-1/3}}{1-\mu^{1/3}}  \le\mu /3$$ since $\mu^{\ell-1}n\ge n_{\ell-1}\ge m'$, and so the lemma follows.
\endproof

\section{Near-optimal decomposition} \label{sec:near optimal}

The goal of this section is to prove the following lemma. Given a vortex, it finds a `near-optimal' decomposition in a graph of sufficiently large minimum degree. The proof proceeds using an `iterative' absorption approach.

\begin{lemma} \label{lem:near optimal}
Let $\delta:=\max\Set{\delta_F^{0+},\delta_F^{vx}}$. Assume that $1/m \ll \mu,1/|F|$. Let $G$ be an $F$-degree-divisible graph with $\delta(G) \ge (\delta + 3\mu)|G|$ and let $U_0 \supseteq U_1 \supseteq \dots \supseteq U_\ell$ be a $(\delta+4\mu,\mu,m)$-vortex in $G$. Then there exists a subgraph $H_\ell$ of $G[U_\ell]$ such that $G-H_\ell$ is $F$-decomposable.
\end{lemma}

The proof of Lemma~\ref{lem:near optimal} relies on a number of further tools. Before we start proving Lemma~\ref{lem:near optimal}, we show how it implies Theorem~\ref{thm:almost cover}.

\lateproof{Theorem~\ref{thm:almost cover}}
Let $F$ and $\mu$ be given and choose $m'$ sufficiently large. Let $\delta:=\max\Set{\delta_F^{0+},\delta_F^{vx}}$. Now, let $G$ be any $F$-degree-divisible graph with $\delta(G)\ge (\delta+5\mu)n$. If $n\le m'$, then $e(G)\le m'^2$. If $n\ge m'$, then by Lemma~\ref{lem:get vortex}, $G$ contains a $(\delta+4\mu,\mu,m)$-vortex $U_0 \supseteq U_1 \supseteq \dots \supseteq U_\ell$ for some $\lfloor \mu m' \rfloor \le m \le m'$. Then, by Lemma~\ref{lem:near optimal}, there exists an $F$-collection that covers all edges of $G$ except some edges of $G[U_\ell]$, that is, at most $m^2\le m'^2$ edges.
\endproof

\subsection{Bounded covering of edges around a vertex}

Many times in the proof of Lemma~\ref{lem:near optimal}, we will wish to find copies of $F$ which cover all the edges around some vertex, $x$ say. To do so, we will use the definition of $\delta_F^{vx}$, but we will often further wish to ensure that no vertex lies in many of these copies of $F$ (except, of course, $x$). This motivates the following definition and lemma.

Let $\delta^{vx,bd}_F$ be the smallest $\delta\ge 0$ such that for all $\mu > 0$ there exists an $n_0\in\bN$ such that whenever $G$ is a graph on $n\geq n_0$ vertices with $\delta(G)\geq(\delta+\mu)n$, and $x$ is a vertex of $G$ with $gcd(F)\mid d_G(x)$, then $G$ contains a collection $\cF$ of edge-disjoint copies of $F$ such that $\set{xy}{y\in N_G(x)} \In E( \cF)$ and $\Delta(\bigcup \cF -x ) \le n_0$.

Clearly, $\delta^{vx,bd}_F \ge \delta^{vx}_F$. Maybe surprisingly, the following is true.

\begin{lemma} \label{lem:surprise}
For all graphs $F$, we have $\delta^{vx,bd}_F=\delta^{vx}_F$.
\end{lemma}

To prove Lemma~\ref{lem:surprise}, we need the following definition. Given a graph $G$ and a vertex $x\in V(G)$ with $r\mid d_G(x)$, we call a partition $V_1,\dots,V_k$ of $V(G)\sm \Set{x}$ a \defn{$(\delta,m,r)$-splitting with respect to $x$}, if for all $i\in[k]$:
\begin{enumerate}[label=$\bullet$]
\item $r\mid d_G(x,V_i)$;
\item $\delta(G[V_i\cup\Set{x}])\ge \delta |V_i|$;
\item $m/3\le |V_i|\le 2m$.
\end{enumerate}

We shall consider $(\delta,m,r)$-splittings where $m$ is bounded and $k$ is comparatively large. Their existence will follow from the next proposition and the subsequent lemma, which are similar to Proposition~\ref{prop:random subset} and Lemma~\ref{lem:get vortex} in their interplay.

\begin{prop} \label{prop:random partition}
Let $\delta \in [0,1]$ and $1/n\ll 1/r$. Suppose that $G$ is a graph on $n+1$ vertices with $\delta(G) \ge \delta n$ and $x$ is a vertex of $V(G)$ such that $r\mid d_G(x)$. Then $V(G)\sm\Set{x}$ can be partitioned into two sets $V_1,V_2$ such that $r\mid d_G(x,V_i)$ and $d_G(v,V_i) \ge (\delta - 2n^{-1/3})|V_i|$ for every $v \in V(G)$ and $i\in\Set{1,2}$, and $|V_1|\le |V_2|\le |V_1|+2r$.
\end{prop}

\proof
Let $V_1',V_2'$ be a partition of $V(G)\sm\Set{x}$ such that
\begin{enumerate}
\item $|V_1'|\le |V_2'|\le |V_1'|+1$;
\item $d_G(v,V_i')\ge (\delta-n^{-1/3})|V_i'|$ for all $v\in V(G)$ and $i\in[2]$.
\end{enumerate}
That such a partition exists can be seen by choosing $V_1'$ as a random subset of $V(G)\sm\Set{x}$ of size $\lfloor n/2\rfloor$. Applying Lemma~\ref{lem:chernoff} shows that with probability at least $3/4$, $V_1'$ and $V_2':=V(G)\sm(\Set{x}\cup V_1')$ have the desired properties. Then by removing at most $r-1$ neighbours of $x$ from $V_1'$ and adding them to $V_2'$, we can obtain the desired partition $V_1,V_2$.
\endproof

\begin{lemma}  \label{lem:splitting}
Let $\delta \in [0,1]$ and $1/m\ll \mu,1/r$. Suppose that $G$ is a graph on $n+1 \ge m+4r$ vertices with $\delta(G) \ge \delta n$ and $x$ is a vertex of $G$ with $r\mid d_G(x)$. Then $G$ has a $(\delta-\mu,m,r)$-splitting with respect to $x$.
\end{lemma}

\proof
Define $n_j^+:=n2^{-j}+ 4r(1-2^{-j})$ and $n_j^-:=n2^{-j}- 4r(1-2^{-j})$, which need not be integers.
Further, for $j\in \bN$, let $a_j:=n^{-1/3}\sum_{j'\in [j]}2^{(j'-1)/3}$, with $a_0:=0$. Let $\ell:=1+\max\set{j\ge 0}{n_j^+\ge m+4r}$.

Now, suppose that for some $j\in[\ell]$, we have already found a partition $\cP_{j-1}$ of $V(G)\sm \Set{x}$ such that for all $V\in \cP_{j-1}$
\begin{enumerate}[label=(\roman*)$_{j-1}$]
\item $n^-_{j-1}\le |V| \le n^+_{j-1}$;
\item $r\mid d_G(x,V)$;
\item $\delta(G[V\cup\Set{x}])\ge (\delta -3a_{j-1})|V|$.
\end{enumerate}

(Note that we can take $\cP_0:=\Set{V(G)\sm\Set{x}}$.) We now find a refinement $\cP_{j}$ of $\cP_{j-1}$ such that (i)$_j$--(iii)$_j$ hold.

Consider $V\in \cP$ and let $G':=G[V\cup \Set{x}]$. Since $|V|\ge n^-_{\ell-1}\ge n^+_{\ell-1} -8r \ge m-4r$, we can apply Proposition~\ref{prop:random partition} to obtain a partition of $V$ into $V_1,V_2$ such that $r\mid d_{G'}(x,V_i)$ and $d_{G'}(v,V_i) \ge (\delta - 3a_{j-1}-2|V|^{-1/3})|V_i|$ for every $v \in V\cup\Set{x}$ and $i\in\Set{1,2}$, and $|V_1|\le |V_2|\le |V_1|+2r$. It is easy to check that (i)$_j$--(iii)$_j$ hold for $V_1$ and $V_2$.\COMMENT{Have to check $2|V|^{-1/3}\le 3(a_j-a_{j-1})=3n^{-1/3}2^{(j-1)/3}$, which amounts to $|V|\ge 8n/(27\cdot 2^{j-1})$. True since $|V|\ge n^-_{j-1}\ge n2^{-(j-1)}-4r$.}
Hence, refining every $V\in \cP_{j-1}$ in this way yields a partition $\cP_{j}$ such that (i)$_j$--(iii)$_j$ hold.

Observe that $$a_\ell=  n^{-1/3}\frac{2^{\ell/3}-1}{2^{1/3}-1}  \le \frac{(n2^{-(\ell-1)})^{-1/3}}{1-2^{-1/3}} \le \frac{m^{-1/3}}{1-2^{-1/3}}  \le\mu /3$$ since $n2^{-(\ell-1)}\ge n^+_{\ell-1}-4r\ge m$. Thus $\cP_{\ell}$ is the desired splitting as $m/3\le n^-_\ell \le n^+_\ell \le 2m$.
\endproof

\lateproof{Lemma~\ref{lem:surprise}}
Let $r:=gcd(F)$. It is sufficient to show $\delta^{vx,bd}_F \le \delta^{vx}_F$. Let $\mu>0$. Let $n'=n'(\mu/4,F)$ be such that whenever $G'$ is a graph on at least $n'$ vertices with $\delta(G')\ge (\delta^{vx}_F+\mu/4)|G'|$ and $x'$ is a vertex of $G'$ with $gcd(F)\mid d_{G'}(x')$, then $G'$ contains an $F$-collection such that every edge at $x'$ is covered.
Let $1/n_0\ll 1/n',1/r,\mu$.
Let $G$ be any graph on $n\ge n_0$ vertices with $\delta(G)\ge (\delta+\mu)n$ and let $x\in V(G)$ with $r\mid d_G(x)$. We have to find an $F$-collection $\cF$ such that $\cF$ covers all edges at $x$ and $\Delta(\bigcup\cF -x)\le n_0$.

By Lemma~\ref{lem:splitting}, there exists a $(\delta+\mu/2,n_0/2,r)$-splitting $V_1,\dots,V_k$ with respect to $x$. For each $i\in [k]$, let $G_i:=G[V_i\cup\Set{x}]$. Thus, $|G_i|\ge |V_i| \ge n_0/6 \ge n'$, $\delta(G_i)\ge (\delta+\mu/2)|V_i| \ge (\delta+\mu/4)|G_i|$ and $r\mid d_{G_i}(x)$. Hence, by our choice of $n'$, $G_i$ contains an $F$-collection $\cF_i$ such that every edge at $x$ is covered. Let $\cF:=\cF_1\cup\dots\cup\cF_k$. Then, $\Delta(\bigcup\cF -x)\le \max\Set{|V_1|,\dots,|V_k|} \le n_0$.
\endproof

\subsection{Bounded approximate decompositions}

For $\gamma \ge 0$, a \defn{$\gamma$-bounded approximate $F$-decomposition} of an $n$-vertex graph $G$ is a collection $\cF$ of edge-disjoint copies of $F$ contained in $G$ such that $\Delta(G- \bigcup \cF) \le \gamma n$. Let $\delta^{\gamma,bd}_F$ be the smallest $\delta\ge 0$ such that for all $\mu > 0$ there exists an $n_0\in\bN$ such that whenever $G$ is a graph on $n\geq n_0$ vertices with $\delta(G)\geq(\delta+\mu)n$, then $G$ has a $\gamma$-bounded approximate $F$-decomposition.

Trivially, for all $\gamma>0$ we have $\delta_F^{\gamma,bd}\ge \delta_F^{\gamma/2}$, the threshold for a $\gamma$-approximate decomposition.\COMMENT{If $\Delta(H)\le \gamma n$, then $e(H)\le \gamma n^2/2$}
Here, we will build to giving an upper bound for $\delta_F^{\gamma,bd}$ in Lemma~\ref{lem:bound-max-degree}. To find a bounded $\gamma$-approximate $F$-decomposition of a graph $G$ with large minimum degree, we will start by breaking $G$ into a large (but constant) number of edge-disjoint subgraphs which each have a high minimum degree but much fewer vertices than $G$. We then iteratively find approximate decompositions of these subgraphs. In doing so, we track vertices which have high minimum degree in the remainder of some previous approximate decomposition, and ensure these vertices always have small degree in the remainder of the later approximate decompositions. We will use the following lemma to find approximate decompositions where vertices in a specified subset $X$ have low degree in the remainder.

\begin{lemma} \label{lem:take-care-of-bad}
Let $1/n \ll \eta \ll \mu,1/|F|$. Let $\delta:=\max\Set{\delta^{\eta}_F,\delta^{vx}_F}$. Suppose that $G$ is a graph on $n$ vertices with $\delta(G) \ge (\delta+\mu)n$ and that $X$ is a subset of $V(G)$ of size at most $\eta^{1/3} n$. Then there exists a subgraph $H$ of $G$ such that $G-H$ is $F$-decomposable and $Y:=\set{x\in V(G)}{d_H(x)>\sqrt{\eta} n}$ has size at most $4\sqrt{\eta} n$ and does not contain any vertex from $X$.
\end{lemma}

\proof
By Lemma~\ref{lem:surprise}, we may assume that the following holds:
\begin{itemize}
\item[($\ast$)] Whenever $G'$ is a graph on at least $n/2$ vertices with $\delta(G') \ge (\delta+\mu/2)|G'|$ and $x\in V(G')$ with $gcd(F)\mid d_{G'}(x)$, then $G'$ contains a collection $\cF$ of edge-disjoint copies of $F$ such that all edges at $x$ are covered and $\Delta(\bigcup\cF - x) \le \eta^{-1/4}$.
\end{itemize}

Let $x_1,\dots,x_\ell$ be an enumeration of $X$. For $i\in[\ell]$, let $0\le r_i < gcd(F)$ be such that $d_G(x_i)\equiv r_i \mod {gcd(F)}$. Let $A_i$ be a set of $r_i$ vertices in $N_G(x_i)\sm X$.\COMMENT{We may assume that the $A_i$'s are pairwise disjoint.} Let $G_0$ be the graph obtained from $G$ by deleting all edges from $x_i$ to $A_i$, so $gcd(F)\mid d_{G_0}(x_i)$ for all $i\in[\ell]$ and $\delta(G_0)\ge (\delta+3\mu/4)n$. We will now successively find $F$-collections $\cF_i$ such that $\cF_i$ is a collection of edge-disjoint copies of $F$ in $G_i$ covering all edges at $x_i$ and $\Delta(\bigcup\cF_i - x_i) \le \eta^{-1/4}$, where $G_i:=(G_0-\bigcup_{j\in[i-1]}\bigcup\cF_{j})\sm\Set{x_1,\dots,x_{i-1}}$. Suppose that for some $i\in[\ell]$, we have already found $\cF_1,\dots,\cF_{i-1}$. Note that $gcd(F)\mid d_{G_i}(x_i)$.\COMMENT{It certainly holds for $G_0-\bigcup_{j\in[i-1]}\bigcup\cF_{j}$, but in this graph, no more edges between $x_i$ and $x_j$ with $j<i$ are left.} Moreover, $\delta(G_i)\ge (\delta +3\mu/4)n - \eta^{-1/4}(i-1) \ge (\delta + \mu/2)n$. Therefore, by $(\ast)$, there exists a collection $\cF_i$ of edge-disjoint copies of $F$ in $G_i$ such that all edges at $x_i$ are covered and $\Delta(\bigcup\cF_i - x_i) \le \eta^{-1/4}$.

Let $\cF':=\bigcup_{i\in[\ell]}\cF_i$ and $G_{\ell+1}:=(G-\bigcup\cF')\sm\Set{x_1,\dots,x_\ell}$. So $\delta(G_{\ell+1}) \ge (\delta+\mu/2)n$. Let $\cF''$ be an $\eta$-approximate $F$-decomposition of $G_{\ell+1}$. Let $H:=G-\bigcup\cF'-\bigcup\cF''$ and $Y:=\set{x\in V(G)}{d_H(x)>\sqrt{\eta} n}$. Since $d_H(x)< gcd(F)$ for all $x\in X$, we have $Y\cap X=\emptyset$. Finally, $e(H)\le \eta n^2+|X|gcd(F) \le 2\eta n^2$ and $2e(H) \ge |Y|\sqrt{\eta} n$. Hence, $|Y|\le 4\sqrt{\eta} n$.
\endproof

In order to obtain an upper bound for $\delta_F^{\gamma,bd}$, we need a $K_t$-decomposition of $K_s$ for some large $t$ and some even larger $s$. We could apply Wilson's theorem, but we don't need such heavy machinery here, only the following simple proposition. We include a proof for completeness.

\begin{prop} \label{prop:clique dec}
Let $p$ be a prime. Then for every $k\in \bN$, $K_{p^k}$ has a $K_p$-decomposition. 
\end{prop}

\proof
First, we prove that $K_p(p)$ is $K_p$-decomposable. Let $V_1,\dots,V_p$ be the partition of the vertex set into independent sets of size $p$ and let $v_{i1},\dots,v_{ip}$ be an enumeration of $V_i$. We define a set $\cF$ of $K_p$'s as follows. The $p$-tuple $(v_{1{i_1}},\dots,v_{p{i_p}})$ is the vertex set of a copy of $K_p$ in $\cF$ if and only if there exists an $r\in \Set{0,1,\dots,p-1}$ such that $i_{j+1}-i_j \equiv r \mod{p}$ for all $j\in[p]$, where $i_{p+1}:=i_1$. It is easy to see that $\cF$ is a $K_p$-decomposition, as $p$ is prime.\COMMENT{Note that such a tuple is unique once $r$ and some $i_j$ are given. Given any edge $v_{i_1j_1}v_{i_2j_2}$, there is a unique solution to $j_2-j_1\equiv (i_1-i_2)r \mod{p}$ since $i_1\neq i_2$ and $p$ is prime.}

We now prove the statement by induction on $k$. For $k=1$, there is nothing to show. For $k>1$, do the following. Partition the vertices of $K_{p^k}$ into $p^{k-1}$ clusters of size $p$. The edges inside each cluster form a copy of $K_p$, so we can remove them. Consider the complete reduced graph where the clusters are vertices. By induction, this reduced graph has a $K_p$-decomposition. Every copy of $K_p$ in this decomposition corresponds to a copy of $K_p(p)$ in the original graph, which is $K_p$-decomposable by the above. 
\endproof

\begin{lemma} \label{lem:bound-max-degree}
For every $\gamma>0$, $\delta_F^{\gamma,bd} \le \max\Set{\delta_F^{0+},\delta^{vx}_F}$.
\end{lemma}

\proof
Let $1/n \ll \eta \ll 1/s  \ll 1/t \ll \mu,\gamma,1/|F|$ and assume that $t$ is prime and $s$ is a power of $t$.
Thus, by Proposition~\ref{prop:clique dec}, $K_s$ has a $K_t$-decomposition. Let $\delta:=\max\Set{\delta_F^{\eta},\delta^{vx}_F}$ and suppose that $G$ is a graph on $n$ vertices with $\delta(G) \ge (\delta+\mu)n$. We have to show that $G$ has a $\gamma$-bounded-approximate $F$-decomposition.

Let $\cP=\Set{V_1,\dots,V_s}$ be a partition of $V(G)$ with the following properties:
\begin{enumerate}[label=(\roman*)]
\item $|V_i| = (1\pm \eta)n/s$;
\item $d_G(x,V_i) \ge (\delta+2\mu/3)|V_i|$ for all $x\in V(G)$, $i\in[s]$.
\end{enumerate}
To see that such a partition exists, independently for every vertex $x\in V(G)$, choose an index $i\in[s]$ uniformly at random and put $x$ into $V_i$. Apply Lemma~\ref{lem:chernoff} to see that such a random partition satisfies (i) and (ii) with probability at least $3/4$.

Note that $\Delta(G-G[\cP]) \le (1 + \eta)n/s \le \gamma n/2$. Thus it is enough to show that $G[\cP]$ has a $\gamma/2$-bounded-approximate $F$-decomposition.

Let $\Set{T_1,\dots,T_\ell}$ be a $K_t$-decomposition of $K_s$ (where we assume that $V(K_s)=[s]$). Clearly, $\ell \le s^2$. For $i\in [\ell]$, define $G_i:=\bigcup_{jk \in E(T_i)}G[V_j,V_k]$. So the $G_i$ form a decomposition of $G[\cP]$. Moreover, using (i) and (ii), we deduce $$\delta(G_i)\ge (\delta+ 2\mu/3)(1-\eta)(t-1)n/s \ge (\delta+\mu/2)(1+\eta)tn/s \ge (\delta+\mu/2)|G_i|.$$

Start with $BAD:=\emptyset$. For $i=1,\dots,\ell$, do the following: Apply Lemma~\ref{lem:take-care-of-bad} with $G_i$, $\mu/2$ and $BAD\cap V(G_i)$ playing the roles of $G$, $\mu$ and $X$ to obtain a subgraph $H_i$ of $G_i$ such that $G_i-H_i$ is $F$-decomposable, $d_{H_i}(x) \le \sqrt{\eta} |G_i|$ for all $x\in BAD$ and $d_{H_i}(x)> \sqrt{\eta} |G_i|$ for at most $4\sqrt{\eta} |G_i|$ vertices $x\in V(G_i)$. Add all the vertices $x$ with $d_{H_i}(x)> \sqrt{\eta} |G_i|$ to $BAD$. Since $|BAD| \le s^2 4\sqrt{\eta}(1+\eta)tn/s \le \eta^{1/3} (1-\eta)tn/s$ at any time, the conditions of Lemma~\ref{lem:take-care-of-bad} are satisfied each time. Let $H:=\bigcup_{i\in [\ell]}H_i$ and let $x$ be any vertex of $G$. Crucially, $d_{H_i}(x)>\sqrt{\eta}(1+\eta)tn/s$ for at most one $i\in[\ell]$. Therefore, $$d_H(x) \le \ell \sqrt{\eta}(1+\eta)tn/s + (1+\eta)tn/s \le 2s\sqrt{\eta}tn + 2tn/s \le \gamma n/2,$$ as required.
\endproof

\subsection{Covering a pseudorandom remainder}

In proving the main lemma in this section, we will have the following situation. Given a small set $U$ in our graph $G$, we will have found copies of $F$ which cover all the edges in $G-U$ and most of the edges between $V(G)\setminus U$ and $U$. We will wish to find copies of $F$ which cover the remaining edges between $V(G)\setminus U$ and $U$ (while necessarily using some edges in $G[U]$). The following lemma tells us this is possible if our remaining edges satisfy certain pseudorandom conditions.

\begin{lemma} \label{lem:pseudorandom remainder}
Let $1/n \ll \rho \ll \mu,1/|F|$. Let $G$ be a graph on $n$ vertices and let $U$ be a subset of $V(G)$ of size at least $\mu n$. Let $W:=V(G)\sm U$ and let $w_1,\dots,w_p$ be an enumeration of $W$. Suppose there are sets $U_1,\dots,U_p \In U$ with the following properties:
\begin{enumerate}[label=(\roman*)]
\item $gcd(F)\mid d_G(w_i)$ for all $i \in [p]$; \label{eqn:coverUW:divisible}
\item $N_G(w_i)\In U_i$ for all $i\in[p]$; \label{eqn:coverUW:neighbours}
\item $d_G(x,U_i) \ge (\delta_F^{vx}+\mu)|U_i|$ for all $x\in U_i\cup\Set{w_i}$; \label{eqn:coverUW:degrees}
\item $|U_i| \ge \rho |U|/2$; \label{eqn:coverUW:large enough}
\item $|U_i \cap U_j| \le 2\rho^2 |U|$ for all $1\le i < j\le p$; \label{eqn:coverUW:intersection}
\item every $x\in U$ is contained in at most $2\rho n$ $U_i$'s. \label{eqn:coverUW:spread}
\end{enumerate}
Then there exists a subgraph $G_U$ of $G[U]$ such that $G_U\cup G[U,W]$ is $F$-decomposable and $\Delta(G_U)\le \mu^2 |U|$.
\end{lemma}

The proof of Lemma~\ref{lem:pseudorandom remainder} is quite similar to that of Lemma~10.7 in \cite{BKLO}, we include it for completeness. The proof will make use of the following result.

\begin{prop}[Jain, see {\cite[Lemma 8]{R}}] \label{prop:Jain}
Let $X_1, \ldots, X_n$ be Bernoulli random variables such that, for any $i \in [n]$ and any $x_1, \ldots, x_{i-1}\in \{0,1\}$,
 \begin{align*}
\prob{X_i = 1 \mid X_1 = x_1, \ldots, X_{i-1} = x_{i-1}} \leq p.
\end{align*}
Let $B \sim B(n,p)$ and $X:=X_1+\dots+X_n$.  Then $\prob{X \geq a} \leq \prob{B \geq a}$ for any~$a\ge 0$.
\end{prop}

\lateproof{Lemma~\ref{lem:pseudorandom remainder}}
Let $\Delta:=\rho^{-1/4}$. By Lemma~\ref{lem:surprise}, we may assume that the following holds.

\begin{itemize}
\item[($\ast$)] Whenever $G'$ is a graph on at least $\rho |U|/2$ vertices with $\delta(G') \ge (\delta_F^{vx}+\mu/2)|G'|$ and $x$ is a vertex in $G'$ with $gcd(F)\mid d_{G'}(x)$, then $G'$ contains a collection $\cF$ of edge-disjoint copies of $F$ such that all edges at $x$ are covered and $\Delta(\bigcup\cF - x) \le \Delta$. (In other words, there is a spanning subgraph $A$ of $G'-x$ such that $A \cup G'[V(A),\Set{x}]$ is $F$-decomposable and $\Delta(A) \le \Delta$.)
\end{itemize}

We want to find edge-disjoint subgraphs $T_1,\dots,T_p$ in $G[U]$ such that $V(T_i)=U_i$, $T_i\cup G[U_i,\Set{w_i}]$ is $F$-decomposable and $\Delta(T_i) \le \Delta$. Then, $G_U:=T_1\cup\dots\cup T_p$ is the desired subgraph, since $G_U \cup G[U,W] = \bigcup_{i\in [p]}(T_i\cup G[U_i,\Set{w_i}])$ by \ref{eqn:coverUW:neighbours} and $\Delta(G_U) \le \Delta \cdot 2\rho n \le \mu^2 |U|$, using \ref{eqn:coverUW:spread}.

We find $T_1, \dots, T_p$ in turn using a randomised algorithm. Let $t: = \lceil 8 \rho^{3/2} |U| \rceil $ and define $G_j:=G[U_j]$ for all $j\in [p]$.
Suppose that we have already found $T_1, \dots, T_{s-1}$ for some $s \in [p]$. We now define $T_s$ as follows.
Let $H_{s-1} : = \bigcup_{i=1}^{s-1}  T_{i}$ and let $G_s' : = (G - H_{s-1})[U_s]$.
If $ \Delta ( H_{s-1}[U_s] ) >  \Delta \rho^{3/2} n$, then let $A_1, \ldots, A_t$ be empty graphs on $U_s$.
If $\Delta ( H_{s-1}[U_s] ) \le  \Delta \rho^{3/2} n$, then 
\begin{align*}
\delta(G'_s) &\ge \delta(G[U_s])-\Delta(H_{s-1}[U_s]) \geq (\delta_F^{vx}+\mu)|U_s| -  \Delta \rho^{3/2} n \\
             &\ge (\delta_F^{vx}+\mu/2)(|U_s|+1) + (t-1)\Delta,
\end{align*}
by \ref{eqn:coverUW:degrees} and \ref{eqn:coverUW:large enough}. Thus, by ($\ast$), we can find $t$ edge-disjoint subgraphs $A_1,\dots,A_t$ of $G_s'$ which are all suitable candidates for~$T_s$.\COMMENT{Indeed, suppose we have found $A_1,\dots,A_{j-1}$ for some $j\in[t]$, then let $\tilde{G}:=(G'_s-(A_1\cup\dots\cup A_{j-1}))\cup G[U_s,w_s]$. By the above and \ref{eqn:coverUW:degrees}, $\delta(\tilde{G}) \ge (\delta_F^{vx}+\mu/2)|\tilde{G}|$. Moreover, $d_{\tilde{G}}(w_s)=d_G(w_s)\equiv 0 \mod{gcd(F)}$. 
So there exists a spanning subgraph $A_j$ of $G'_s-(A_1\cup\dots\cup A_{j-1})$ such that $A_j\cup G[U_s,w_s]$ is $F$-decomposable and $\Delta(A_j)\le \Delta$.}

In either case, we have found edge-disjoint subgraphs $A_1, \dots, A_t$ of $G_s'$.
Pick $i \in [t]$ uniformly at random and set $T_s : = A_i$.
The lemma follows if the following holds with positive probability:
\begin{align}
\Delta ( H_{p}[U_j] ) \le   \Delta \rho^{3/2} n \text{ for all $j\in[p]$.}	\label{eqn:Krmkey}
\end{align}
To analyse this, for $s,j \in [p]$ and $u \in U_j$, let $Y^{j,u}_s$ be the indicator function of the event $\Set{d_{T_s}(u,U_j)\ge 1}$.
Let $X^{j,u} :  = \sum_{s=1}^p Y^{j,u}_s$.
Note that $d_{T_s}(u,U_j)\le Y^{j,u}_s\Delta$, so $d_{H_p}(u,U_j) \leq \Delta X^{j,u}$.
Therefore to prove \eqref{eqn:Krmkey} it suffices to show that with positive probability, $X^{j,u} \le  \rho^{3/2} n$ for all $j \in [p]$ and $u \in U_j$. 

Fix $j \in [p]$ and $u \in U_j$.
Let $S^{j,u}$ be the set of indices $s \ne j$ such that $u \in U_{s}$.
By~\ref{eqn:coverUW:spread}, $ |S^{j,u}|  \le 2 \rho n$. 
Note that $Y_s^{j,u} = 0 $ for all $s \notin S^{j,u} \cup \{j\}$.
So 
\begin{align}
X^{j,u} \le 1 + \sum_{s \in S^{j,u}} Y^{j,u}_{s}. \label{eqn:Xju}
\end{align}
Let $s_1, \dots, s_{|S^{j,u}|}$ be the enumeration of $S^{j,u} $ such that $s_b < s_{b+1}$ for all $b \in [|S^{j,u}|-1]$.
Note that  $d_{G_{s_b}} ( u, U_j )  \le |U_j \cap U_{s_b}| \le 2 \rho^2 |U|$ by~\ref{eqn:coverUW:intersection}.
So at most $2 \rho^2 |U|$ of the subgraphs $A_{i}$ that we picked in $G'_{s_b}$ contain an edge incident to~$u$ in~$G_j$. 
This implies that for all $y_1,\dots,y_{b-1}\in \{0,1\}$ and all $b \in [|S^{j,u}|] $,
\begin{align*}
	\mathbb{P} ( Y^{j,u}_{s_b} = 1 \mid Y^{j,u}_{s_1} = y_1, \dots, Y^{j,u}_{s_{b-1}} = y_{b-1})  \le 
	\frac{2 \rho^2 |U|}{ t } \le \frac {\rho^{1/2}} {4}.
\end{align*}
Let $B \sim B( |S^{j,u}|  , \rho^{1/2}/4 )$.
By \eqref{eqn:Xju}, Proposition~\ref{prop:Jain} and the fact that $|S^{j,u}|\le 2\rho n$ we have that
\begin{align*}
\mathbb{P}( X^{j,u} >  \rho^{3/2} n ) 
& \le
\mathbb{P}( \sum_{s \in S^{j,u}} Y^{j,u}_{s}  >  3 \rho^{3/2} n/4 ) 
\le \mathbb{P}( B >  3 \rho^{3/2} n/4  )
\\
& \le \mathbb{P}( | B - \mathbb{E}(B) | >  \rho^{3/2} n/4 )
\le 2 e^{-  \rho^{2} n/16},
\end{align*} 
where the last inequality holds by Lemma~\ref{lem:chernoff}. Since there are at most $n^2$ pairs $(j, u)$, there is a choice of $T_1, \ldots, T_p$ such that $X^{j,u} \le  \rho^{3/2} n$ for all $j \in [p]$ and all $u \in U_j$, proving the claim. 
\endproof

\subsection{Proof of Lemma~\ref{lem:near optimal}}

Before proving the main tool of this section (from which Lemma~\ref{lem:near optimal} will follow simply by induction), we need one final proposition. Given a subset $R$ of a graph $G$ with certain properties, and any sparse subgraph $H$ of $G-R$, it allows copies of $F$ to be found in $G$ which cover the edges of $H$ without covering any other edges in $G-R$.

\begin{prop} \label{prop:cover-sparse-graph}
Let $1/n \ll \gamma \ll \mu,1/|F|$. Let $G$ be a graph on $n$ vertices and let $V(G)=L\cupdot R$ such that $|R| \ge \mu n$ and $d_G(x,R)\ge (\delta_F^e+\mu)|R|$ for all $x\in V(G)$. Let $H$ be any subgraph of $G[L]$ such that $\Delta(H) \le \gamma n$. Then there exists a subgraph $A$ of $G$ such that $A[L]$ is empty, $A\cup H$ is $F$-decomposable and $\Delta(A)\le \mu^2 |R|$.
\end{prop}

\proof
Let $e_1,\dots,e_m$ be an enumeration of $E(H)$. We will find edge-disjoint copies $F_1,\dots,F_m$ in $G$ such that $F_i$ contains $e_i$ and $V(F_i)\cap L =V(e_i)$. Suppose we have already found $F_1,\dots,F_{j-1}$ for some $j\in [m]$. Let $G_{j-1}:=F_1\cup \dots \cup F_{j-1}$ and suppose that $\Delta(G_{j-1}) \le \sqrt{\gamma}n+|F|$. Let $BAD:=\set{x\in V(G)}{d_{G_{j-1}}(x)>\sqrt{\gamma} n}$. Note that for all $x\in L$, $d_{G_{j-1}}(x)\le |F|\Delta(H) \le \sqrt{\gamma} n$,\COMMENT{Since $V(F_i)\cap L =V(e_i)$, $x$ is only contained in $F_i$'s when $e_i$ is incident to $x$} so $BAD\cap L=\emptyset$. We have $$e(G_{j-1}) \le e(F)e(H) \le |F|^2\Delta(H)n \le |F|^2 \gamma n^2.$$ On the other hand, $2e(G_{j-1}) \ge |BAD|\sqrt{\gamma} n$. Thus, $|BAD| \le 2|F|^2\sqrt{\gamma}n \le \mu|R|/2$. Let $G':=(G-G_{j-1})[(R\sm BAD) \cup V(e_j)]$. Observe that $\delta(G') \ge (\delta_F^e +\mu/4)|G'|$, so there exists a copy $F_j$ of $F$ in $G'$ that contains $e_j$. Moreover, since $F_j$ does not contain any vertex of $BAD$,\COMMENT{Because $BAD\cap L=\emptyset$} we have ensured $\Delta(G_j)\le \sqrt{\gamma}n+|F|$ for the next step. Finally, $A:=\bigcup_{i\in [m]}(F_i-e_i)$ is the desired subgraph.
\endproof

We are now ready to prove the main tool that will enable us to prove Lemma~\ref{lem:near optimal} by induction.

\begin{lemma} \label{lem:cover-down}
Let $\delta:=\max\Set{\delta_F^{0+},\delta_F^{vx}}$ and $1/n \ll \mu,1/|F|$. Let $G$ be a graph on $n$ vertices and $U\In V(G)$ with $|U|=\lfloor \mu n \rfloor$.
Suppose that $\delta(G) \ge (\delta+2\mu)n$ and for all $x\in V(G)$, $d_G(x,U)\ge (\delta + \mu)|U|$. Then, if $gcd(F)\mid d_G(x)$ for all $x\in V(G)\sm U$, there exists an $F$-collection $\cF$ in $G$ such that every edge in $G-G[U]$ is covered, and $\Delta(\bigcup\cF[U]) \le \mu^2 |U|/4$.
\end{lemma}

Our strategy is as follows. Since $U$ is relatively small, we know that $G-G[U]$ still has high minimum degree. Therefore, we can obtain an approximate decomposition that uses no edges inside $U$, but covers almost all edges outside $U$. Before doing this, we set aside two sparse subgraphs $R'$ and $R''$ of $G[U,V(G)\sm U]$ with pseudorandom properties. Letting $H$ be the leftover of the approximate decomposition, we use $R'$ and some edges of $G[U]$ to cover all edges in $H[V(G)\sm U]$ using Proposition~\ref{prop:cover-sparse-graph}. Finally, we combine $H[U,V(G)\sm U]$ and the leftover of $R'$ with $R''$. Since $R''$ is relatively dense (compared to $R'$ and $H$) and has pseudorandom properties, we can cover all these edges using Lemma~\ref{lem:pseudorandom remainder}.

\lateproof{Lemma~\ref{lem:cover-down}} Choose new constants $\gamma,\xi,\rho>0$ such that $1/n\ll \gamma \ll \xi \ll \rho \ll \mu,1/|F|$. Let $W:=V(G)\sm U$. 
We will first choose suitable graphs $R'$ and $R''$ which we will put aside for later use. Let $k:=\lceil \xi^{-1} \rceil$ and $K:=\binom{k+1}{2}$. 

Let $V_1,\dots,V_K$ be a partition of $U$ with the following properties:
\begin{align}
d_G(x,V_i) &\ge (\delta+\mu/2)|V_i| \quad \text{for all $x\in V(G)$ and $i\in[K]$};  \label{eqn:prereserve:degrees} \\
|U|/2K &\le |V_i| \le 2|U|/K. \label{eqn:prereserve:size}
\end{align}
To see that such a partition exists, independently for every vertex $u\in U$, choose an index $i\in[K]$ uniformly at random and put $u$ into $V_i$. Apply Lemma~\ref{lem:chernoff} to see that such a random partition has the desired properties with probability at least $3/4$.

Split $W$ arbitrarily into $k$ sets $W_1,\dots,W_k$ as evenly as possible and let $G_W^1,\dots,G_W^K$ be an enumeration of the $K$ graphs of the form $G[W_i]$ or $G[W_i,W_j]$. Thus, $G[W] = \bigcup_{i\in[K]}G_W^i$ and $|G_W^i|\le 2(|W|/k+1) \le 2\xi n$ for all $i\in[K]$.

For every $i\in[K]$, let $R_i:=G[V_i,V(G_W^i)]$. Let $R':=R_1\cup\dots\cup R_K$. Note that $d_{R'}(u) \le |V(G_W^i)| \le 2\xi n$ for all $u\in U$ and $d_{R'}(w) \le k \cdot 2 |U|/K \le 4\xi n$\COMMENT{$K\ge k^2/2$} for all $w\in W$, so 
\begin{align}
\Delta(R') \le 4\xi n.
\end{align}

Let $G':=G-R'$. So $d_{G'}(x,U) \ge (\delta+3\mu/4)|U|$ for all $x\in V(G)$.
Let $p:=|W|$ and let $U_1',\dots,U_p'$ be subsets of $U$ with the following properties:
\begin{enumerate}[label=(\alph*)]
\item $d_{G'}(x,U_i') \ge (\delta+\mu/2)|U_i'|$ for all $x\in V(G)$ and $i\in[p]$;   \label{eqn:reserve1:degrees} 
\item $\rho |U|/2 \le |U_i'| \le 2\rho |U|$ for all $i\in[p]$;                              \label{eqn:reserve1:size} 
\item $|U_i' \cap U_j'| \le 3\rho^2 |U|/2$ for all $1\le i < j\le p$;  \label{eqn:reserve1:intersection} 
\item each $u\in U$ is contained in at most $3\rho p/2$ of the $U_i'$.           \label{eqn:reserve1:spread}
\end{enumerate}
That these subsets exist can again be seen by a probabilistic argument. Indeed, for every pair $(u,i) \in U\times [p]$, include $u$ in $U_i'$ with probability $\rho$ independently of all other pairs. Applying Lemma~\ref{lem:chernoff} shows that the random sets $U_1',\dots,U_p'$ satisfy the desired properties with probability at least $3/4$.

Let $w_1,\dots,w_p$ be an enumeration of $W$ and let $R'':=\bigcup_{i\in [p]}G'[U_i',\Set{w_i}]$. By \ref{eqn:reserve1:size} and \ref{eqn:reserve1:spread}, $\Delta(R'') \le \max\Set{2\rho |U|,3\rho p/2}\le 2\rho n$.

Let $G'':=G-G[U]-R'-R''$.
Observe that $\delta(G'')\ge (\delta+\mu/2)n$. We now apply Lemma~\ref{lem:bound-max-degree} to find an approximate decomposition of $G''$. More precisely, by Lemma~\ref{lem:bound-max-degree}, there exists a subgraph $H$ of $G''$ such that $G''-H$ has an $F$-decomposition $\cF_1$ and $\Delta(H)\le \gamma n$. Let $H_W:=H[W]$ and $H_{UW}:=H[U,W]$.

Next, we want to cover the edges of $H_W$ using $R'$. Recall that $G_W^1,\dots,G_W^K$ is a decomposition of $G[W]$. For all $i\in [K]$, let $H_i:=H_W\cap G_W^i$ and $G_i:=G[V_i]\cup R_i \cup H_i$. So the $H_i$ decompose $H_W$. Note that $V(G_i)=V_i\cup V(G_W^i)$ and thus $\mu\xi^2 n/10 \le |V_i| \le |G_i| \le 3\xi n$, implying that $|V_i|\ge \xi^2 |G_i|$.\COMMENT{$|G_i|\le 2\mu n/K +2\xi n\le 3\xi n$; $|G_i|\ge |V_i|\ge |U|/2K \ge \mu\xi^2 n/10\ge \xi^2|G_i|$} Moreover, $d_{G_i}(x,V_i) \ge (\delta+\xi^2)|V_i|$ for all $x \in V(G_i)$ by \eqref{eqn:prereserve:degrees} and our choice of $R_i$. Since $\delta\ge \delta_F^{vx}\ge \delta_F^e$ and $\Delta(H_i) \le \gamma n \le \sqrt{\gamma}|G_i|$, we can apply Proposition~\ref{prop:cover-sparse-graph} with $\xi^2$ and $\sqrt{\gamma}$ playing the roles of $\mu$ and $\gamma$ to obtain a subgraph $A_i$ of $G_i$ such that $A_i\cup H_i$ is $F$-decomposable, $A_i[V(G_i)\sm V_i]$ is empty and $\Delta(A_i)\le \xi^4|V_i|$.
Let $A:=A_1\cup\dots\cup A_K$. So $A\cup H_W$ has an $F$-decomposition $\cF_2$ and $\Delta(A) \le \xi n$.

We now want to cover the remaining edges of $H_{UW}\cup R'$ using $R''$. 
Let $G''':=G-\bigcup\cF_1-\bigcup \cF_2$. Note that $G'''[W]$ is empty. For every $i\in [p]$, let $U_i'':=N_{G'''}(w_i)\sm U_i'$. Hence, $|U_i''|\le \Delta(H_{UW}) + \Delta(R') \le (\gamma +4\xi) n$. Let $U_i:=U_i'\cup U_i''$. We want to check the conditions of Lemma~\ref{lem:pseudorandom remainder} for $G'''$ and $U_1,\dots,U_p$, with $\mu/4$ playing the role of $\mu$. Conditions~\ref{eqn:coverUW:divisible}, \ref{eqn:coverUW:neighbours} and \ref{eqn:coverUW:large enough} clearly hold. To see that \ref{eqn:coverUW:degrees} holds, let $i\in [p]$ be arbitrary and consider first $u\in U_i$. Since $G'''[U]=G-\bigcup \cF_2$, we have
\begin{align*}
d_{G'''}(u,U_i) &\ge d_{G}(u,U_i') - \Delta(A) \overset{\ref{eqn:reserve1:degrees}}{\ge} (\delta + \mu/2)|U_i'| -\xi n \ge (\delta + \mu/4)|U_i|.
\end{align*}
Secondly, $d_{G'''}(w_i,U_i) \ge d_{R''}(w_i,U_i) = d_{G'}(w_i,U_i') \overset{\ref{eqn:reserve1:degrees}}{\ge} (\delta+\mu/2)|U_i'| \ge (\delta+\mu/4)|U_i|$.

For \ref{eqn:coverUW:intersection}, observe that $$|U_i\cap U_j| \overset{\ref{eqn:reserve1:intersection}}{\le} 3\rho^2 |U|/2 + 2(\gamma+4\xi)n \le 2\rho^2|U|$$ for all $1\le i <j \le p$.\COMMENT{$|U_i\cap U_j|=|(U_i'\cap U_j')| +| U_i''|+ | U_j''|$}

Finally, note that every $u\in U$ is contained in at most $\Delta(H_{UW}) + \Delta(R') \le (\gamma +4\xi) n$ of the sets $U_i''$. Combining this with \ref{eqn:reserve1:spread}, we conclude that $u$ is contained in at most $3\rho p/2 + (\gamma +4\xi) n \le 2\rho n$ of the $U_i$'s, establishing \ref{eqn:coverUW:spread}. Thus, by applying Lemma~\ref{lem:pseudorandom remainder}, we obtain a subgraph $G_U$ of $G'''[U]$ such that $G_U\cup G'''[U,W]$ has an $F$-decomposition $\cF_3$ and $\Delta(G_U) \le (\mu/4)^2 |U|$.

Let $\cF:=\cF_1\cup\cF_2\cup\cF_3$. By construction, $\cF$ covers every edge of $G-G[U]$, and $\Delta(\cF[U])\le \Delta(A)+\Delta(G_U) \le \mu^2 |U|/4$.
\endproof

We can finally deduce Lemma~\ref{lem:near optimal} by inductively applying Lemma~\ref{lem:cover-down}.

\lateproof{Lemma~\ref{lem:near optimal}}
If $\ell = 0$, then we can put $H_\ell := G$. Therefore, let us assume $\ell \ge 1$. We prove the following stronger statement by induction on~$\ell$.

\medskip \emph{Let $G$ be an $F$-degree-divisible graph with $\delta(G) \ge (\delta + 3\mu)|G|$ and let $U_1$ be a subset of $V(G)$ of size $\lfloor \mu |G| \rfloor$ such that $d_G(x,U_1) \ge (\delta + 7\mu/2)|U_1|$ for all $x \in V(G)$. Let $U_1 \supseteq U_2 \supseteq \dots \supseteq U_\ell$ be a $(\delta+4\mu,\mu,m)$-vortex in $G[U_1]$. Then there exists a subgraph $H_\ell$ of $G[U_\ell]$ such that $G-H_\ell$ is $F$-decomposable.}
\medskip

If $\ell=1$, then Lemma~\ref{lem:cover-down} applied to $G$ and $U_1$ yields a subgraph $H_1$ of $G[U_1]$ such that $G-H_1$ is $F$-decomposable. So let us assume that $\ell \ge 2$ and that the claim holds for $\ell-1$. Let $G':=G-G[U_2]$. Note that $\delta(G') \ge (\delta + 2\mu)|G'|$ and $d_{G'}(x,U_1)\ge (\delta+\mu)|U_1|$. Furthermore, $d_{G'}(x)=d_G(x)$ and thus $gcd(F)\mid d_{G'}(x)$ for all $x\in V(G') \sm U_1$. By Lemma~\ref{lem:cover-down}, there exists an $F$-collection $\cF$ in $G'$ that covers all edges outside $G'[U_1]$ and satisfies $\Delta(\bigcup\cF[U_1])\le \mu^2 |U_1|/4$. Let $G'':=G[U_1]-\bigcup\cF$. So $G''$ is $F$-degree-divisible and $\delta(G'') \ge (\delta + 3\mu)|G''|$.\COMMENT{$d_{G''}(x)\ge d_G(x,U_1)-\mu^2 |U_1|/4 \ge (\delta+7\mu/2)|U_1|-\mu^2 |U_1|/4 \ge (\delta+3\mu) |U_1|$} Moreover, $U_2$ is a subset of $V(G'')$ of size $\lfloor \mu |G''| \rfloor$ such that $d_{G''}(x,U_2) \ge (\delta + 4\mu)|U_2| - \Delta(\bigcup\cF[U_1]) \ge (\delta + 7 \mu/2) |U_2|$ for all $x \in V(G'')$. Finally, $U_2\supseteq \dots \supseteq U_\ell$ is a $(\delta+4\mu,\mu,m)$-vortex in $G''[U_2]$, since $G''[U_2]=G[U_2]$. By induction, there exists a subgraph $H_\ell$ of $G''[U_\ell]$ such that $G''-H_\ell$ has an $F$-decomposition $\cF'$. But now, $\cF \cup \cF'$ is an $F$-decomposition of $G-H_\ell$, completing the proof.
\endproof

\section{Regularity} \label{sec:regularity}

In this section, we state Szemer\'edi's regularity lemma and related tools. Let $G$ be a graph with two disjoint sets of vertices $A,B \In V(G)$. The \defn{density} of $G[A,B]$ is then defined as $\alpha_G(A,B):=e_G(A,B)/(|A||B|)$. Given $\eps>0$, we call $G[A,B]$ \defn{$\eps$-regular} if for all sets $X \In A$ and $Y \In B$ with $|X| \ge \eps|A|$ and $|Y| \ge \eps|B|$, we have $|\alpha_G(X,Y)-\alpha_G(A,B)|<\eps$.

\begin{fact} \label{fact:regularity}
Let $G[A,B]$ be $\eps$-regular with density $\alpha$ and let $Y\In B$ with $|Y|\ge \eps|B|$. Then all but at most $\eps|A|$ vertices of $A$ have at least $(\alpha-\eps)|Y|$ neighbours in~$Y$.  
\end{fact}

\begin{fact} \label{fact:slicing}
Let $G[A,B]$ be $\eps$-regular with density $\alpha$ and for some $c>\eps$, let $A'\In A$ and $B'\In B$ with $|A'|\ge c|A|$ and $|B'|\ge c|B|$. Then $G[A',B']$ is $2\eps/c$-regular with density $\alpha\pm \eps$.  
\end{fact}

\begin{lemma}[Regularity lemma] \label{lem:regularity}
For all $\eps>0$ and $k_0\in \bN$, there exists a $k_0'=k_0'(\eps,k_0)$ such that for all $\alpha \in [0,1]$ the following holds. Let $G$ be a graph on $n\ge k_0'$ vertices and $W_1,\dots,W_\ell$ a partition of $V(G)$ with $\ell\le k_0$. Then there exist a partition $V_0,V_1,\dots,V_k$ of $V(G)$ and a spanning subgraph $G'$ of $G$ satisfying the following:
\begin{enumerate}[label=(R\arabic*)]
\item $k_0 \le k \le k_0'$;
\item $|V_0| \le \eps n$;
\item $|V_1|=\dots=|V_k|$;
\item $d_{G'}(x) \ge d_G(x)-(\alpha+\eps)n$ for all $x \in V(G)$;
\item $G'[V_i]$ is empty for all $i\in [k]$;
\item for all $1\le i < j \le k$, $G'[V_i,V_j]$ is either $\eps$-regular with density at least $\alpha$ or empty;
\item for all $i\in [k]$ and $j\in[\ell]$, $V_i\cap W_j=\emptyset$ or $V_i\In W_j$.
\end{enumerate}
\end{lemma}

Given a graph $G$ and a partition $V_1,\dots,V_k$ of $V(G)$, we associate a so-called \defn{reduced graph $R$} with this partition, where $R$ has vertex set $[k]$ and $ij\in E(R)$ if and only if $G[V_i,V_j]$ is non-empty. We refer to the function $\sigma\colon V(G) \rightarrow V(R)$ such that $x\in V_{\sigma(x)}$ for all $x\in V(G)$ as the \defn{cluster function}. In the setting of Lemma~\ref{lem:regularity}, we slightly abuse notation and say that $R$ is the reduced graph of $V_1,\dots,V_k$ if $R$ is the reduced graph of $V_1,\dots,V_k$ with respect to $G'[V(G)\sm V_0]$.

\begin{prop} \label{prop:reduced graph}
Suppose that $G$ is a graph on $n$ vertices with $\delta(G) \ge \delta n$ and that $G'$ and $V_0,\dots,V_k$ satisfy (R1)--(R6). Then $\delta(R)\ge (\delta -\alpha -2\eps)k$, where $R$ is the reduced graph of $V_1,\dots,V_k$.
\end{prop}\COMMENT{
\proof
Let $i\in [k]$ and $L:=|V_i|$. Then, $$\delta(G')L \le e_{G'}(V_i,V(G)\sm V_i) \le L(|V_0|+d_R(i)L).$$ But $\delta(G')\ge (\delta - \alpha -\eps)n$. Hence, $(\delta - \alpha -\eps)n \le \eps n+d_R(i)L$. Since $n/L\ge k$, it follows that $d_R(i)\ge (\delta - \alpha -2\eps)k$.
\endproof
}

We will often use the following embedding lemma which is straightforward to prove and has become known as the `key lemma'. 

\begin{lemma}[Key lemma] \label{lem:key lemma}
Let $1/n \ll \eps \ll \alpha \ll 1/m$. Let $G$ be a graph such that 
\begin{enumerate}[label=$\bullet$]
\item $V(G)=V_1\cupdot \dots \cupdot V_k$;
\item for all $i\in[k]$, $|V_i|=n$ and $G[V_i]$ is empty;
\item for all $1\le i< j\le k$, $G[V_i,V_j]$ is either $\eps$-regular with density at least $\alpha$ or empty.
\end{enumerate}
Let $R$ be the reduced graph of $V_1,\dots,V_k$.
Let $H$ be a graph of order at most $m$ and suppose that $\psi\colon H \rightarrow R$ is a homomorphism. Moreover, let $(C_x)_{x\in V(H)}$ be candidate sets such that $C_x\In V_{\psi(x)}$ and $|C_x|\ge \alpha n$ for all $x\in V(H)$. Then there exists an embedding $\phi\colon H\rightarrow G$ such that $\phi(x)\in C_x$ for all $x\in V(H)$.
\end{lemma}

The following definition will be crucial for our embeddings.

We call $G[A,B]$ \defn{weakly-$(\alpha,\eps)$-super-regular} if
\begin{itemize}
\item $G[A,B]$ is $\eps$-regular with density at least $\alpha$;
\item for all $a\in A$, $d_G(a,B)\ge \alpha|B|$ or $d_G(a,B)=0$;
\item for all $b\in B$, $d_G(b,A)\ge \alpha|A|$ or $d_G(b,A)=0$.
\end{itemize}

The next proposition implies that we can turn any $\eps$-regular pair into a weakly-super-regular pair by deleting a small number of edges. This will allow us to simultaneously turn all regular pairs of a regularity partition obtained by the regularity lemma into weakly-super-regular pairs (which is impossible for the more standard notion of super-regularity as one needs to delete vertices in that case).
 
\begin{prop} \label{prop:clean regularity pair}
Let $G[A,B]$ be $\eps$-regular with density at least $\alpha$ and assume that $|A|=|B|=:m$. Then $G[A,B]$ can be made weakly-$(\alpha-2\sqrt{\eps},4\sqrt{\eps})$-super-regular by deleting at most $(\alpha-\eps) m$ edges at every vertex. 
\end{prop}

\COMMENT{
\proof
Let $G$ be the underlying bipartite graph. We may assume that $\alpha=d_G(A,B)$ and $2\eps \le \alpha \le 1$. Let $Z^{bad}:=\set{z\in A\cup B}{d_G(z) < (\alpha-\eps)m}$. Then, $A^{bad}:=A \cap Z^{bad}$ and $B^{bad}:=B \cap Z^{bad}$ are both of size at most $\eps m$. Let $G'$ be obtained from $G$ by deleting all edges incident to a vertex in $Z^{bad}$. So if $z \in Z^{bad}$, then $d_{G-G'}(z)=d_G(z)<(\alpha-\eps)m$, and if $z \notin Z^{bad}$, then $d_{G-G'}(z) \le \eps m$. We claim that $G'[A,B]$ is weakly-$(\alpha-2\sqrt{\eps},4\sqrt{\eps})$-super-regular. Let $\alpha':=d_{G'}(A,B)$. Clearly, if $z \in Z^{bad}$, then $d_{G'}(z)=0$. On the other hand, if $z \notin Z^{bad}$, then $d_{G'}(z) \ge d_G(z) - \eps m \ge (\alpha-2\eps) m \ge (\alpha-2\sqrt{\eps}) m$. Now, let $X \In A$ and $Y \In B$ with $|X|,|Y| \ge 4\sqrt{\eps}m$. We have $e_{G-G'}(X,Y) \le e(X,B^{bad}) + e(Y,A^{bad}) \le \sqrt{\eps}|X||Y|$ and thus $d_{G'}(X,Y) \ge d_G(X,Y) - \sqrt{\eps} \ge \alpha - 2\sqrt{\eps}$. In particular, $\alpha - 2\sqrt{\eps} \le \alpha' \le \alpha$. Hence, we have $|d_{G'}(X,Y) - \alpha'| \le |d_{G'}(X,Y) - \alpha| + |\alpha - \alpha'| \le 4\sqrt{\eps}$.
\endproof
}

Many of our constructions will be carried out in a graph where we have found a regularity partition and `cleaned' the edges between classes using Proposition~\ref{prop:clean regularity pair}. In order to describe such graphs, we use the following definition. Given a graph $G$, we call a partition $V_1,\dots,V_k$ of $V(G)$ an \defn{$(\alpha,\eps,k)$-partition} of $G$, if 
\begin{enumerate}[label=(P\arabic*)]
\item $|V_i| = (1\pm \eps)|G|/k$;
\item $G[V_i]$ is empty for every $i\in [k]$;
\item for all $1\le i < j \le k$, $G[V_i,V_j]$ is either weakly-$(\alpha,\eps)$-super-regular or empty.
\end{enumerate}

We will often use the fact that if $V_1,\dots,V_k$ is an $(\alpha,\eps,k)$-partition of $G$ and $G'$ is a spanning subgraph of $G$ such that $d_{G'}(x,V_i)\ge d_G(x,V_i) - \eps^2 |V_i|$ for all $x\in V(G)$ and $i\in[k]$, then $V_1,\dots,V_k$ is an $(\alpha -2\eps,3\eps,k)$-partition of $G'$.\COMMENT{(P1) and (P2) clearly hold. Consider $i,j$ and the pair $(V_i,V_j)$. For $x\in V_i$, if $d_G(x,V_i)=0$, then $d_{G'}(x,V_i)=0$. Otherwise, $d_{G'}(x,V_i) \ge (\alpha-\eps^2) |V_i|$. Let $\alpha_{ij}:=e_{G}(V_i,V_j)/|V_i||V_j|$ and $\alpha_{ij}':=e_{G'}(V_i,V_j)/|V_i||V_j|$. Consider $X\In V_i$ and $Y\In V_j$ with $|X|\ge 3\eps|V_i|$ and $|Y| \ge 3\eps|V_j|$. Then, $e_{G'}(X,Y)\le e_G(X,Y) \le (\alpha_{ij}+\eps)|X||Y|$ and $e_{G'}(X,Y) \ge (\alpha_{ij}-\eps)|X||Y| - |X|\eps^2|V_j| \ge (\alpha_{ij}-2\eps)|X||Y|$. In particular, $\alpha_{ij}-2\eps\le \alpha_{ij}'\le \alpha_{ij}$. Hence $(\alpha_{ij}'-2\eps)|X||Y|\le e_{G'}(X,Y) \le (\alpha_{ij}'+3\eps)|X||Y|$.}

Similarly, if $\eps\le 1/2$ and $V_i'\In V_i$ with $|V_i\sm V_i'|\le \eps^2|V_i|$ for all $i\in[k]$, then $V_1',\dots,V_k'$ is an $(\alpha-\eps,3\eps,k)$-partition of $G[V_1'\cup \dots \cup V_k']$.\COMMENT{Let $G':=G[V_1'\cup \dots \cup V_k']$. (P2) clearly holds. For (P1) we check $|V_i'|\ge(1-\eps^2)|V_i| \ge (1-\eps)(1-\eps^2)|G|/k\ge (1-3\eps)|G'|/k$ and $|V_i'|\le (1+\eps)|G'|/(1-\eps^2)k \le (1+3\eps)|G'|/k$. (P3) Consider $i,j$. Let $x\in V_i'$. If $d_{G'}(x,V_j')>0$, then $d_{G'}(x,V_j')\ge (\alpha-\eps^2)|V_j| \ge (\alpha-\eps)|V_j'|$. Slicing lemma with $c=(1-\eps^2)$ gives that $(V_i',V_j')$ is $2\eps/(1-\eps^2)\le 3\eps$-regular with density at least $\alpha-\eps$.}

\section{The general decomposition theorem} \label{sec:general dec}

In Section~\ref{sec:near optimal} we saw how to find a near-optimal $F$-decomposition of a given graph $G$ which covers all but a bounded number of edges of $G$. As mentioned in Section~\ref{sec:sketch}, our goal is to deal with the leftover edges using absorbers.
Given graphs $F$ and $H$, an \defn{$F$-absorber for $H$} is a graph $A$ such that $V(H)\In V(A)$ is independent in $A$ and both $A$ and $A\cup H$ have $F$-decompositions. For $H$ to have an $F$-absorber it is clearly necessary that $H$ is at least $F$-divisible. Clearly, if $G$ was $F$-divisible in the beginning, then the leftover graph $H$ obtained from a near-optimal $F$-decomposition of $G$ is still $F$-divisible. The strategy is thus to find an $F$-absorber $A$ for every possible leftover graph $H$ and remove all these absorbers before covering almost all edges of $G$. The union of these absorbers will then allow us to deal with any leftover from the near-optimal decomposition. In this section, we prove a general decomposition result, that is, $G$ has an $F$-decomposition whenever the minimum degree of $G$ is large enough to guarantee (i) an approximate decomposition, (ii) a covering of all edges at one vertex, (iii) an absorber for any bounded size subgraph (see Theorem~\ref{thm:general dec}).

In order to minimise the minimum degree we require to find absorbers, we will make sure that the possible leftover graphs respect an $(\alpha,\eps,k)$-partition. We will achieve this by applying the regularity lemma to $G$ first and then `cleaning' $G$ with respect to the obtained partition.
Call $F$ \defn{$\delta$-absorbing} if the following is true:
\begin{itemize}
\item[] Let $1/n \ll 1/k_0',\eps \ll \alpha, 1/b \ll 1/m, \mu ,1/|F|$ and suppose that $G$ is a graph on $n$ vertices with $\delta(G) \ge (\delta + \mu) n$ which has an $(\alpha,\eps,k)$-partition for some $k \le k_0'$, and that $H$ is any $F$-divisible subgraph of $G$ of order at most $m$. Then $G$ contains an $F$-absorber for $H$ of order at most $b$. 
\end{itemize}

The focus of Sections~\ref{sec:compressions}, \ref{sec:transform}, \ref{sec:switchers} and \ref{sec:absorbers} will be to find for a given graph $F$ the minimal $\delta$ such that $F$ is $\delta$-absorbing.

\begin{theorem} \label{thm:general dec}
Let $F$ be a graph. Let $\delta \ge \max\Set{\delta_F^{0+},\delta_F^{vx}}$, and assume that $F$ is $\delta$-absorbing. Then $\delta_F \le \delta$.
\end{theorem}

\proof
Let $1/n \ll 1/k_0' \ll \eps \ll 1/k_0 , \alpha, 1/b \ll 1/m' \ll \mu,1/|F|$. Let $G$ be an $F$-divisible graph on $n$ vertices with $\delta(G) \ge (\delta + 10\mu) n $. We need to show that $G$ has an $F$-decomposition.

Set $U_0:=V(G)$. Using Lemma~\ref{lem:chernoff}, it is easy to see that there exists a subset $U_1$ of size $\lfloor \mu n \rfloor$ such that for all $x\in V(G)$,
\begin{align}
d_G(x,U_1) &\ge (\delta + 9\mu)|U_1|,\label{eq:degreesU1} \\
d_G(x,U_0\sm U_1) &\ge (\delta + 9\mu)|U_0\sm U_1|.\label{eq:backdegreesU1}
\end{align}

Apply the regularity lemma (Lemma~\ref{lem:regularity}) to $G[U_1]$ to obtain a partition of $U_1$ into sets $V_0,V_1,\dots,V_k$ and a spanning subgraph $G'$ of $G[U_1]$ such that 
\begin{enumerate}[label=(R\arabic*)]
\item $k_0 \le k \le k_0'$;
\item $|V_0| \le \eps |U_1|$;
\item $|V_1|=\dots=|V_k|$;
\item $d_{G'}(x) \ge d_{G[U_1]}(x)-(\alpha+\eps)|U_1|$ for all $x \in U_1$;
\item $G'[V_i]$ is empty for all $i\in [k]$;
\item for all $1\le i < j \le k$, $G'[V_i,V_j]$ is either $\eps$-regular with density at least $\alpha$ or empty.
\end{enumerate}

Applying Proposition~\ref{prop:clean regularity pair} to every pair $G'[V_i,V_j]$ of density at least $\alpha$ yields a spanning subgraph $G_{cl}$ of $G'$ such that
\begin{enumerate}[label=(R\arabic*$'$)]
\item $d_{G_{cl}}(x) \ge d_{G[U_1]}(x)-2\alpha|U_1|$ for all $x \in U_1$; \label{cleaned:degrees}
\item $G_{cl}[V_i]$ is empty for all $i\in [k]$; \label{cleaned:inside}
\item for all $1\le i < j \le k$, $G_{cl}[V_i,V_j]$ is either weakly-$(\alpha-2\sqrt{\eps},4\sqrt{\eps})$-super-regular or empty. \label{cleaned:super-regular}
\end{enumerate}

Let $H:=G[U_1]-G_{cl}$.\COMMENT{The order is now irrelevant. We can first apply Proposition~\ref{prop:cover-sparse-graph} to absorb $G[U_1]-G_{cl}$, or do it after finding absorbers, as long as it happens before the near-optimal.}
 Note that $\Delta(H) \le 2\alpha|U_1| \le 2\alpha n$ and $\delta\ge \delta_F^{vx} \ge \delta_F^e$. So by (\ref{eq:backdegreesU1}) and Proposition~\ref{prop:cover-sparse-graph}, we can obtain a subgraph $A$ of $G$ such that $A[U_1]$ is empty, $A\cup H$ is $F$-decomposable and $\Delta(A) \le \mu^4 n$.

By (\ref{eq:degreesU1}) and \ref{cleaned:degrees},
\begin{align}
\delta(G_{cl}) \ge (\delta + 8\mu)|U_1|.\label{cleaned degree}
\end{align}

Let $U_2$ be a subset of $U_1\sm V_0$ of size $\lfloor \mu |U_1|\rfloor$ satisfying the following properties (consider a random choice of $U_2$ and apply Lemma~\ref{lem:chernoff}):
\begin{enumerate}[label=(\roman*)]
\item $d_{G_{cl}}(x,U_2) \ge (\delta + 7\mu)|U_2|$ for all $x\in U_1$; \label{cond:equal-cutting 1}
\item $|V_i\sm U_2| = (1-\mu \pm \eps)|V_i|$ for all $i \in [k]$; \label{cond:equal-cutting 2}
\item $d_{G_{cl}}(x,V_i\sm U_2) \ge (1-2\mu) d_{G_{cl}}(x,V_i)$ for all $x\in U_1$, $i\in [k]$. \label{cond:equal-cutting 3}
\end{enumerate}
Finally, let $U_2 \supseteq U_3 \supseteq \dots \supseteq U_\ell$ be a $(\delta+6\mu,\mu,m)$-vortex in $G_{cl}[U_2]$ for some $\lfloor \mu m' \rfloor \le m \le m'$, which exists by Lemma~\ref{lem:get vortex}. We want to find an absorber for every $F$-divisible subgraph of $G_{cl}[U_\ell]$. We will find these in the following graph $G_{abs}$: For every $i\in [k]$, let $V_i':=(V_i\sm U_2)\cup (V_i \cap U_\ell)$, and let $$G_{abs}:=G_{cl}[V_1'\cup\dots\cup V_k'].$$ Note that, since $U_\ell \In U_2 \In U_1\sm V_0$, we have that $U_\ell \In V(G_{abs})$. We claim that
\begin{itemize}
\item[(iv)] $V_1',\dots,V_k'$ is an $(\alpha/2,10\sqrt{\eps},k)$-partition of $G_{abs}$;
\item[(v)] $\delta(G_{abs}) \ge (\delta+6\mu)|G_{abs}|$.
\end{itemize}

To verify (v), recall \eqref{cleaned degree} and note that $|U_1\sm(V_1'\cup\dots\cup V_k')| \le |U_2|+|V_0| \le 2\mu |U_1|$. 

Now, we check (iv). $G_{abs}[V_i']$ is clearly empty for every $i\in [k]$, so condition (P2) is satisfied. By \ref{cond:equal-cutting 2}, we have $|V_i'|= (1-\mu \pm 2\eps)|V_i|$. Now, (R3) implies that $||V_i'|-|V_j'||\le 4\eps |U_1|/k$ for all $1\le i<j\le k$ and hence $|V_i'|=(1 \pm 5\eps)|G_{abs}|/k$, so (P1) is satisfied.

In order to establish condition (P3), consider $1\le i < j\le k$. If $G_{cl}[V_i,V_j]$ is empty, then $G_{abs}[V_i',V_j']$ is also empty. So let us assume that $G_{cl}[V_i,V_j]$ is weakly-$(\alpha-2\sqrt{\eps},4\sqrt{\eps})$-super-regular. We need to show that $G_{abs}[V_i',V_j']$ is weakly-$(\alpha/2,10\sqrt{\eps})$-super-regular. By Fact~\ref{fact:slicing},\COMMENT{Can assume $|V_i'|\ge 0.8|V_i|$ and then apply Slicing lemma with $c=0.8$.} $G_{abs}[V_i',V_j']$ is $10\sqrt{\eps}$-regular with density at least $\alpha-6\sqrt{\eps} \ge \alpha/2$. Let $x\in V_i'$ and suppose that $d_{G_{abs}}(x,V_j') >0$. Then, $d_{G_{cl}}(x,V_j) \ge (\alpha-2\sqrt{\eps})|V_j|$. Using \ref{cond:equal-cutting 3}, we can check that $$d_{G_{abs}}(x,V_j') \ge d_{G_{cl}}(x,V_j\sm U_2) \ge (1-2\mu)d_{G_{cl}}(x,V_j) \ge (1-2\mu)(\alpha-2\sqrt{\eps})|V_j| \ge \alpha|V_j'|/2.$$
This proves (iv).

We will now find absorbers for all possible leftover graphs inside $U_\ell$. Let therefore $H_1,\dots,H_s$ be an enumeration of all $F$-divisible spanning subgraphs of $G_{cl}[U_\ell]$. We want to find edge-disjoint subgraphs $A_1,\dots,A_s$ in $G_{abs}$ such that for all $i\in[s]$, $A_i[U_\ell]$ is empty and $A_i$ is an $F$-absorber for $H_i$ of order at most $b$. Suppose that for some $j\in[s]$, we have already chosen $A_1,\dots,A_{j-1}$. Let $\tilde{G}_j$ be the graph obtained from $G_{abs}$ by deleting all edges inside $U_\ell$ except those of $H_j$, and all edges of $A_1,\dots,A_{j-1}$. Since we deleted at most $m+sb \le m+2^{\binom{m}{2}}b$ edges at every vertex, (v) implies that $\delta(\tilde{G}_j) \ge (\delta+5\mu)|\tilde{G}_j|$. Also, (iv) implies that $V_1',\dots,V_k'$ is an $(\alpha/4,30\sqrt{\eps},k)$-partition of $\tilde{G}_j$. Since $F$ is $\delta$-absorbing, $\tilde{G}_j$ contains an $F$-absorber $A_j$ for $H_j$ of order at most $b$ and so that $A_j[U_\ell]$ is empty.

Let $A^\ast := A_1\cup\dots\cup A_s$.
Let $G_{app}:=G-(A\cup H)-A^\ast$. Observe that 
\begin{align}
G_{app}[U_2]=(G-H)[U_2]=G_{cl}[U_2]. \label{eq:inner is the same}
\end{align}
We want to apply Lemma~\ref{lem:near optimal} to $G_{app}$. Note that
\begin{align}
\Delta(A^\ast \cup A \cup H) \le 2\mu^4 n \le \mu |U_2|.\label{absorbers removed}
\end{align}
Clearly, $G_{app}$ is $F$-divisible and $\delta(G_{app}) \ge (\delta + 8 \mu)n$ by \eqref{absorbers removed}. We claim that $U_0 \supseteq U_1 \supseteq \dots \supseteq U_\ell$ is a $(\delta + 4\mu, \mu, m)$-vortex in $G_{app}$. Conditions (V1)--(V3) hold by construction. Moreover, for $i\ge 3$, we have $d_{G_{app}}(x,U_i) \ge (\delta + 6\mu)|U_i|$ for all $x \in U_{i-1}$ by (\ref{eq:inner is the same}). So let $i\in \Set{1,2}$. For all $x \in U_{i-1}$, we get
\begin{align*}
d_{G_{app}}(x,U_i) \overset{\eqref{absorbers removed}}{\ge} d_G(x,U_i) - \mu |U_2| \ge (\delta + 7\mu)|U_i| - \mu |U_2| \ge (\delta + 6\mu)|U_i|,
\end{align*}
where we use (\ref{eq:degreesU1}) if $i=1$ and \ref{cond:equal-cutting 1} if $i=2$.
Thus, by Lemma~\ref{lem:near optimal}, there exists a subgraph $H_\ell^\ast$ of $G_{app}[U_\ell]$ such that $G_{app}-H_\ell^\ast$ is $F$-decomposable. In particular, $H_\ell^\ast$ is $F$-divisible. Crucially, by (\ref{eq:inner is the same}), $G_{app}[U_\ell]=G_{cl}[U_\ell]$, so $H_\ell^\ast = H_{s'}$ for some $s'\in [s]$. Since $A_{s'}\cup H_{s'}$ has an $F$-decomposition and $A_i$ has an $F$-decomposition for every $i\in[s]\sm\Set{s'}$, we conclude that $A^\ast \cup H_\ell^\ast$ is $F$-decomposable. Therefore, $$G=(A\cup H) \cup A^\ast \cup (G_{app}-H_\ell^\ast) \cup H_\ell^\ast$$ has an $F$-decomposition.
\endproof

\section{Models and compressions} \label{sec:compressions}

In order to prove an upper bound on $\delta_F$ using Theorem~\ref{thm:general dec}, we must be able to find $F$-absorbers for rather arbitrary subgraphs of a large graph $G$. Suppose that $G'$ is a subgraph of $G$ and we know that $H$ is an $F$-absorber for $G'$. The \Erd--Stone theorem tells us that $\delta(G) \ge (1-1/(\chi(H)-1)+\mu)|G|$ guarantees a copy of $H$ in $G$. However, for $H$ to fulfill its purpose, we need it to be rooted at $V(G')$, which is more difficult to achieve. In the following, we will introduce the concepts that enable us to embed rooted graphs in an efficient way. 

\subsection{Models and labellings}

We define a \defn{model} to be a pair $(H,U)$ where $H$ is a graph and $U$ is an independent set in $H$. The vertices of $U$ are called \defn{roots}. Generally speaking, we want to embed $H$ into a large graph $G$ in such a way that the roots are mapped to a prescribed destination. More formally,
given a model $(H,U)$, a graph $G$ and a map $\Lambda\colon U\rightarrow 2^{V(G)}$,
an \emph{embedding of $(H,U)$ into $G$ respecting $\Lambda$} is an injective homomorphism $\phi\colon H \to G$ such that $\phi(u)\in \Lambda(u)$ for all $u\in U$. Clearly, a necessary condition for the existence of an embedding respecting $\Lambda$ is that there exist distinct $(v_u)_{u\in U}\in V(G)$ with $v_u\in \Lambda(u)$. If this is satisfied, then we call $\Lambda$ a \defn{$G$-labelling of $U$}. All our labellings will be of the following form. There will be a set $U_1\In U$ such that $|\Lambda(u)|=1$ for all $u\in U_1$ (and thus $\Lambda(u)\neq \Lambda(u')$ for all distinct $u,u'\in U_1$). The labels of the elements of $U\sm U_1$ will be large subsets of clusters of a regularity partition. Furthermore, most often we will have $U_1=U$.

Call a model $(H,U)$ \defn{$\delta$-embeddable}, if for all $\mu>0$, there exists an $n_0\in\bN$ such that whenever $G$ is a graph on $n\ge n_0$ vertices with $\delta(G)\ge (\delta+\mu)n$ and $\Lambda$ is any $G$-labelling of $U$, there exists an embedding of $H$ into $G$ respecting $\Lambda$. As described below, the degeneracy of $H$ rooted at $U$ can be used to give a simple bound on the values of $\delta$ for which $(H,U)$ is $\delta$-embeddable.
Here, for a graph $K$, the \defn{degeneracy of $K$ rooted at $X\In V(K)$} is the smallest $d\in \bN\cup\Set{0}$ such that there exists an ordering $v_1,\dots,v_{|K|-|X|}$ of the vertices of $V(K)\sm X$ such that for all $i\in[|K|-|X|]$, $$d_K(v_i,X\cup\set{v_j}{1\le j<i})\le d.$$
Observe that if $H$ has degeneracy at most $d$ rooted at $U$ for some $d\in \bN$, then $(H,U)$ is $(1-1/d)$-embeddable. Indeed, if $G$ is a graph with $\delta(G)\ge (1-1/d+\mu)|G|$, then every set of $d$ vertices has many common neighbours. Hence, any set of $|U|$ vertices in $G$ can be extended to a copy of $H$ by embedding the vertices of $V(H)\sm U$ one by one.

As we shall discuss in the next subsection, this simple bound is usually not sufficient for our purposes. Before that, we will prove the following lemma which is used in Sections~\ref{sec:lower bounds} and~\ref{sec:deltafa}. It says that if the minimum degree of $G$ is sufficiently large to embed a copy of $(H,U)$ according to any given labelling, then we can in fact embed many copies edge-disjointly into $G$, provided that the respective labellings are not overly restrictive. 

\begin{lemma} \label{lem:finding}
Let $(H,U)$ be a $\delta$-embeddable model for some $\delta\in[0,1]$. Let $1/n\ll \mu \ll 1/|H|$ and let $G$ be a graph on $n$ vertices with $\delta(G)\ge (\delta+\mu)n$. Suppose that $\Lambda_1,\dots,\Lambda_m$ are $G$-labellings of $U$ such that $m\le \mu^4 n^2$ and for all $v\in V(G)$, $|\set{j\in[m]}{v\in \bigcup Im(\Lambda_j)}| \le \mu^2 n$. Then there exist embeddings $\phi_1,\dots,\phi_m$ of $H$ into $G$ such that $\phi_j$ respects $\Lambda_j$ for all $j\in[m]$, $\phi_1(H),\dots,\phi_m(H)$ are edge-disjoint, and $\Delta(\bigcup_{j\in [m]}\phi_j(H))\le \mu n$.
\end{lemma}

\proof
We may assume that $|\Lambda_j(u)|=1$ for all $j\in[m]$ and $u\in U$. Let $R_j:=\bigcup_{u\in U}\Lambda_j(u)$. We will find $\phi_1,\dots,\phi_m$ one by one. For $j\in [m]$ and $v\in V(G)$, define $root(v,j):=|\set{j'\in[j]}{v\in R_{j'}}|$. Suppose that for some $j\in[m]$, we have already defined $\phi_1,\dots,\phi_{j-1}$ such that
\begin{align}
d_{G_{j-1}}(v)\le \mu^2 n + (root(v,j-1)+1)|H| \label{eq:bad vertices}
\end{align}
for all $v\in V(G)$, where $G_{j-1}:=\bigcup_{j'\in[j-1]}\phi_{j'}(H)$. We now want to define $\phi_j$ such that (\ref{eq:bad vertices}) holds with $j$ replaced by $j+1$. Let $BAD:=\set{v\in V(G)}{d_{G_{j-1}}(v)>\mu^2 n}$. Note that $2e(G_{j-1})\ge |BAD|\mu^2 n$ and $e(G_{j-1})\le me(H)\le \mu^4e(H)n^2$. Thus, $|BAD|\le 2\mu^2 e(H) n \le \mu n/4$. Let $BAD':=BAD\sm R_{j}$ and define $G':=(G-G_{j-1})[V(G)\sm BAD']$. Since $\Delta(G_{j-1})\le \mu n/4$, we have $\delta(G')\ge (\delta+\mu/2)n$. Note that $\Lambda$ is a $G'$-labelling of $(H,U)$. Since $(H,U)$ is $\delta$-embeddable, there exists an embedding $\phi_j$ of $H$ into $G'$ respecting $\Lambda_j$. It remains to show that (\ref{eq:bad vertices}) holds with $j$ replaced by $j+1$. Let $v\in V(G)$. If $v\notin BAD$, then $d_{G_j}(v)\le d_{G_{j-1}}(v)+|H|\le \mu^2 n+|H|$, so assume that $v\in BAD$. If $v\in BAD'$, then $d_{\phi_j(H)}(v)=0$ and hence $d_{G_j}(v)\le d_{G_{j-1}}(v)$. Finally, if $v\in BAD\sm BAD'$, then $v\in R_{j}$, implying that $root(v,j)=root(v,j-1)+1$. Thus, $d_{G_j}(v)\le d_{G_{j-1}}(v)+|H|\le \mu^2 n + (root(v,j)+1)|H|$.
\endproof

\subsection{Compressions}

As observed in the previous section, the rooted degeneracy of a model $(H,U)$ translates into a bound on the minimum degree of a graph $G$ which ensures the embeddability of $(H,U)$. However, this bound is usually too large for our purposes. We will improve on this by considering the above degeneracy approach in the reduced graph setting. This involves the notion of \emph{compressions} which we now discuss informally. 

Let $(H,U)$ be a model. Let $G$ be a graph and $\Lambda$ a $G$-labelling of $U$, and as an example assume that $U=\Set{u_1,u_2}$ and $\Lambda(u_i)=\Set{v_i}$. Suppose that $R$ is the reduced graph of some regularity partition of $G$ with cluster function $\sigma$. The key lemma tells us that if there exists a homomorphism $H\rightarrow R$, then $H$ can be embedded into $G$, but the corresponding embedding does not necessarily respect the given labelling $\Lambda$. It will be useful to have an embedding $\phi$ of $H$ such that, though $\phi$ might not respect $\Lambda$,
it is close to respecting $\Lambda$ (in the sense that at least every vertex is embedded in the correct cluster). More precisely, letting $x_1',x_2'\in V(R)$ be such that $\sigma(v_i)=x_i'$, we require $\sigma(\phi(u_i))=x_i'$.
Now, if $\psi'\colon H\rightarrow R$ was a homomorphism such that $\psi'(u_i)=x_i'$, then the key lemma again would give us such an embedding.

In order to investigate the embedding properties of $(H,U)$ without referring to the host graph and labelling, we consider the following intermediate graph: By considering only $(H,U)$, we find a graph $K$ such that $\psi\colon H\rightarrow K$ is a homomorphism, $\Set{x_1,x_2}$ is independent in $K$ and $K$ has low degeneracy rooted at $\Set{x_1,x_2}$, where $x_i:=\psi(u_i)$ are distinct. Then, coming back to the above situation, we can use a degeneracy argument to find a copy of $K$ in $R$ rooted at $\Set{x_1',x_2'}$, i.e., a homomorphism $\psi''\colon K\rightarrow R$ such that $\psi''(x_i)=x_i'$. We can then take $\psi':=\psi''\circ \psi$ to find a homomorphism from $H$ into $R$.
In finding this homomorphism, we have used the degeneracy of $K$, not $H$. We can easily pick such a graph $K$ with degeneracy at most $\chi(H)+|U|-1$, and sometimes even lower. In many cases this is much lower than the degeneracy of $H$, reducing the minimum degree required.
Loosely speaking, our gain has come from the fact that we `compressed' the original embedding problem to the reduced graph level.

Note that we assumed $\Set{x_1,x_2}$ is independent. If $x_1'x_2'\in E(R)$, then we could allow $K$ to contain the edge $x_1x_2$, possibly allowing a graph $K$ with lower degeneracy to be used. Similarly, if $x_1'=x_2'$, we can construct $K$ with $x_1=x_2$, which may allow us to use a graph $K$ with lower degeneracy.\COMMENT{In order to exploit these advantages while avoiding explicit reference to $G$ and $R$, we consider the possible scenarios separately. More precisely, assuming that $J$ is either an isolated vertex, two isolated vertices or an edge, we find $K$ and $\psi$ such that $\psi\colon H\rightarrow K$ is a homomorphism with $K[\psi(U)]=J$ and $K$ has low degeneracy rooted at $\psi(U)$. If $U$ is larger than in our example, we also need to keep track of how the root vertices are distributed among $V(J)$.}

The following definitions formalise the above discussion.
Given a set $U$, a pair $(J,f)$ is called a \defn{root-compression of $U$}, if $J$ is a graph and $f\colon U \rightarrow V(J)$ is a surjective map. 
Further, $(J,f,K,\psi)$ is called a \defn{compression of $(H,U)$} if
\begin{enumerate}[label=(C\arabic*)]
\item $K$ is a graph and $J$ is an induced subgraph of $K$; \label{arr:shapes}
\item $(J,f)$ is a root-compression of $U$; \label{arr:roots}
\item $\psi\colon H\rightarrow K$ is a homomorphism such that $\psi{\restriction_U}=f$. \label{arr:hom}
\end{enumerate}

We say that $(J,f,K,\psi)$ is a \defn{$d$-compression} if $K$ has degeneracy at most $d$ rooted at $V(J)$.
For our purposes, $J$ will be a rather simple graph, e.g.~a short cycle or path, whereas $K$ may have a more complex structure. It is thus often convenient to ignore $K$ itself and only record its degeneracy. Accordingly, $(H,U)$ is called \defn{$d$-compressible with respect to $(J,f)$}, if there exist $K$ and $\psi$ such that $(J,f,K,\psi)$ is a $d$-compression of $(H,U)$.
When referring to a compression $(J,f,K,\psi)$, $J$ and $f$ are technically redundant since $J=K[\psi(U)]$ and $f=\psi{\restriction_U}$. However, as indicated above, the root-compression $(J,f)$ is often the essential part of $(J,f,K,\psi)$ and we want to be able to refer to it directly.

Our embedding approach is roughly as follows (this is formalised in Lemma~\ref{lem:rooted embedding}). Assume that $(H,U)$ is $d$-compressible with respect to $(J,f)$. Let $G$ be a graph with $\delta(G) \ge (1-1/d+\mu)|G|$ and let $R$ be the reduced graph of some regularity partition of $G$. Instead of seeking an embedding approach which works for all possible $G$-labellings $\Lambda$ of $U$, we limit ourselves to `admissible' labellings. Informally, this means that $J$ is a subgraph of $R$ such that every label $\Lambda(u)$ is contained in the cluster of $f(u)$. By definition, there exist $K$ and $\psi$ such that $(J,f,K,\psi)$ is a $d$-compression of $(H,U)$. As $R$ inherits the minimum degree ratio of $G$ and since $K$ has degeneracy at most $d$ rooted at $V(J)$, we can extend $J$ to a copy of $K$ in $R$. Then, using the key lemma, we can embed $H$ into $G$ such that $u\in U$ is mapped into the cluster of $f(u)$. This embedding is close to what we desire in the sense that the image of $u$ and the label of $u$ are already in the same cluster. Finally, we will modify this embedding to an embedding respecting $\Lambda$.\COMMENT{Note that we face a trade-off between $J$ and $d$ here. Since $d$ directly relates to the minimum degree that we have to presume in $G$, the ultimate goal is to minimise $d$. Assume that the order of $J$ is fixed. Then, intuitively, the more edges we have in $J$, the less we have to attach to $J$ to obtain a $K$ which meets our requirements, while keeping $d$ small.
However, the denser $J$ is, the less labellings are admissible, possibly rendering our approach vacuous.}
It remains to make more precise when a labelling is admissible. We can then prove our embedding lemma.

For this, we say that a $G$-labelling $\Lambda$ of $U$ \defn{respects} the root-compression $(J,f)$, if $\Lambda(u)=\Lambda(u')$ for all $u,u'\in U$ with $f(u)=f(u')$. In this case, the function $\Lambda_J\colon V(J)\rightarrow 2^{V(G)}$, where $\Lambda_J(x):=\Lambda(u)$ for any $u\in f^{-1}(x)$, is well-defined.

The following definition collects a set of restrictions that we put on a $G$-labelling $\Lambda$. They will later enable us to find an embedding of a given model $(H,U)$ into $G$ respecting $\Lambda$. Note that the conditions do not involve $H$, but depend on the chosen root-compression $(J,f)$ of $U$.  
\begin{defin} \label{def:admissible}
Let $U$ be a set with root-compression $(J,f)$, $G$ a graph and $\Lambda$ a $G$-labelling of $U$. Let $U_1:=\set{u\in U}{|\Lambda(u)|=1}$ and $U_2:=U\sm U_1$. Moreover, let $J_1:=f(U_1)$ and $J_2:=V(J)\sm J_1$. We call $\Lambda$ \underline{$(\alpha,\eps,k)$-admissible}, if $\Lambda$ respects $(J,f)$ and $G$ has an $(\alpha,\eps,k)$-partition\COMMENT{This can be weakened in the sense that in (P3) it would be enough if $G[V_i,V_j]$ is $\eps$-regular with density at least $d$ or empty.} $V_1,\dots,V_k$ with reduced graph $R$ such that the following hold:
\begin{enumerate}[label=(D\arabic*)]
\item there exists a homomorphism $j\colon J\rightarrow R$ such that $\Lambda_J(x) \In V_{j(x)}$ for all $x\in V(J)$;\label{adm:hom}
\item $|\Lambda_J(x)|\ge \alpha|V_{j(x)}|$ for all $x\in J_2$;\label{adm:label size}
\item for every $x\in V(J)$, $d_G(W_x,V_{j(x)}) \ge \alpha |V_{j(x)}|$, where $W_x:=\bigcup_{y\in N_J(x,J_1)}\Lambda_J(y)$;\label{adm:conform}
\item $\Lambda_J(x) \In N_G(W_x,V_{j(x)})$ for all $x\in J_2$.\label{adm:label loc}
\end{enumerate}
\end{defin}

In order to prove the main lemma of this section (i.e.~Lemma~\ref{lem:rooted embedding}), we need the following simple result.

\begin{prop} \label{prop:jumping}
Let $\alpha\in(0,1)$. Let $G$ be a graph and let $V_1,\dots,V_k$ be a partition of $V(G)$ such that $|V_i|\ge 1/\alpha$ for all $i\in[k]$. Let $W\In V(G)$. Then there exists a spanning subgraph $G'$ of $G$ such that $d_{G'}(x)\ge d_G(x)- \alpha 2^{2^{|W|}+1}|G|$ for every $x\in V(G)$ and $d_{G'}(W',V_i)\notin (0, \alpha |V_i|)$ for every subset $W'\In W$ and every $i\in[k]$.
\end{prop}

\proof
Let $n:=|G|$ and let $W_1,\dots,W_s$ be an enumeration of the subsets of $W$. For each $j\in[s]$, let $\beta_j:=2^{s-j}$, $\gamma_j:=\beta_1+\dots+\beta_j$ and $\gamma_0:=0$. Suppose that for some $j\in[s]$, we have found a subgraph $G_{j-1}$ by deleting at most $\gamma_{j-1}\alpha n$ edges at every vertex such that $d_{G_{j-1}}(W_\ell,V_i) \notin (0, \beta_{j-1}\alpha |V_i|)$ for every $i\in[k]$ and $\ell\in[j-1]$. Consider $W_j$ and fix some $w\in W_j$. Let $I\In [k]$ contain the indices $i$ for which $d_{G_{j-1}}(W_j,V_i) < \beta_j \alpha |V_i|$. Let $G_j$ be obtained from $G_{j-1}$ by deleting for every $i\in I$ the edges from $w$ to $N_{G_{j-1}}(W_j,V_i)$. Clearly, for all $i\in[k]$, $d_{G_{j}}(W_j,V_i) \notin (0, \beta_{j}\alpha |V_i|)$, and $G_j$ is obtained from $G_{j-1}$ by deleting at most $\beta_j\alpha n$ edges at every vertex. It remains to show that for all $\ell\in[j-1]$ and $i\in[k]$, we have $d_{G_j}(W_\ell,V_i)\notin (0,\beta_j\alpha |V_i|)$. Suppose that $\ell\in[j-1]$, $i\in[k]$ and $d_{G_j}(W_\ell,V_i)>0$. If $w\in W_\ell$, then $d_{G_j}(W_\ell,V_i)\ge d_{G_{j-1}}(W_\ell,V_i)-\beta_j \alpha |V_i| \ge (\beta_{j-1}-\beta_j)\alpha |V_i| =\beta_j\alpha |V_i|$. If $w\notin W_\ell$, then we have $d_{G_j}(W_\ell,V_i)\ge d_{G_{j-1}}(W_\ell,V_i)-1\ge (\beta_{j-1}-\beta_j)\alpha |V_i| =\beta_j\alpha |V_i|$ as well.

Let $G':=G_s$. Note that $\beta_s=1$ and $\gamma_s\le 2^{s+1}$.\COMMENT{Of course $2^{2^{|W|}+1}$ could be improved tremendously here.}
\endproof

\begin{lemma} \label{lem:rooted embedding}
Let $1/n \ll 1/k_0', \eps \ll \alpha \ll 1/t,\mu$. Let $(H,U)$ be a $d$-compressible model with respect to the root-compression $(J,f)$ such that $|H|\le t$. Suppose that $G$ is a graph on $n$ vertices with $\delta(G)\ge (1-1/d+\mu)n$ and $\Lambda$ is an $(\alpha,\eps,k)$-admissible $G$-labelling of $U$ for some $k\le k_0'$. Then there exists an embedding of $(H,U)$ into $G$ respecting $\Lambda$.
\end{lemma}

\proof
Let $U_1,U_2,J_1,J_2,V_1,\dots,V_k,R,j,(W_x)_{x\in V(J)}$ be as in Definition~\ref{def:admissible}. Let $\sigma\colon V(G)\rightarrow V(R)$ be the cluster function of $V_1,\dots,V_k$. Let $K$ and $\psi$ be such that $(J,f,K,\psi)$ is a $d$-compression of $(H,U)$. We may assume that $|K|\le t$.\COMMENT{We can delete vertices of $K\sm Im(\psi)$. This does not increase degeneracy. Finally, if $\psi$ is surjective, then $|K|\le |H|$.} We define $W_x$ also for $x\in V(K)\sm V(J)$, that is, $W_x:=\bigcup_{y\in N_K(x,J_1)}\Lambda_J(y)$.

Let $W:=\bigcup_{u\in U_1}\Lambda(u)$ and apply Proposition~\ref{prop:jumping} to obtain a spanning subgraph $G'$ of $G$ such that $\delta(G') \ge (1-1/d+\mu/2)n$ and $d_{G'}(W',V_i)\notin (0, \alpha |V_i|)$ for every subset $W'\In W$ and every $i\in[k]$. For every $x\in J_1$, let $v_x\in V(G)$ be such that $\Lambda_J(x)=\Set{v_x}$. So $W=\set{v_x}{x\in J_1}$. For every $x\in J_2$, pick any vertex $v_x \in V_{j(x)}$.

We are going to define a homomorphism $\hat{\xi}\colon K\to R$ such that $\hat{\xi}{\restriction_{V(J)}}=j$ and $N_G ( W_x , V_{ \hat{\xi}(x) }) \neq \emptyset $ for all $x \in V(K) \sm V(J)$. First we define a homomorphism $\xi\colon (K-J) \rightarrow G'$ as follows.
For every $x\in V(J)$, let $\xi(x):=v_x$. Let $x_1,\dots,x_\ell$ be an ordering of the vertices of $V(K)\sm V(J)$ such that for every $i\in[\ell]$, we have $|N^{<}_K(x_i)|\le d$, where $N^{<}_K(x_i):= N_K(x_i)\cap (V(J)\cup \set{x_j}{j<i})$. We define $\xi(x_i)$ one by one. Suppose that for some $i\in[\ell]$, we have already defined $\xi(x_1),\dots,\xi(x_{i-1})$. Since $\delta(G')\ge (1-1/d+\mu/2)n$, we know that $N_{G'}(\xi(N_K^{<}(x_i)))\neq \emptyset$, so we can pick $\xi(x_i)$ from this set.
Note that for all $i \in [\ell]$, we have $W_{x_i} \subseteq \xi (N_K^<(x_i))$\COMMENT{If $z\in W_{x_i}$, then $\exists y\in N_K(x_i,J_1):z=v_y$. So $y\in N_K^<(x_i)$ and $\xi(y)=v_y=z$, so $z\in \xi(N_K^<(x_i))$} and thus
\begin{equation} \label{eq:hello}
\xi(x_i) \in N_{G'}(W_{x_i}).
\end{equation}
We now obtain a homomorphism $\hat{\xi}\colon K \rightarrow R$, where we let $\hat{\xi}(x):=\sigma(\xi(x))$ for all $x\in V(K)$. Note that $\hat{\xi}(x)=j(x)$ for all $x\in V(J)$,\COMMENT{$\hat{\xi}(x)=\sigma(\xi(x))=\sigma(v_x)=j(x)$} therefore, $\hat{\xi}$ is indeed a homomorphism.\COMMENT{Let $xy\in E(K)$. If $xy\in E(J)$, then $\hat{\xi}(x)\hat{\xi}(y)\in E(R)$ because $j$ is a homomorphism. If $xy\notin E(J)$, then $\xi(x)\xi(y)\in E(G')$ and hence $\hat{\xi}(x)\hat{\xi}(y)\in E(R)$. Here, we use that the partition is cleaned, that is, there are only edges between $V_i$ and $V_j$ if there is a respective edge in $R$.}

For every vertex $x\in V(K)$, we will construct a set $Z_x\In V(G)$ such that (Z1)--(Z5) below hold.
These sets $Z_x$ will then be suitable for an application of the key lemma.
\begin{enumerate}[label=(Z\arabic*)]
\item The $Z_x$'s are pairwise disjoint and disjoint from $W$;
\item $|Z_x|\ge \alpha n/2kt$ for all $x\in V(K)$;
\item if $xy\in E(K)$, then $G[Z_x,Z_y]$ is $\sqrt{\eps}$-regular with density at least $\alpha /2$;
\item for all $x\in J_2$, $Z_x\In\Lambda_J(x)$;
\item for all $xy\in E(K)$ with $x\in J_1$, $Z_y\In N_{G}(v_x)$.
\end{enumerate}

In order to achieve this, we first define a set $Z'_x$ for every $x\in V(K)$ such that 
\begin{enumerate}[label=(Z\arabic*$'$)]
\item $Z'_x\In V_{\hat{\xi}(x)}$ and $|Z'_x|\ge \alpha |V_{\hat{\xi}(x)}|$ for all $x\in V(K)$; \label{eq:Z1'}
\item for all $x\in J_2$, $Z'_x=\Lambda_J(x)$; \label{eq:Z2'}
\item for all $xy\in E(K)$ with $x\in J_1$, $Z'_y\In N_{G}(v_x)$. \label{eq:Z3'}
\end{enumerate}

For $x\in J_2$, we let $Z_x':=\Lambda_J(x)$. For $x\in V(K)\sm J_2$, we let $Z_x':=N_G(W_x,V_{\hat{\xi}(x)})$. Thus, \ref{eq:Z2'} holds by definition.

We now check \ref{eq:Z1'}. Recall that $\hat{\xi}(x)=j(x)$ for all $x\in V(J)$, so we clearly have $Z'_x\In V_{\hat{\xi}(x)}$ for all $x\in V(K)$. Now, if $x\in J_1$, then $|Z'_x|\ge \alpha |V_{\hat{\xi}(x)}|$ holds by \ref{adm:conform}. If $x\in J_2$, then it holds by \ref{adm:label size}. So let $x\in V(K)\sm V(J)$. Since $W_x\In W$, we have that $d_{G'}(W_x,V_{\hat{\xi}(x)}) \notin (0,\alpha |V_{\hat{\xi}(x)}|)$. But $\xi(x)\in N_{G'}(W_x)$ by~(\ref{eq:hello}) and $\xi(x) \in V_{\hat{\xi}(x)}$ by definition of $\hat{\xi}$. Thus $N_{G'}(W_x,V_{\hat{\xi}(x)})$ is non-empty and we deduce that $|Z_x'|=d_{G}(W_x,V_{\hat{\xi}(x)})\ge \alpha |V_{\hat{\xi}(x)}|$.

We continue with checking \ref{eq:Z3'}. Let $xy\in E(K)$ and $x\in J_1$. Note that $v_x\in W_y$ by definition of $W_y$. Moreover, note that $Z_y'\In N_G(W_y)$, which follows directly from the definition of $Z_y'$ if $y\notin J_2$, and from \ref{adm:label loc} if $y\in J_2$. So $Z'_y\In N_{G}(v_x)$ follows from that.

It is now relatively easy to obtain $Z_x$ from $Z_x'$. Since there are $|K|\le t$ sets and $|W|\le t$, we can choose a subset $Z_x\In Z_x'$ for every $x\in V(K)$ such that they satisfy (Z1) and $|Z_x|\ge |Z_x'|/t-1 \ge 2\alpha |V_{\hat{\xi}(x)}|/3t$. So (Z2), (Z4) and (Z5) are also satisfied. Finally, since $\hat{\xi}$ is a homomorphism, $(V_{\hat{\xi}(x)},V_{\hat{\xi}(y)})$ is $\eps$-regular with density at least $\alpha$ whenever $xy\in E(K)$. Fact~\ref{fact:slicing} then implies (Z3).\COMMENT{applied with $c:=2\alpha/3t>\eps$. $Z_x\In V_{\hat{\xi}(x)}$ and $|Z_x|\ge c|V_{\hat{\xi}(x)}|$. $2\eps/c \le \sqrt{\eps}$ and $\alpha-\eps \ge \alpha/2$.}

Since $\psi\colon H\rightarrow K$ is a homomorphism by (C3), we can apply the key lemma (Lemma~\ref{lem:key lemma}) to obtain an injective homomorphism $\phi'\colon H \rightarrow G$ such that $\phi'(h)\in Z_{\psi(h)}$ for all $h\in V(H)$. Note that $\phi'(u)\in V_{j(f(u))}$ for every $u\in U$.\COMMENT{$\phi'(u)\in Z_{\psi(u)}\In Z'_{\psi(u)}\In V_{\hat{\xi}(\psi(u))} = V_{j(f(u))}$}
We define $\phi\colon H \rightarrow G$ by taking $\phi(u):=v_{f(u)}$ for all $u\in U_1$ and $\phi(h):=\phi'(h)$ for all $h\in V(H)\sm U_1$. We claim that $\phi$ is an embedding of $H$ into $G$ respecting $\Lambda$. The map $\phi$ is injective because of (Z1). For $u\in U_1$, we have $\phi(u)\in \Lambda(u)$ by definition. For $u\in U_2$, we have $\phi(u)=\phi'(u)\in Z_{\psi(u)} = Z_{f(u)} \In \Lambda(u)$ by (Z4). Finally, let $uh\in E(H)$ with $u\in U_1$ and $h\notin U_1$. Then $Z_{\psi(h)}\In N_G(v_{f(u)})$ by (Z5) and $\phi(h)=\phi'(h)\in Z_{\psi(h)}$. Hence, $\phi(u)\phi(h)\in E(G)$, completing the proof.
\endproof

\subsection{Attaching models}

The previous lemma provides us with a tool to embed models respecting given labellings. In general, we are interested in the minimum degree threshold at which this is possible. Note that the condition on $\delta(G)$ in Lemma~\ref{lem:rooted embedding} is governed by the degeneracy of the model $(H,U)$ with respect to $(J,f)$. Later on, we will try to find models with good, that is, low-degenerate, compressions.

We conclude this section by collecting some tools that will enable us to build those models in a modular way.

\begin{fact} \label{fact:model shrink}
Let $(H,U)$ be $d$-compressible with respect to $(J,f)$ and let $\beta\colon J\rightarrow J'$ be a surjective homomorphism. Then $(H,U)$ is $d$-compressible with respect to $(J',f')$, where $f'(u):=\beta(f(u))$ for all $u\in U$.
\end{fact}

\proof
Since $\beta$ is surjective, $(J',f')$ is a valid root-compression of $U$. Let $(J,f,K,\psi)$ be a $d$-compression of $(H,U)$. We may assume that $J'$ and $K$ are vertex-disjoint. Let $K'$ be the graph obtained from the union of $J'$ and $K\sm V(J)$ by adding an edge between $x'\in V(J')$ and $y\in V(K)\sm V(J)$ if there exists $x\in V(J)$ such that $\beta(x)=x'$ and $xy\in E(K)$. Clearly, $J'$ is an induced subgraph of $K'$. Moreover, the degeneracy of $K'$ rooted at $J'$ is at most $d$ since $K'\sm V(J')=K\sm V(J)$ and $d_{K'}(y,V(J'))\le d_{K}(y,V(J))$ for all $y\in V(K)\sm V(J)$, so we can take the same order of the vertices. Define $\psi'\colon V(H)\rightarrow V(K')$ as follows. If $\psi(v)\notin V(J)$, let $\psi'(v):=\psi(v)$. If $\psi(v)\in V(J)$, let $\psi'(v):=\beta(\psi(v))$. Clearly, $\psi'(u)=f'(u)$ for all $u\in U$ and $\psi$ is a homomorphism.
\endproof

Note that whenever $(H,U)$ and $(S,W)$ are models such that $V(H)\cap V(S) =W$, then $(H\cup S,U)$ is a model too.

\begin{prop} \label{prop:model attaching}
Let $(H,U)$ and $(S,W)$ be models such that $V(H)\cap V(S) =W$. Assume that $(J,f,K,\psi)$ is a $d$-compression of $(H,U)$ and that $(S,W)$ is $d$-compressible with respect to $(J_S,f_S)$. Suppose that $\beta\colon J_S\rightarrow K$ is a homomorphism that satisfies $\beta(f_S(w))=\psi(w)$ for all $w\in W$. Then there exists a $d$-compression $(J,f,K',\psi')$ of $(H\cup S,U)$ such that $K\In K'$ and $\psi'{\restriction_{V(H)}}=\psi$.
\end{prop}

\proof
By Fact~\ref{fact:model shrink}, we can assume that $J_S$ is a subgraph of $K$ satisfying $f_S(w)=\psi(w)$ for all $w\in W$, that is, $\beta$ is the identity. (Indeed, define $J':=\beta(J_S)$ and $f'(w):=\beta(f_S(w))$ for all $w\in W$. Then $J'$ is a subgraph of $K$ and $f'(w)=\beta(f_S(w))=\psi(w)$, and $(S,W)$ is $d$-compressible with respect to $(J',f')$ by Fact~\ref{fact:model shrink}.)

So let $(J_S,f_S,K_S,\psi_S)$ be a $d$-compression of $(S,W)$ and assume that $V(K_S)\cap V(K)=V(J_S)$. We can then take $K':=K\cup K_S$. So $J$ is an induced subgraph of $K'$ and $K'$ has degeneracy at most $d$ rooted at $V(J)$. Since $\psi_S(w)=f_S(w)=\psi(w)$ for all $w\in W$, we can define $\psi'(x):=\psi(x)$ for all $x\in V(H)$ and $\psi'(x):=\psi_S(x)$ for all $x\in V(S)$ in order to obtain a homomorphism $\psi'\colon H\cup S \rightarrow K'$.
\endproof

The fact that $K\In K'$ and $\psi'{\restriction_{V(H)}}=\psi$ in Proposition~\ref{prop:model attaching} allows us to attach several models to a given initial model without interference.

\begin{cor} \label{cor:model attaching}
Let $(H,U)$ be a model with $d$-compression $(J,f,K,\psi)$ and let $(S_1,W_1),\dots,(S_t,W_t)$ be models such that $V(S_i)\cap V(S_j)=W_i\cap W_j$ for all $1\le i<j\le t$. For every $i\in[t]$, suppose that
\begin{enumerate}[label=(\roman*)]
\item $V(H)\cap V(S_i)=W_i$;
\item $(S_i,W_i)$ is $d$-compressible with respect to $(J_i,f_i)$;
\item $\beta_i\colon J_i\rightarrow K$ is a homomorphism satisfying $\beta_i(f_i(w))=\psi(w)$ for all $w\in W_i$.
\end{enumerate}
Then there exists a $d$-compression $(J,f,K',\psi')$ of $(H\cup S_1 \cup \dots \cup S_t,U)$ such that $K\In K'$ and $\psi'{\restriction_{V(H)}}=\psi$.
\end{cor}\COMMENT{
\proof
Assume we have a attached $S_1,\dots,S_{i-1}$ such that $(J,f,K',\psi')$ is a $d$-compression of $(H\cup S_1 \cup \dots \cup S_{i-1},U)$, where $K\In K'$ and $\psi'{\restriction_{V(H)}}=\psi$. So $\beta_i\colon J_i \rightarrow K'$ is a homomorphism satisfying $\beta_i(f_i(w))=\psi(w)=\psi'(w)$ for all $w\in W_i$ since $W_i\In V(H)$. Moreover, $V(H\cup S_1 \cup \dots \cup S_{i-1})\cap V(S_i)=W_i\cup \bigcup_{j<i} (W_i\cap W_j)=W_i$. By the above Proposition, there exists a $d$-compression $(J,f,K'',\psi'')$ of $(H\cup S_1 \cup \dots \cup S_i,U)$ such that $K'\In K''$ and $\psi''{\restriction_{V(H\cup S_1 \cup \dots \cup S_{i-1})}}=\psi'$, so $K\In K'\In K''$ and $\psi''{\restriction_{V(H)}}=\psi$.
\endproof
}

\section{Transformers} \label{sec:transform}

Given two vertex-disjoint graphs $H$ and $H'$, a graph $T$ is called an \defn{$(H,H')_F$-transformer} if both $H\cup T$ and $H'\cup T$ have $F$-decompositions and $V(H)\cup V(H')\In V(T)$ is independent in $T$. 
Transformers in this sense were introduced in \cite{BKLO} as building blocks for absorbers. For two graphs $H$ and $H'$, write $H\rightsquigarrow H'$ if there exists an edge-bijective homomorphism from $H$ to $H'$. When constructing good $F$-absorbers, a crucial step is to have a good $(H,H')_F$-transformer whenever $H\rightsquigarrow H'$. In terms of how to build absorbers out of such transformers, the main ideas in \cite{BKLO} are essentially sufficient for our purposes. We will discuss this as briefly as possible in Section~\ref{sec:absorbers}. However, in order to achieve our goals, we need more sophisticated transformers. This will be our focus in this and the next section.

We call a graph $F$ \defn{$\delta$-transforming} if the following holds:
\begin{itemize}
\item[] Let $1/n \ll 1/k_0', \eps \ll \alpha, 1/b \ll 1/m, \mu,1/|F|$ and suppose that $G$ is a graph on $n$ vertices with $\delta(G) \ge (\delta + \mu) n$ which has an $(\alpha,\eps,k)$-partition for some $k\le k_0'$, and suppose that $H$ and $H'$ are vertex-disjoint subgraphs of $G$ of order at most $m$, where $H$ is $gcd(F)$-regular and $H\rightsquigarrow H'$. Then $G$ contains an $(H,H')_F$-transformer of order at most $b$.
\end{itemize}

In Section~\ref{sec:absorbers}, we will see that $F$ being $\delta$-transforming implies that $F$ is $\delta$-absorbing (see Lemma~\ref{lem:transf2abs}).
We will build transformers out of so-called switchers. Let $S$ be a graph and $U$ an independent set in $S$. Let $E_1,E_2$ be two disjoint sets of edges on~$U$. We call $S$ an \defn{$(E_1,E_2)_F$-switcher} if both $S\cup E_1$ and $S\cup E_2$ are $F$-decomposable. We will mostly use `cycle switchers' and `double-star switchers'. In the first case, $E_1$ and $E_2$ are the two perfect matchings forming an even cycle. In the second case, $E_1$ and $E_2$ are two stars with the same leaves, but distinct centers.\COMMENT{doesnt make sense if $r=1$, but I think we can ignore that here} 

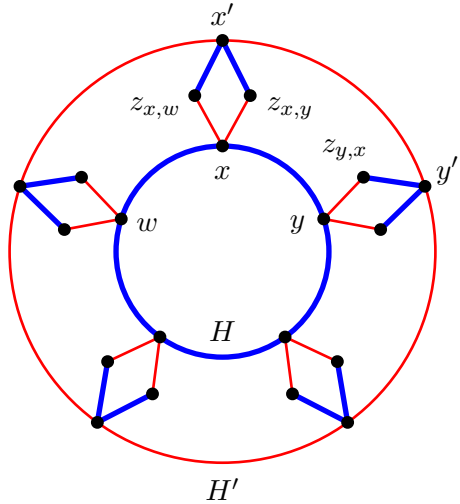
\begin{figure}[htbp]
\centering
\begin{tikzpicture}[scale=0.7, line width = 2pt]

	\draw[blue] (0,0) circle (2);
	\draw[red, line width = 1pt] (0,0) circle (4);
	
	\begin{scope}[rotate = 90]
		\foreach \i in {0,72,...,288}
		{
				\draw [red,line width = 1pt] (\i+10:3)--(\i:2)--(\i-10:3);
				\draw [blue] (\i+10:3)--(\i:4)--(\i-10:3);
				
				\filldraw[fill=black] (\i:2) circle (2pt);
				\filldraw[fill=black] (\i:4) circle (2pt);
				\filldraw[fill=black] (\i+10:3) circle (2pt);
				\filldraw[fill=black] (\i-10:3) circle (2pt);
		}
		
		\end{scope}
		
		\node  at (0,-1.5)  {$H$};
		\node  at (0,-4.5)  {$H'$};
		
		\node  at (0,1.5)  {$x$};
		\node  at (0,4.5)  {$x'$};
		
		\node  at (20:1.5)  {$y$};
		\node  at (20:4.5)  {$y'$};
		
		\node  at (160:1.5)  {$w$};
		
		\node at ($(115:3)$) {$z_{x,w}$};
		\node at ($(65:3)$) {$z_{x,y}$};
		
		\node at ($(40:3)$) {$z_{y,x}$};
		
\end{tikzpicture}

\caption{Sketch of a $(C_5,C_5)_F$-transformer built from five $C_6$-switchers and five $K_{2,2}$-switchers.}
\label{fig:C5transformer}
\end{figure}

We now briefly sketch how to build transformers out of these switchers (details are given in Lemma~\ref{lem:switch2transform}). For the sake of simplicity, suppose that $H'$ is a vertex-disjoint copy of $H$ and that $H$ is $r$-regular, where $r:=gcd(F)$. For any $x\in V(H)$, let $x'$ denote its copy in $H'$. We will build the desired $(H,H')_F$-transformer by introducing $r$ new vertices for every $x\in V(H)$ and joining them to $x$ and $x'$. We can then pair up these `middle' vertices with the $r$ neighbours of $x$ in $H$. A number of $C_6$-switchers will now allow us to translate the transforming task between $H$ and $H'$ into a switching task between two stars with $r$ common leaves. For example, let $xy\in E(H)$ and assume that $z_{x,y}$ is the middle vertex between $x$ and $x'$ associated with $y$. Similarly, assume that $z_{y,x}$ is the middle vertex between $y$ and $y'$ associated with $x$. Then let $E_{H}:=\Set{xy,z_{x,y}x',z_{y,x}y'}$ and $E_{H'}:=\Set{x'y',z_{x,y}x,z_{y,x}y}$. A $C_6$-switcher will now allow us to either cover $E_{H}$ or $E_{H'}$ with edge-disjoint copies of $F$. Doing this simultaneously for all edges of $H$, we can either cover (i) $E(H)$ together with all edges of the form $z_{x,y}x'$ or (ii) $E(H')$ together with all edges of the form $z_{x,y}x$ (see Figure~\ref{fig:C5transformer}).
For every $x\in V(H)$, a $K_{2,r}$-switcher will then enable us to cover all edges of the form $z_{x,y}x$ in case of~(i), or all edges of the form $z_{x,y}x'$ in case of~(ii).

In order to find the desired switchers in a graph $G$, we will use Lemma~\ref{lem:rooted embedding} and therefore have to equip switchers with compressions. Since compressions involve a fair amount of notation, we introduce the following conventions. When mentioning the cycle $C_\ell$, we usually assume that $V(C_\ell)=\Set{\mathfrak{c}_1,\dots,\mathfrak{c}_\ell}$ and $E(C_{\ell})=\set{\mathfrak{c}_i\mathfrak{c}_{i+1}}{i\in[\ell-1]}\cup\Set{\mathfrak{c}_1\mathfrak{c}_\ell}$. Similarly, $\Set{\mathfrak{p}_1,\dots,\mathfrak{p}_{\ell+1}}$ denotes the vertex set of the path $P_\ell$ with edge set $E(P_\ell)=\Set{\mathfrak{p}_1\mathfrak{p}_2,\dots,\mathfrak{p}_{\ell}\mathfrak{p}_{\ell+1}}$.

Very often, we want to switch between the two perfect matchings of an even cycle. For the sake of readability, $S$ is called a \defn{$(u_1,\dots,u_{2\ell})_F$-switcher} if $S$ is an $(E_1,E_2)_F$-switcher, where
\begin{align*}
E_1 := \Set{u_1u_2,u_3u_4,\dots,u_{2\ell-1}u_{2\ell}} \ \ \ \text{and} \ \ \
E_2 := \Set{u_2u_3,\dots,u_{2\ell-2}u_{2\ell-1},u_{2\ell}u_1}.
\end{align*}

\begin{figure}[htbp]

\begin{minipage}{0.49\textwidth} 

\centering
\begin{tikzpicture}[scale=0.55]

\fill[fill=black!30]
		(120 : 4.75 ) arc (120:240:4.75) -- (240 : 3.25 ) arc (240:120:3.25)--cycle;

\fill[fill=black!30]
		(240 : 4.75 ) arc (240:360:4.75) -- (360 : 3.25 ) arc (360:240:3.25)--cycle;
		
\fill[fill=black!30]
		(0 : 4.75 ) arc (0:120:4.75) -- (120 : 3.25 ) arc (120:0:3.25)--cycle;
	
		\foreach \i in {0,120,240}
		{
				\fill [fill=black!30,rotate=\i] (0,0.75) rectangle (4,-0.75);
				\draw [fill = white] (\i:4) circle (1);
		}
		\draw [fill = white] (0,0) circle (1);
		
		\filldraw[fill=black] (120:4) circle (2pt);
		\filldraw[fill=black] (240:4) circle (2pt);
		\filldraw[fill=black] (0,0.5) circle (2pt);
		\filldraw[fill=black] (0,0) circle (2pt);
		\filldraw[fill=black] (0,-0.5) circle (2pt);

		\begin{scope}[red,line width = 1pt]
		\draw (120:4) -- (0,0.5);
		\draw (120:4) -- (0,0);
		\draw (120:4) -- (0,-0.5);
		\end{scope}
		
		\begin{scope}[blue,line width = 1pt]
		\draw (240:4) -- (0,0.5);
		\draw (240:4) -- (0,0);
		\draw (240:4) -- (0,-0.5);
		\end{scope}
		
		\node  at (0.4,0.5)  {$u_1$};
		\node  at (0.4,0)  {$u_2$};
		\node  at (0.4,-0.5)  {$u_3$};
		\node  at (120:4.4)  {$u_4$};
		\node  at (240:4.4)  {$u_5$};
		\node  at (0:5.4)  {$\mathfrak{a}$};
		\node  at (120:5.4)  {$\mathfrak{p}_1$};
		\node  at (0:-1.4)  {$\mathfrak{p}_2$};
		\node  at (240:5.4)  {$\mathfrak{p}_3$};

\end{tikzpicture}

\end{minipage}
\hfill
\begin{minipage}{0.49\textwidth}

\centering
\begin{tikzpicture}[scale=0.5]
		
				\fill [fill=black!30] (-4,4.75) rectangle (4,3.25);
				\fill [fill=black!30] (-4,-4.75) rectangle (4,-3.25);
				\fill [fill=black!30] (4.75,4) rectangle (3.25,-4);
				\fill [fill=black!30] (-4.75,4) rectangle (-3.25,-4);
				
				\foreach \i in {45,135,225,315}
		{
				\fill [fill=black!30,rotate=\i] (0,0.75) rectangle (5.6,-0.75);
		}
		
				\draw [fill = white] (0,0) circle (1);
				\draw [fill = white] (4,4) circle (1);
				\draw [fill = white] (-4,4) circle (1);
				\draw [fill = white] (4,-4) circle (1);
				\draw [fill = white] (-4,-4) circle (1);
				
				\filldraw[fill=black] (4,4) circle (2pt);
				\filldraw[fill=black] (-4,4) circle (2pt);
				\filldraw[fill=black] (4,-4) circle (2pt);
				\filldraw[fill=black] (-4,-4) circle (2pt);
				
		\begin{scope}[red,line width = 1pt]
		\draw (4,4) -- (-4,4);
		\draw (4,-4) -- (-4,-4);
		\end{scope}
		
		\begin{scope}[blue,line width = 1pt]
		\draw (4,4) -- (4,-4);
		\draw (-4,4) -- (-4,-4);
		\end{scope}
			
		\node  at (-4.5,4)  {$u_1$};
		\node  at (4.5,4)  {$u_2$};
		\node  at (4.5,-4)  {$u_3$};
		\node  at (-4.5,-4)  {$u_4$};
		
		\node  at (0:1.4)  {$\mathfrak{a}$};
		\node  at (-5,5)  {$\mathfrak{c}_1$};
		\node  at (5,5)  {$\mathfrak{c}_2$};
		\node  at (5,-5)  {$\mathfrak{c}_3$};
		\node  at (-5,-5)  {$\mathfrak{c}_4$};	
\end{tikzpicture}

\end{minipage}
\caption{A $3$-compressible $(K_{2,3})_F$-switcher with augmentation $\{\mathfrak{p}_1\mathfrak{p}_3\}$ and a $4$-compressible $(C_{4})_F$-switcher with augmentation $\emptyset$. The shaded areas indicate where the edges of the switchers may lie.}
\label{fig:switchers}
\end{figure}
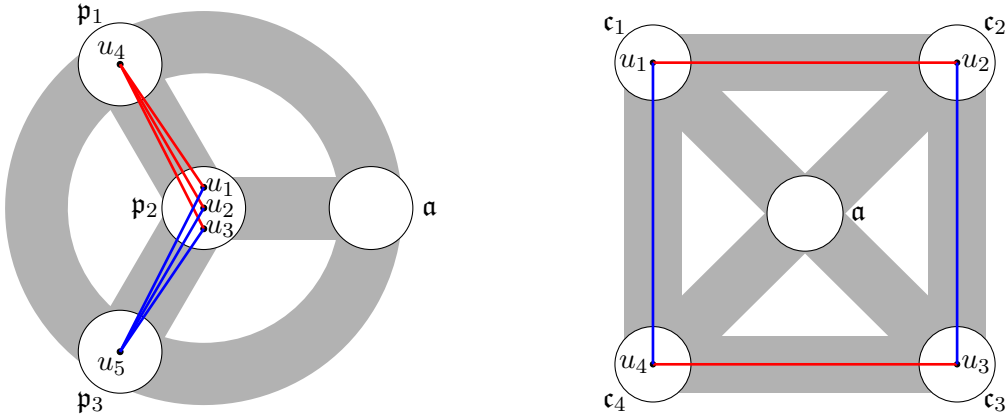

Let $\ell\ge 2$, $d\ge 0$ and suppose that $Aug \In E(\overline{C_{2\ell}})$.
A \defn{$d$-compressible $(C_{2\ell})_F$-switcher with augmentation $Aug$} is a model $(S,\Set{u_1,\dots,u_{2\ell}})$ satisfying the following properties (see Figure~\ref{fig:switchers}):
\begin{itemize}
\item it is $d$-compressible with respect to the root-compression $(C_{2\ell} \cup Aug,f)$, where
\item $f(u_i):=\mathfrak{c}_i$ for all $i\in[2\ell]$;
\item $S$ is a $(u_1,\dots,u_{2\ell})_F$-switcher.
\end{itemize}
The set $Aug$ may be viewed as unwanted, and we will omit saying `with augmentation~$\emptyset$'. When using switchers to build transformers, we must in fact have $Aug=\emptyset$ (see Lemma~\ref{lem:switch2transform}). However, when constructing switchers in Section~\ref{sec:switchers}, we will first obtain compressions where $Aug\neq\emptyset$ and then perform reductions to achieve that $Aug=\emptyset$.

Similarly, let $r\ge 1$, $d\ge 0$ and suppose that $Aug \In E(\overline{P_2})$. A \defn{$d$-compressible $(K_{2,r})_F$-switcher with augmentation $Aug$} is a model $(S,\Set{u_1,\dots,u_{r+2}})$ satisfying the following properties (see Figure~\ref{fig:switchers}):
\begin{itemize}
\item it is $d$-compressible with respect to the root-compression $(P_2\cup Aug,f)$, where
\item $f(u_i):=\mathfrak{p}_2$ for all $i\in[r]$, $f(u_{r+1}):=\mathfrak{p}_1$ and $f(u_{r+2}):=\mathfrak{p}_3$;
\item $S$ is a $(\set{u_{r+1}u_i}{i\in [r]},\set{u_{r+2}u_i}{i\in [r]})_F$-switcher.
\end{itemize}
Note that the existence of a $(K_{2,r})_F$-switcher implies that $gcd(F)\mid r$. Also note that though as graphs $C_4$ and $K_{2,2}$ are isomorphic, a $(C_{4})_F$-switcher and a $(K_{2,2})_F$-switcher are two different concepts according to the above definitions.

\begin{lemma} \label{lem:switch2transform}
Let $F$ be any graph and $r:=gcd(F)$. Let $d\ge 2$ and assume that there exists a $d$-compressible $(C_6)_F$-switcher and a $d$-compressible $(K_{2,r})_F$-switcher. Then $F$ is $(1-1/d)$-transforming.
\end{lemma}

\proof
Let $U^\star=\Set{u^\star_1,\dots,u^\star_{r+2}}$ and $U^\circ=\Set{u^\circ_1,\dots,u^\circ_6}$ be sets. Define the function $f^\star\colon U^\star \rightarrow V(P_2)$ by $f^\star(u^\star_{r+1}):=\mathfrak{p}_1$, $f^\star(u^\star_{r+2}):=\mathfrak{p}_3$ and $f^\star(u^\star_i):=\mathfrak{p}_2$ for all $i\in[r]$. Moreover, let $f^\circ\colon U^\circ \rightarrow V(C_6)$ be defined as $f^\circ(u^\circ_i):=\mathfrak{c}_i$ for all $i\in[6]$.

By our assumption, there exists a $(\set{u^\star_{r+1}u^\star_i}{i\in [r]},\set{u^\star_{r+2}u^\star_i}{i\in [r]})_F$-switcher $S^\star$ such that $(S^\star,U^\star)$ is a $d$-compressible model with respect to $(P_2,f^\star)$.

Moreover, there exists a $(u^\circ_1,u^\circ_2,u^\circ_3,u^\circ_4,u^\circ_5,u^\circ_6)_F$-switcher $S^\circ$ such that $(S^\circ,U^\circ)$ is a $d$-compressible model with respect to $(C_6,f^\circ)$.

Let $t:=\max\Set{|S^\star|,|S^\circ|}$.
Let $1/n \ll 1/k_0',\eps \ll \alpha, 1/b \ll 1/m,\mu,1/|F|$. So we may assume $\alpha, 1/b \ll 1/t$ since $t$ only depends on $F$.\COMMENT{and among all possible $S$'s, we can choose those which minimise $t$.} Suppose that $G$ is any graph on $n$ vertices with $\delta(G) \ge (1-1/d + \mu) n$ which has an $(\alpha,\eps,k)$-partition $V_1,\dots,V_k$ for some $k\le k_0'$, and suppose that $H$ and $H'$ are vertex-disjoint subgraphs of $G$ of order at most $m$, where $H$ is $r$-regular and $\phi$ is an edge-bijective homomorphism from $H$ to $H'$. For a vertex $x\in V(H)$, to enhance readability, we will sometimes write $x'$ for $\phi(x)$.
We need to show that $G$ contains an $(H,H')_F$-transformer of order at most $b$.

\medskip
\noindent\emph{Step 1: Setup}
\medskip

Let $R$ be the reduced graph of $V_1,\dots,V_k$ and $\sigma\colon V(G)\rightarrow V(R)$ the associated cluster function. Consider any vertex $x\in V(H)$. Note that $|N_G(x)\cap N_G(x')|\ge 2\mu n$.\COMMENT{Only place where we need $d\ge 2$.} Hence, there exists an index $\tau(x)\in[k]$ such that $$|N_G(x)\cap N_G(x') \cap V_{\tau(x)}| \ge 2\mu |V_{\tau(x)}|.$$ Let $\hat{V}_{\tau(x)}:=N_G(x)\cap N_G(x') \cap V_{\tau(x)}$.
Let $y$ be a neighbour of $x$ in $H$. Since $xy,x'y'\in E(G)$, the pairs $(V_{\sigma(x)},V_{\sigma(y)})$ and $(V_{\sigma(x')},V_{\sigma(y')})$ must be weakly-$(\alpha,\eps)$-super-regular. Moreover, $D_{x,y}:=N_G(y)\cap V_{\sigma(x)}$ has size at least $\alpha|V_{\sigma(x)}|$ and $D'_{x,y}:=N_G(y')\cap V_{\sigma(x')}$ has size at least $\alpha|V_{\sigma(x')}|$.

Our goal in this step is to find for every $x\in V(H)$ a set $\Lambda_x\In V(G)$ with the following properties:
\begin{enumerate}[label=(L\arabic*)]
\item $\Lambda_x \In \hat{V}_{\tau(x)}$;
\item $|\Lambda_x| \ge \mu|V_{\tau(x)}|$;
\item for all $z\in \Lambda_x$ and $y\in N_H(x)$, $d_G(z,D_{x,y}) \ge \alpha^2 |V_{\sigma(x)}|/2$ and $d_G(z,D'_{x,y}) \ge \alpha^2 |V_{\sigma(x')}|/2$.
\end{enumerate}

Consider $x\in V(H)$. By the definition of $\tau(x)$, the pairs $(V_{\sigma(x)},V_{\tau(x)})$ and $(V_{\sigma(x')},V_{\tau(x)})$ are weakly-$(\alpha,\eps)$-super-regular. Since $|\hat{V}_{\tau(x)}|\ge 2\mu |V_{\tau(x)}|$, we deduce that the pairs $(V_{\sigma(x)},\hat{V}_{\tau(x)})$ and $(V_{\sigma(x')},\hat{V}_{\tau(x)})$ are $\eps/\mu$-regular with density at least $\alpha-\eps$.
By Fact~\ref{fact:regularity},\COMMENT{and since $\alpha\ge \eps/\mu$} we know that there is a set $\Lambda_x\In \hat{V}_{\tau(x)}$ containing all but at most $2r\eps|\hat{V}_{\tau(x)}|/\mu$ vertices of $\hat{V}_{\tau(x)}$ such that every $z\in \Lambda_x$ has at least $(\alpha-\eps-\eps/\mu)|D_{x,y}|$ neighbours in $D_{x,y}$ and at least $(\alpha-\eps-\eps/\mu)|D'_{x,y}|$ neighbours in $D'_{x,y}$ for every $y\in N_H(x)$. Checking that $|\Lambda_x| \ge (1-2r\eps/\mu)|\hat{V}_{\tau(x)}| \ge \mu|V_{\tau(x)}|$ and $d_G(z,D_{x,y}) \ge (\alpha-\eps-\eps/\mu)|D_{x,y}| \ge \alpha^2 |V_{\sigma(x)}|/2$ and equally $d_G(z,D'_{x,y}) \ge \alpha^2 |V_{\sigma(x')}|/2$ for all $z\in \Lambda_x$ and $y\in N_H(x)$ confirms that $\Lambda_x$ satisfies (L1)--(L3).

\medskip
\noindent\emph{Step 2: Finding switchers}
\medskip

We first find the double-star switchers. This means that for every $x\in V(H)$, we want to find a subgraph $S^\star_x$ of $G$ and a set $Z_x\In V(G)$ such that
\begin{enumerate}[label=(S\arabic*$^\star$)]
\item $|Z_x|=r$ and $Z_x\In \Lambda_x$;\label{middle vertices}
\item $(S^\star_x)_{x\in V(H)}$, $(E^+_x \cup E^-_x)_{x\in V(H)}$ and $G[V(H)\cup V(H')]$ are edge-disjoint, where $E_x^+:=\set{xz}{z\in Z_x}$ and $E_x^-:=\set{x'z}{z\in Z_x}$; \label{item:star-switching}
\item $|S^\star_x|\le t$;
\item $S^\star_x$ is an $(E^+_x,E^-_x)_F$-switcher.\label{star-switch}
\end{enumerate}

We can find them one by one using Lemma~\ref{lem:rooted embedding}. Suppose that for some subset $Y\In V(H)$, we have already found $S^\star_y$ and $Z_y$ for all $y\in Y$ and now want to define $S^\star_x$ and $Z_x$ for $x\in V(H)\sm Y$. Let $G_0$ be the subgraph consisting of all edges of $(S^\star_y)_{y\in Y}$, $(E^+_y \cup E^-_y)_{y\in Y}$ and $G[V(H)\cup V(H')]$. Let $G':=G-G_0$. Note that since $\Delta(G_0) \le tm+rm+2m$, $V_1,\dots,V_k$ is an $(\alpha/2,3\eps,k)$-partition for $G'$.\COMMENT{see end of Section~\ref{sec:regularity}}

Define a $G'$-labelling $\Lambda$ of $U^\star$ as follows: $\Lambda(u^\star_{r+1}):=\Set{x}$,  $\Lambda(u^\star_{r+2}):=\Set{x'}$,  $\Lambda(u^\star_i):=\Lambda_x$ for all $i\in[r]$. It remains to check that $\Lambda$ is $(\alpha/2,3\eps,k)$-admissible. Clearly, $\Lambda$ respects $(P_2,f^\star)$. 
Let $j(\mathfrak{p}_1):=\sigma(x)$, $j(\mathfrak{p}_2):=\tau(x)$ and $j(\mathfrak{p}_3):=\sigma(x')$. Then \ref{adm:hom} holds. \ref{adm:label size} holds because $|\Lambda_{P_2}(\mathfrak{p}_2)|=|\Lambda_x|\ge \mu|V_{\tau(x)}| \ge \alpha |V_{j(\mathfrak{p}_2)}|/2$ by (L2), where $\Lambda_{P_2}$ is as in Definition~\ref{def:admissible}. For \ref{adm:conform}, note that $W_{\mathfrak{p}_1}=W_{\mathfrak{p}_3}=\emptyset$, where $W_{\mathfrak{p}_i}$ is as in Definition~\ref{def:admissible}, so \ref{adm:conform} holds trivially in this case. Moreover, $N_G(W_{\mathfrak{p}_2},V_{j(\mathfrak{p}_2)})=N_G(\Set{x,x'},V_{\tau(x)})=\hat{V}_{\tau(x)}$ and so $d_G(W_{\mathfrak{p}_2},V_{j(\mathfrak{p}_2)})\ge \mu |V_{j(\mathfrak{p}_2)}|$. \ref{adm:label loc} holds by (L1).

Hence, by Lemma~\ref{lem:rooted embedding}, there exists an embedding $\rho$ of $(S^\star,U^\star)$ into $G'$ respecting $\Lambda$. Let $S^\star_x:=\rho(S^\star)$ and $Z_x:=\rho(\set{u^\star_i}{i\in[r]})$. Then $S^\star_x$ and $Z_x$ satisfy (S1$^\star$)--(S4$^\star$).

We now find the cycle switchers. For this, we associate the vertices of $Z_x$ with the neighbours of $x$ in $H$, that is, we assume $Z_x=\set{z_{x,y}}{y\in N_H(x)}$.\COMMENT{$|N_H(x)|=|Z_x|$, so take any bijection and define $z_{x,y}$ as the image of $y$.} For every edge $xy\in E(H)$, we want to find a subgraph $S^\circ_{xy}$ of $G$ such that
\begin{enumerate}[label=(S\arabic*$^\circ$)]
\item $(S^\circ_{xy})_{xy\in V(H)}$, $(S^\star_x)_{x\in V(H)}$, $(E^+_x \cup E^-_x)_{x\in V(H)}$ and $G[V(H)\cup V(H')]$ are edge-disjoint;
\item $|S^\circ_{xy}|\le t$;
\item $S^\circ_{xy}$ is an $(x,y,z_{y,x},y',x',z_{x,y})$-switcher.\label{cycle-switch}
\end{enumerate}

Again, we find them one by one using Lemma~\ref{lem:rooted embedding}. Suppose that for some subset $Y\In E(H)$, we have already found $S^\star_{e}$ for all $e\in Y$, and that we now want to define $S^\star_{xy}$ for $xy\in E(H)\sm Y$. Let $G_0$ be the subgraph consisting of all edges of $(S^\circ_e)_{e\in Y}$, $(S^\star_x)_{x\in V(H)}$, $(E^+_x \cup E^-_x)_{x\in V(H)}$ and $G[V(H)\cup V(H')]$ and let $G':=G-G_0$. Since $\Delta(G_0) \le trm/2 + tm + rm + 2m$, we have that $V_1,\dots,V_k$ is an $(\alpha/2,3\eps,k)$-partition for $G'$.

Define a $G'$-labelling $\Lambda$ of $U^\circ$ as follows: $\Lambda(u^\circ_1):=\Set{x}$, $\Lambda(u^\circ_2):=\Set{y}$, $\Lambda(u^\circ_3):=\Set{z_{y,x}}$, $\Lambda(u^\circ_4):=\Set{y'}$, $\Lambda(u^\circ_5):=\Set{x'}$, $\Lambda(u^\circ_6):=\Set{z_{x,y}}$.
Trivially, $\Lambda$ respects $(C_6,f^\circ)$. 
Define $j(\mathfrak{c}_1):=\sigma(x)$, $j(\mathfrak{c}_2):=\sigma(y)$, $j(\mathfrak{c}_3):=\tau(y)$, $j(\mathfrak{c}_4):=\sigma(y')$, $j(\mathfrak{c}_5):=\sigma(x')$ and $j(\mathfrak{c}_6):=\tau(x)$. Then \ref{adm:hom} holds. Moreover, \ref{adm:label size} and \ref{adm:label loc} hold trivially since $|\Lambda(u^\circ_i)|=1$ for all $i\in[6]$. We will now check that \ref{adm:conform} holds. Note that $W_{\mathfrak{c}_1}=\Set{z_{x,y},y}$ and so $$d_G(W_{\mathfrak{c}_1},V_{j(\mathfrak{c}_1)})=d_G(\Set{z_{x,y},y},V_{\sigma(x)}) = d_G(z_{x,y},D_{x,y}) \ge\alpha^2 |V_{\sigma(x)}|/2$$ by \ref{middle vertices} and (L3). The same applies to $\mathfrak{c}_i$ with $i\in \Set{2,4,5}$. Note that $W_{\mathfrak{c}_3}=\Set{y,y'}$ and so $$d_G(W_{\mathfrak{c}_1},V_{j(\mathfrak{c}_3)})=d_G(\Set{y,y'},V_{\tau(y)}) = |\hat{V}_{\tau(y)}| \ge 2\mu |V_{\tau(y)}|$$ by definition of $\tau(y)$. The same applies to $\mathfrak{c}_6$. Thus, $\Lambda$ is $(\alpha^2/2,3\eps,k)$-admissible.

Hence, by Lemma~\ref{lem:rooted embedding}, there exists an embedding $\rho$ of $(S^\circ,U^\circ)$ into $G'$ respecting $\Lambda$. Let $S^\circ_{xy}:=\rho(S^\circ)$. Then (S1$^\circ$)--(S3$^\circ$) are satisfied.

\medskip
\noindent\emph{Step 3: Transforming}
\medskip

We can now define the desired $(H,H')_F$-transformer.
Observe that \ref{item:star-switching} implies that
\begin{align}
\bigcup_{x\in V(H)}E_x^+ &= \bigcup_{xy\in E(H)}\Set{xz_{x,y},yz_{y,x}}, \label{eqn:star edges 1}\\
\bigcup_{x\in V(H)}E_x^- &= \bigcup_{xy\in E(H)}\Set{x'z_{x,y},y'z_{y,x}}. \label{eqn:star edges 2}
\end{align}

Let $$T:=\bigcup_{xy\in E(H)}S^\circ_{xy} \cup \bigcup_{x\in V(H)}(S^\star_x \cup E^+_x \cup E^-_x).$$ We claim that $T$ is the desired $(H,H')_F$-transformer. By construction, $T$ has order at most $trm/2+mt\le b$ and $T[V(H)\cup V(H')]$ is empty.
Finally and most importantly, 
$$T\cup H \overset{\eqref{eqn:star edges 2}}{=}\bigcup_{xy\in E(H)}(S^\circ_{xy} \cup \Set{xy,x'z_{x,y},y'z_{y,x}}) \cup \bigcup_{x\in V(H)}(S^\star_x \cup E^+_x)$$ is $F$-decomposable by \ref{cycle-switch} and \ref{star-switch}. Similarly,
$$T\cup H' \overset{\eqref{eqn:star edges 1}}{=}\bigcup_{xy\in E(H)}(S^\circ_{xy} \cup \Set{x'y',xz_{x,y},yz_{y,x}}) \cup \bigcup_{x\in V(H)}(S^\star_x \cup E^-_x)$$ is $F$-decomposable.
\endproof

\section{Constructing switchers} \label{sec:switchers}

In the previous section, we saw how the problem of finding an $(H,H')_F$-transformer for two rather arbitrary graphs $H,H'$ can be reduced to the problem of constructing well-compressible switchers. We will now describe such switchers. The following `discretisation lemma' is key to narrowing the value of $\delta_F$ down for Theorem~\ref{thm:main}(ii)--(iii). 

The main idea of the lemma is as follows. Suppose that $d\in \bN$ and $\delta_F=1-1/d-\mu$ for some $\mu>0$ and we aim to construct a $(d-1)$-compressible $(E_1,E_2)_F$-switcher $S$. We therefore have to show that both $S\cup E_1$ and $S\cup E_2$ are $F$-decomposable. We will achieve this by simply using the definition of $\delta_F$. More precisely, we will consider an arbitrary large graph $S$ such that both $S\cup E_1$ and $S\cup E_2$ are $F$-divisible and $\delta(S)\ge (\delta_F+\mu/2)|S|$. Then, both $S\cup E_1$ and $S\cup E_2$ must be $F$-decomposable. Moreover, we can also ensure that $S$ is $d$-partite. The location of a vertex in one of the $d$ classes then naturally defines a homomorphism $\psi\colon S\to K_d$, which will ensure that the switcher $S$ is $(d-1)$-compressible.

\begin{lemma}[Discretisation lemma] \label{lem:discretisation}
Let $F$ be any graph and let $r:=gcd(F)$ and $d\in \bN$. Suppose that $\delta_F < 1-1/d$. Then the following assertions hold.
\begin{enumerate}[label=(\roman*)]
\item There exists a $(d-1)$-compressible $(C_4)_F$-switcher with augmentation $\Set{\mathfrak{c}_1\mathfrak{c}_3,\mathfrak{c}_2\mathfrak{c}_4}$.
\item There exists a $(d-1)$-compressible $(K_{2,r})_F$-switcher with augmentation $\Set{\mathfrak{p}_1\mathfrak{p}_3}$.
\end{enumerate}
\end{lemma}

\proof
Let $\mu:=1-1/d-\delta_F$, $d^\circ:=\max\Set{d,4}$ and $d^\star:=\max\Set{d,3}$. By definition of $\delta_F$, there exists $n_0\in \bN$ such that whenever $G$ is an $F$-divisible graph on $n\ge n_0$ vertices with $\delta(G) \ge (\delta_F+\mu/2)n$, then $G$ is $F$-decomposable. Let $s\in \bN$ be such that $1/s\ll 1/n_0,1/|F|,1/d,\mu$ and such that $s$ is divisible by $e(F)$ and $gcd(F)$.\COMMENT{$s\ge \max\Set{n_0/d^\circ,2(r+1)|F|/\mu,(r+1)|F|+1}$ is enough}

To prove (i), let $K^\circ$ be the complete graph on $\Set{\mathfrak{c}_1,\dots,\mathfrak{c}_{d^\circ}}$. Clearly, $K^\circ$ has degeneracy at most $d-1$ rooted at $\Set{\mathfrak{c}_1,\dots,\mathfrak{c}_4}$.
Let $G^\circ$ be a complete $d^\circ$-partite graph with vertex classes $V_1,\dots,V_{d^\circ}$ of size $s$ each. Then, $G^\circ$ is $F$-divisible and $|G^\circ|\ge n_0$.\COMMENT{$e(G)=s^2\binom{d^\circ}{2}\equiv 0 \mod{e(F)}$ and $d_G(x)=(d^\circ-1)s\equiv 0\mod{gcd(F)}$ for all $x\in V(G)$.} For all $i\in[4]$, pick some $u_i\in V_i$ and define $f^\circ(u_i):=\mathfrak{c}_i$. Note that $\chi(F)\le d \le d^\circ$ since trivially $\delta_F\ge 1-1/(\chi(F)-1)$. Since $s$ is sufficiently large, we can easily find edge-disjoint copies $F_1,\dots,F_4$ of $F$ in $G^\circ$ such that $u_1u_3\in E(F_1)$, $u_1u_4\in E(F_2)$, $u_2u_3\in E(F_3)$ and $u_2u_4\in E(F_4)$, but $u_1u_2,u_3u_4\notin E(F_1\cup\dots\cup F_4)$. Define $S^\circ:=G^\circ-\Set{u_1u_2,u_3u_4}-(F_1\cup\dots\cup F_4)$. Let $\psi^\circ\colon S^\circ\rightarrow K^\circ$ be such that $\psi^\circ(x)=\mathfrak{c}_i$ if and only if $x\in V_i$. Hence, $\psi^\circ$ is a homomorphism. Moreover, $\Set{u_1,\dots,u_4}$ is independent in $S^\circ$, so $(K^\circ[\Set{\mathfrak{c}_1,\dots,\mathfrak{c}_4}],f^\circ,K^\circ,\psi^\circ)$ is a $(d-1)$-compression of $(S^\circ,\Set{u_1,\dots,u_4})$.
It remains to show that $S^\circ$ is a $(u_1,u_2,u_3,u_4)_F$-switcher. But $S^\circ\cup\Set{u_1u_2,u_3u_4}=G^\circ-(F_1\cup\dots\cup F_4)$ is $F$-divisible and thus $S^\circ\cup\Set{u_2u_3,u_4u_1}$ is $F$-divisible as well. Finally, $$\delta(S^\circ\cup\Set{u_1u_2,u_3u_4}) \ge (d^\circ-1)s -4|F| \ge (1-1/d^\circ-\mu/2)sd^\circ\ge (\delta_F+\mu/2)|S^\circ|$$ and $\delta(S^\circ\cup\Set{u_2u_3,u_4u_1})=\delta(S^\circ\cup\Set{u_1u_2,u_3u_4})$. So since $|S^\circ|=|G^\circ|\ge n_0$, both $S^\circ\cup\Set{u_1u_2,u_3u_4}$ and $S^\circ\cup\Set{u_2u_3,u_4u_1}$ are $F$-decomposable.

To prove (ii), let $K^\star$ be the complete graph on $\Set{\mathfrak{p}_1,\dots,\mathfrak{p}_{d^\star}}$. Clearly, $K^\star$ has degeneracy at most $d-1$ rooted at $\Set{\mathfrak{p}_1,\mathfrak{p}_2,\mathfrak{p}_3}$.
Let $G^\star$ be a complete $d^\star$-partite graph with vertex classes $V_1,\dots,V_{d^\star}$ of size $s$ each. Then, $G^\star$ is $F$-divisible and $|G^\star|\ge n_0$.\COMMENT{$e(G)=s^2\binom{d^\star}{2}\equiv 0 \mod{e(F)}$ and $d_G(x)=(d^\star-1)s\equiv 0\mod{gcd(F)}$ for all $x\in V(G)$.} Let $u_{r+1}\in V_1$, $u_1,\dots,u_r\in V_2$ and $u_{r+2}\in V_3$. Define $f^\star(u_i):=\mathfrak{p}_2$ for all $i\in[r]$, $f^\star(u_{r+1}):=\mathfrak{p}_1$, and $f^\star(u_{r+2}):=\mathfrak{p}_3$. Let $E^+:=\set{u_{r+1}u_i}{i\in[r]}$ and $E^-:=\set{u_{r+2}u_i}{i\in[r]}$.
Let $F_1,\dots,F_{r+1}$ be edge-disjoint copies of $F$ in $G^\star$ such that $u_{r+2}u_i\in E(F_i)$ for all $i\in[r+1]$ and $E^+ \cap E(F_1\cup\dots\cup F_{r+1})=\emptyset$.
Define $S^\star:=G^\star-E^+-(F_1\cup\dots\cup F_{r+1})$. Let $\psi^\star\colon S^\star\rightarrow K^\star$ be such that $\psi^\star(x)=\mathfrak{p}_i$ if and only if $x\in V_i$. Hence, $\psi^\star$ is a homomorphism. Moreover, $\Set{u_1,\dots,u_{r+2}}$ is independent in $S^\star$, so $(K^\star[\Set{\mathfrak{p}_1,\mathfrak{p}_2,\mathfrak{p}_3}],f^\star,K^\star,\psi^\star)$ is a $(d-1)$-compression of $(S^\star,\Set{u_1,\dots,u_{r+2}})$.
It remains to show that $S^\star$ is an $(E^+,E^-)_F$-switcher. But $S^\star\cup E^+=G^\star-(F_1\cup\dots\cup F_{r+1})$ is $F$-divisible and thus $S^\star\cup E^-$ is $F$-divisible as well.
Finally, $$\delta(S^\star\cup E^+) \ge (d^\star-1)s -(r+1)|F| \ge (1-1/d^\star-\mu/2)sd^\star\ge (\delta_F+\mu/2)|S^\star|,$$ and similarly, $\delta(S^\star\cup E^-)\ge (\delta_F+\mu/2)|S^\star|$. Therefore, both $S^\star\cup E^+$ and $S^\star\cup E^-$ are $F$-decomposable.
\endproof

Recall that in order to apply Lemma~\ref{lem:switch2transform}, we require a $(C_6)_F$-switcher and a $(K_{2,gcd(F)})_F$-switcher with no augmentations, whilst the above lemma outputs switchers with augmentations. In the following, we will carry out a sequence of reductions which will provide us with the switchers required for Lemma~\ref{lem:switch2transform}. Roughly speaking, in each reduction, we assume that we have a $d$-compressible switcher $S'$ with some augmentation(s). We then construct a $d$-compressible switcher $S$ by combining several copies of $S'$ in such a way that $S$ has fewer augmentations than $S'$. In order to ensure that $S$ is still $d$-compressible, we will appeal to Corollary~\ref{cor:model attaching}. Under rather natural assumptions, it allows us to attach $d$-compressible models to an existing model without increasing degeneracy.

\begin{lemma}\label{lem:star-switcher reduction}
Let $F$ be any graph, $d\ge 3$ and $r\in \bN$. Assume that there exists a $d$-compressible $(K_{2,r})_F$-switcher with augmentation $\Set{\mathfrak{p}_1\mathfrak{p}_3}$. Then there also exists a $d$-compressible $(K_{2,r})_F$-switcher (with augmentation $\emptyset$).
\end{lemma}

We prove the lemma as follows. First we add a new vertex $\mathfrak{a}$ connected to $\mathfrak{p}_1$, $\mathfrak{p}_2$ and $\mathfrak{p}_3$. We then obtain the desired switcher by combining two switchers with the underlying augmented paths $\mathfrak{p}_1\mathfrak{p}_2\mathfrak{a}$ and $\mathfrak{a}\mathfrak{p}_2\mathfrak{p}_3$.

\proof
Let $U=\Set{u_1,\dots,u_{r+2}}$ and let $S'$ be the graph on $U\cup\Set{w}$ with edge set $E^w:=\set{wu_i}{i\in[r]}$. Define $f\colon U \rightarrow V(P_2)$ as $f(u_i):=\mathfrak{p}_2$ for all $i\in[r]$, $f(u_{r+1}):=\mathfrak{p}_1$ and $f(u_{r+2}):=\mathfrak{p}_3$. Moreover, define $\psi$ such that $\psi{\restriction_U}=f$ and $\psi(w):=\mathfrak{a}$, where $\mathfrak{a}$ is a new vertex.
Let $K$ be the graph with $V(K)=V(P_2)\cup\Set{\mathfrak{a}}$ and $E(K)=E(P_2)\cup \Set{\mathfrak{a}\mathfrak{p}_1,\mathfrak{a}\mathfrak{p}_2,\mathfrak{a}\mathfrak{p}_3}$.
Clearly, $(P_2,f,K,\psi)$ is a $3$-compression of $(S',U)$. Set $E^+:=\set{u_{r+1}u_i}{i\in[r]}$ and $E^-:=\set{u_{r+2}u_i}{i\in[r]}$. 

Define $f^+\colon (U\sm \Set{u_{r+2}})\cup\Set{w}\rightarrow V(P_2)$ as $f^+(u_i):=\mathfrak{p}_2$ for all $i\in[r]$, $f^+(u_{r+1}):=\mathfrak{p}_1$ and $f^+(w):=\mathfrak{p}_3$.
Since there exists a $d$-compressible $(K_{2,r})_F$-switcher with augmentation $\Set{\mathfrak{p}_1\mathfrak{p}_3}$, there exists an $(E^+,E^w)_F$-switcher $S^+$ such that $(S^+,(U\sm \Set{u_{r+2}})\cup\Set{w})$ is $d$-compressible with respect to $(P_2\cup\Set{\mathfrak{p}_1\mathfrak{p}_3},f^+)$.

Similarly, there exists an $(E^-,E^w)_F$-switcher $S^-$ such that $(S^-,(U\sm \Set{u_{r+1}})\cup\Set{w})$ is $d$-compressible with respect to $(P_2\cup\Set{\mathfrak{p}_1\mathfrak{p}_3},f^-)$, where $f^-(u_i):=\mathfrak{p}_2$ for all $i\in[r]$, $f^-(u_{r+2}):=\mathfrak{p}_1$ and $f^-(w):=\mathfrak{p}_3$. We can also assume that $V(S^+)\cap V(S') = (U\sm \Set{u_{r+2}})\cup\Set{w}$, $V(S^-)\cap V(S') = (U\sm \Set{u_{r+1}})\cup\Set{w}$ and $V(S^+)\cap V(S^-) = \Set{w,u_1,\dots,u_r}$.

Now define $\beta^+\colon P_2\cup\Set{\mathfrak{p}_1\mathfrak{p}_3} \rightarrow K$ as $\beta^+(\mathfrak{p}_1):=\mathfrak{p}_1$, $\beta^+(\mathfrak{p}_2):=\mathfrak{p}_2$ and $\beta^+(\mathfrak{p}_3):=\mathfrak{a}$. Analogously, $\beta^-\colon P_2\cup\Set{\mathfrak{p}_1\mathfrak{p}_3} \rightarrow K$ is defined as $\beta^-(\mathfrak{p}_1):=\mathfrak{p}_3$, $\beta^-(\mathfrak{p}_2):=\mathfrak{p}_2$ and $\beta^-(\mathfrak{p}_3):=\mathfrak{a}$. Then, $\beta^+(f^+(v))=\psi(v)$ for all $v\in (U\sm \Set{u_{r+2}})\cup\Set{w}$ and $\beta^-(f^-(v))=\psi(v)$ for all $v\in (U\sm \Set{u_{r+1}})\cup\Set{w}$. Let $S:=S'\cup S^+ \cup S^-$. Hence, by Corollary~\ref{cor:model attaching}, $(S,U)$ is $d$-compressible with respect to $(P_2,f)$. Finally, since $S\cup E^+=(S^-\cup E^w)\cup(S^+\cup E^+)$ and $S\cup E^-=(S^+\cup E^w)\cup(S^-\cup E^-)$, $(S,U)$ is an $(E^+,E^-)_F$-switcher.
\endproof

Since the definitions of the homomorphisms $\beta$ and functions $f$ are usually natural and clear from the context, we will often omit the corresponding details in future applications of Corollary~\ref{cor:model attaching}. 
We can now combine several $(C_4)_F$-switchers into a  $(C_6)_F$-switcher.

\begin{lemma} \label{lem:C4-to-C6 reduction}
Let $F$ be any graph and $d\ge 3$. Assume that there exists a $d$-compressible $(C_4)_F$-switcher. Then there also exists a $d$-compressible $(C_6)_F$-switcher.
\end{lemma}

\proof
Let $U=\Set{u_1,\dots,u_6}$ and let $S'$ be the graph on $U\cup \Set{w_1,w_2}$ with edge set $\Set{u_1w_1,u_5w_1,u_2w_2,u_4w_2,w_1w_2}$. Define $f\colon U\rightarrow V(C_6)$ by $f(u_i):=\mathfrak{c}_i$ for $i\in[6]$. Moreover, let $\psi$ be defined such that $\psi{\restriction_U}=f$, $\psi(w_1):=\mathfrak{a}_1$, and $\psi(w_2):=\mathfrak{a}_2$, where $\mathfrak{a}_1,\mathfrak{a}_2$ are new vertices.
Let $K$ be the graph on $V(C_6)\cup \Set{\mathfrak{a}_1,\mathfrak{a}_2}$ with edge set $E(C_6)\cup\Set{\mathfrak{c}_1\mathfrak{a}_1,\mathfrak{c}_5\mathfrak{a}_1,\mathfrak{c}_2\mathfrak{a}_2,\mathfrak{c}_4\mathfrak{a}_2,\mathfrak{a}_1\mathfrak{a}_2}$. Observe that $(C_6,f,K,\psi)$ is a $3$-compression of $(S',U)$.

By Corollary~\ref{cor:model attaching} and our assumption, we can attach graphs $S_1,\dots,S_4$ to $S'$ such that
\begin{itemize}
\item $S_1$ is a $(u_1,w_1,u_5,u_6)_F$-switcher;
\item $S_2$ is a $(w_2,u_2,u_3,u_4)_F$-switcher;
\item $S_3$ is a $(u_1,u_2,w_2,w_1)_F$-switcher;
\item $S_4$ is a $(w_2,u_4,u_5,w_1)_F$-switcher;
\item $(S'\cup S_1 \cup \dots \cup S_4,U)$ is $d$-compressible with respect to $(C_6,f)$.
\end{itemize}
Let $S:=S'\cup S_1 \cup \dots \cup S_4$. It is easy to check that $S$ is a $(u_1,\dots,u_6)_F$-switcher. For example, $S\cup\Set{u_1u_2,u_3u_4,u_5u_6}$ can be decomposed into $S_1\cup \Set{u_1w_1,u_5u_6}$, $S_2\cup \Set{u_2w_2,u_3u_4}$, $S_3\cup \Set{u_1u_2,w_1w_2}$ and $S_4\cup \Set{u_4w_2,u_5w_1}$, which are all $F$-decomposable.
\endproof

\begin{lemma} \label{lem:C_4-switcher reduction}
Let $F$ be any graph and $d\ge 4$. Assume that there exists a $d$-compressible $(C_4)_F$-switcher with augmentation $\Set{\mathfrak{c}_1\mathfrak{c}_3,\mathfrak{c}_2\mathfrak{c}_4}$. Then there also exists a $d$-compressible $(C_4)_F$-switcher (with augmentation $\emptyset$).
\end{lemma}

\proof
Let $U=\Set{u_1,\dots,u_4}$ and define $f\colon U\rightarrow V(C_4)$ as $f(u_i):=\mathfrak{c}_i$ for $i\in [4]$. Let $w$ be a new vertex and let $\hat{S}$ be the graph on $U\cup \Set{w}$ with edge set $\Set{u_1w,wu_3}$. Define $\psi(u_i):=f(u_i)$ for $i\in[4]$ and $\psi(w):=\mathfrak{a}$, where $\mathfrak{a}$ is a new vertex.

We will first show that there exists a $d$-compressible $(C_4)_F$-switcher with augmentation $\Set{\mathfrak{c}_1\mathfrak{c}_3}$. To this end, let $K$ be the graph on $V(C_4)\cup\Set{\mathfrak{a}}$ with edge set $E(C_4)\cup \Set{\mathfrak{c}_1\mathfrak{c}_3} \cup \set{\mathfrak{c}_i\mathfrak{a}}{i\in[4]}$. Note that $(C_4\cup \Set{\mathfrak{c}_1\mathfrak{c}_3},f,K,\psi)$ is a $4$-compression of $(\hat{S},U)$. By Corollary~\ref{cor:model attaching} and our assumption, we can attach graphs $S_1,S_2$ to $\hat{S}$ such that
\begin{itemize}
\item $S_1$ is a $(u_1,u_2,u_3,w)_F$-switcher;
\item $S_2$ is a $(u_1,w,u_3,u_4)_F$-switcher;
\item $(S,U)$ is $d$-compressible with respect to $(C_4\cup \Set{\mathfrak{c}_1\mathfrak{c}_3},f)$, where $S:=\hat{S}\cup S_1 \cup S_2$.
\end{itemize}
Clearly, $S$ is a $(u_1,u_2,u_3,u_4)_F$-switcher.

We now conclude that there exists a $d$-compressible $(C_4)_F$-switcher without augmentation. Let $K'$ be the graph on $V(C_4)\cup\Set{\mathfrak{a}}$ with edge set $E(C_4)\cup \set{\mathfrak{c}_i\mathfrak{a}}{i\in[4]}$. Note that $(C_4,f,K',\psi)$ is a $4$-compression of $(\hat{S},U)$. By the above and Corollary~\ref{cor:model attaching}, we can attach $S_1',S_2'$ to $\hat{S}$ such that
\begin{itemize}
\item $S_1'$ is a $(u_1,u_2,u_3,w)_F$-switcher;
\item $S_2'$ is a $(u_1,w,u_3,u_4)_F$-switcher;
\item $(S',U)$ is $d$-compressible with respect to $(C_4,f)$, where $S':=\hat{S}\cup S_1' \cup S_2'$.
\end{itemize}
Then, $S'$ is a $(u_1,u_2,u_3,u_4)_F$-switcher.
\endproof

We will now gather the remaining building blocks to show that every graph $F$ is $(1-1/(\chi(F)+1))$-absorbing. Recall from Lemmas~\ref{lem:switch2transform} and \ref{lem:C4-to-C6 reduction} that in order to do this we have to construct a $(\chi(F)+1)$-compressible $(C_4)_F$-switcher and a $(\chi(F)+1)$-compressible $(K_{2,gcd(F)})_F$-switcher. In fact, we will only describe a $(C_4)_F$-switcher and then derive the $(K_{2,gcd(F)})_F$-switcher via a further reduction.

\begin{lemma} \label{lem:C4-switcher}
Let $F$ be any graph and $\chi:=\chi(F)$. There exists a $(\chi+1)$-compressible $(C_4)_F$-switcher.
\end{lemma}

\proof
Let $\ell:=|F|-1$ and suppose that $V(F)=\Set{f_0,\dots,f_\ell}$ and $f_0f_\ell\in E(F)$. Let $F':=F-f_\ell$.
Let $z_{0,0},\dots,z_{\ell-1,\ell-1}$ be $\ell^2$ new vertices. We define copies of $F'$ on them as follows: For every $i\in\Set{0,\dots,\ell-1}$, let $F^+_i$ be a copy of $F'$ on $\Set{z_{i,0},\dots,z_{i,\ell-1}}$ such that $z_{i,j}$ plays the role of $f_{i\oplus j}$ for all $j\in\Set{0,\dots,\ell-1}$, where $i\oplus j$ denotes addition modulo~$\ell$. Similarly, for every $j\in\Set{0,\dots,\ell-1}$, let $F^-_j$ be a copy of $F'$ on $\Set{z_{0,j},\dots,z_{\ell-1,j}}$ such that $z_{i,j}$ plays the role of $f_{i\oplus j}$ for all $i\in\Set{0,\dots,\ell-1}$. Note that the graphs $F^+_0,\dots,F^+_{\ell-1},F^-_0,\dots,F^-_{\ell-1}$ are all edge-disjoint. Let $S'$ be the graph obtained from the union of all these graphs by adding two new vertices $u_1,u_3$ and joining both of them to every $z_{i,j}$ with $f_{i\oplus j}\in N_F(f_\ell)$.
We claim that $S'$ has two natural $F$-decompositions $\cF_1,\cF_2$. To see this, for $k\in \Set{1,3}$, define
\begin{align*}
E^+_i(u_k) &:=\set{u_kz_{i,j}}{f_{i\oplus j}\in N_F(f_\ell), \, j\in \Set{0,\dots,\ell-1}}
\end{align*}
for all $i\in\Set{0,\dots,\ell-1}$ and 
\begin{align*}
E^-_j(u_k) &:=\set{u_kz_{i,j}}{f_{i\oplus j}\in N_F(f_\ell), \, i\in \Set{0,\dots,\ell-1}}
\end{align*}
for all $j\in\Set{0,\dots,\ell-1}$. Note that, for each $k\in\Set{1,3}$, $\bigcup_{i=0}^{\ell-1}E_i^+(u_k)=\bigcup_{j=0}^{\ell-1}E_j^-(u_k)$ consists of all edges at $u_k$ in $S'$. Thus,
$$S'=\bigcup_{i=0}^{\ell-1}(F_i^+\cup E_i^+(u_1))\cup (F_i^-\cup E_i^-(u_3))=\bigcup_{i=0}^{\ell-1}(F_i^+\cup E_i^+(u_3))\cup (F_i^-\cup E_i^-(u_1)).$$

Since for all $k\in\Set{1,3}$ and $i\in\Set{0,\dots,\ell-1}$, both $F_i^+\cup E_i^+(u_k)$ and $F_i^-\cup E_i^-(u_k)$ form a copy of $F$, this shows that $S'$ has two (natural) $F$-decompositions. 

In order to obtain the desired switcher, we make some minor modifications to~$S'$. Note that the vertex $z_{0,0}$ is contained in the copies $F^+_0$ and $F^-_0$, and $u_1z_{0,0},u_3z_{0,0}\in E(S')$ because $f_0f_\ell\in E(F)$. Let $S$ be the graph obtained from $S'-z_{0,0}$ as follows: Add new vertices $u_2,u_4$ and add all edges from $u_2$ to $N_{F^+_0}(z_{0,0})$ and all edges from $u_4$ to $N_{F^-_0}(z_{0,0})$. We claim that $S$ is a $(u_1,u_2,u_3,u_4)_F$-switcher. Clearly, $\Set{u_1,\dots,u_4}$ is independent in $S$. Let $\tilde{F}^+_0:=S[\Set{u_2,z_{0,1},\dots,z_{0,\ell-1}}]$ and $\tilde{F}^-_0:=S[\Set{u_4,z_{1,0},\dots,z_{\ell-1,0}}]$. Further, let
\begin{align*}
\tilde{E}^+_0(u_1) &:=E^+_0(u_1)\sm\Set{u_1z_{0,0}}; \\
\tilde{E}^-_0(u_1) &:=E^-_0(u_1)\sm\Set{u_1z_{0,0}}; \\
\tilde{E}^+_0(u_3) &:=E^+_0(u_3)\sm\Set{u_3z_{0,0}}; \\
\tilde{E}^-_0(u_3) &:=E^-_0(u_3)\sm\Set{u_3z_{0,0}}.
\end{align*}
Crucially, $\tilde{F}^+_0\cup \tilde{E}^+_0(u_1) \cup \Set{u_1u_2}$, $\tilde{F}^-_0\cup \tilde{E}^-_0(u_3) \cup \Set{u_3u_4}$, $\tilde{F}^+_0\cup \tilde{E}^+_0(u_3) \cup \Set{u_2u_3}$ and $\tilde{F}^-_0\cup \tilde{E}^-_0(u_1) \cup \Set{u_1u_4}$ are all copies of $F$. Moreover,
\begin{align*}
S &=(\tilde{F}^+_0 \cup \tilde{E}^+_0(u_1) )\cup (\tilde{F}^-_0 \cup \tilde{E}^-_0(u_3)) \cup \bigcup_{i=1}^{\ell-1}[(F^+_i\cup E^+_i(u_1)) \cup (F^-_i \cup E^-_i(u_3))] \\
  &=(\tilde{F}^+_0\cup \tilde{E}^+_0(u_3)) \cup (\tilde{F}^-_0 \cup \tilde{E}^-_0(u_1) ) \cup \bigcup_{i=1}^{\ell-1}[(F^+_i\cup E^+_i(u_3))\cup (F^-_i\cup E^-_i(u_1) )]. 
\end{align*}
With the above, it follows immediately that $S$ is a $(u_1,u_2,u_3,u_4)_F$-switcher.

Let $g\colon \Set{u_1,\dots,u_4}\rightarrow V(C_4)$ defined by $g(u_i):=\mathfrak{c}_i$ for $i\in[4]$. It remains to show that the model $(S,\Set{u_1,\dots,u_4})$ with root-compression $(C_4,g)$ is $(\chi+1)$-compressible.
Let $K$ be the graph obtained from $C_4$ by adding a complete graph on new vertices $\mathfrak{c}_5,\dots,\mathfrak{c}_{\chi+2}$ and all edges between $\Set{\mathfrak{c}_1,\dots,\mathfrak{c}_4}$ and $\Set{\mathfrak{c}_5,\dots,\mathfrak{c}_{\chi+2}}$. This way, $C_4$ is an induced subgraph of $K$ and $K$ has degeneracy at most $\chi+1$ rooted at $V(C_4)$. Let $c$ be a colouring of $F$ with colour set $\Set{\mathfrak{c}_1,\mathfrak{c}_2}\cup\Set{\mathfrak{c}_5\dots,\mathfrak{c}_{\chi+2}}$ such that $c(f_\ell)=\mathfrak{c}_1$ and $c(f_0)=\mathfrak{c}_2$. We define $\psi\colon V(S)\rightarrow V(K)$ as follows: Take $\psi(u_i):=\mathfrak{c}_i$ for all $i\in[4]$. For every vertex $z_{i,j}\in V(S)$, we let $\psi(z_{i,j}):=c(f_{i\oplus j})$. It is easy to see that $\psi$ is a homomorphism. For example, if $z_{i,j}z_{i',j'}\in E(S)$, then we must either have $i=i'$ and $z_{i,j}z_{i',j'}\in E(F_i^+)$ or $j=j'$ and $z_{i,j}z_{i',j'}\in E(F_j^-)$. In both cases, $f_{i\oplus j}f_{i'\oplus j'}\in E(F)$, implying $c(f_{i\oplus j})\neq c(f_{i'\oplus j'})$. Hence,  $\psi(z_{i,j})\neq \psi(z_{i',j'})$ and thus $\psi(z_{i,j})\psi(z_{i',j'})\in E(K)$. If $u_1z_{i,j}\in E(S)$, then $f_{i \oplus j}\in N_F(f_\ell)$ and so $\psi(u_1)=\mathfrak{c}_1=c(f_\ell)\neq c(f_{i\oplus j}) = \psi(z_{i,j})$. If $u_2z_{i,j}\in E(S)$, then $i=0$ and $f_j\in N_F(f_0)$ and so $\psi(u_2)=\mathfrak{c}_2=c(f_0)\neq c(f_{i\oplus j}) = \psi(z_{i,j})$.\COMMENT{
Since $K[\Set{\mathfrak{c}_1,\mathfrak{c}_2}\cup\Set{\mathfrak{c}_5\dots,\mathfrak{c}_{\chi+2}}]$ is a clique, $\psi(x)\neq \psi(y)$ is enough to have $\psi(x)\psi(y)\in E(K)$. Note that if $\psi(u_3)$ was $\mathfrak{c}_1$ and $\psi(u_4)$ was $\mathfrak{c}_2$, we were thereby done. But since, $N_K(\mathfrak{c}_1)=N_K(\mathfrak{c}_3)$ and $N_K(\mathfrak{c}_2)=N_K(\mathfrak{c}_4)$, we can as well take $\psi(u_3)$ and $\psi(u_4)$ as we did.}
Together with a similar argument for $u_3$ and $u_4$, this completes the proof.
\endproof

Using the $(C_4)_F$-switcher constructed above, we are able to construct a $(K_{2,d_F(v)})_F$-switcher for every $v\in V(F)$ (see Lemma~\ref{lem:K2r-switcher}). The following lemma implies that this suffices to obtain the desired $(K_{2,r})_F$-switcher.

\begin{lemma} \label{lem:degree2gcd-switcher reduction}
Let $F$ be any graph, $d\ge 0$ and $r\in \bN$ such that $gcd(F)\mid r$. Suppose that for all $v\in V(F)$, there exists a $d$-compressible $(K_{2,d_F(v)})_F$-switcher. 
Then there also exists a $d$-compressible $(K_{2,r})_F$-switcher.
\end{lemma}

\proof
Let $V_1$ and $V_2$ be multisubsets of $V(F)$ such that $r+n_1 = n_2$, where $n_1:=\sum_{v\in V_1}d_F(v)$ and $n_2:=\sum_{v\in V_2}d_F(v)$.

Let $U=\Set{u_1,\dots,u_{r+2}}$ and let $W$ be a set of $n_1$ new vertices. Define $f(u_i):=\mathfrak{p}_2$ for all $i\in[r]$, $f(u_{r+1}):=\mathfrak{p}_1$ and $f(u_{r+2}):=\mathfrak{p}_3$. Let $S'$ be the graph on $U\cup W$ with edge set $\set{u_{r+1}w,u_{r+2}w}{w\in W}$. Define $\psi$ such that $\psi{\restriction_U}=f$ and $\psi(w):=\mathfrak{p}_2$ for all $w\in W$. So $(P_2,f,P_2,\psi)$ is a $0$-compression of $(S',U)$. 

Let $E^+:=\set{u_{r+1}u_i}{i\in[r]}$ and $E^-:=\set{u_{r+2}u_i}{i\in[r]}$. In order to obtain an $(E^+,E^-)_F$-switcher, we partition $W$ into subsets $(U_v)_{v\in V_1}$ and partition $W\cup \Set{u_1,\dots,u_r}$ into subsets $(U_v)_{v\in V_2}$ such that $|U_v|=d_F(v)$ for all $v\in V_1\cup V_2$.

By Corollary~\ref{cor:model attaching} and our assumption, we can attach graphs $(S_v)_{v\in V_1\cup V_2}$ to $S'$ such that $S_v$ is a $(\set{u_{r+1}z}{z\in U_v},\set{u_{r+2}z}{z\in U_v})_F$-switcher and $(S,U)$ is $d$-compressible with respect to $(P_2,f)$, where $S:=S'\cup \bigcup_{v\in V_1\cup V_2}S_v$. Finally, observe that $$S\cup E^+=\bigcup_{v\in V_2}(S_v\cup\set{u_{r+1}z}{z\in U_v}) \cup \bigcup_{v\in V_1}(S_v\cup\set{u_{r+2}z}{z\in U_v})$$ is $F$-decomposable, and by symmetry, $S\cup E^-$ is also $F$-decomposable.
\endproof

\begin{lemma} \label{lem:K2r-switcher}
Let $F$ be a graph, let $r:=gcd(F)$ and $\chi:=\chi(F)$. There exists a $(\chi+1)$-compressible $(K_{2,r})_F$-switcher.
\end{lemma}

\proof
By the previous lemma, we can assume that $r=d_F(v)$ for some $v\in V(F)$. So let $w_1,\dots,w_r$ be an enumeration of $N_F(v)$. Moreover, fix a colouring $c$ of $F$ with colour set $[\chi]$ and $c(v)=1$.

Let $U=\Set{u_1,\dots,u_{r+2}}$ be a set of new vertices and let $S'$ be the graph on $U\cup V(F-v)$ with $E(S'):=D^+\cup D^- \cup E(F-v)$, where $D^+:=\set{u_{r+1}w_i}{i\in[r]}$ and $D^-:=\set{u_{r+2}w_i}{i\in[r]}$. Define $f\colon U\rightarrow V(P_2)$ as $f(u_i):=\mathfrak{p}_2$ for all $i\in[r]$, $f(u_{r+1}):=\mathfrak{p}_1$ and $f(u_{r+2}):=\mathfrak{p}_3$. Define $\psi$ such that $\psi{\restriction_{U}}=f$ and $\psi(x):=\mathfrak{a}_{c(x)}$ for all $x\in V(F-v)$, where $\mathfrak{a}_1,\dots,\mathfrak{a}_\chi$ are new vertices.

Let $K$ be the graph obtained from $P_2$ by adding a clique on $\mathfrak{a}_1,\dots,\mathfrak{a}_\chi$ and all edges between $\Set{\mathfrak{p}_1,\mathfrak{p}_3}$ and $\Set{\mathfrak{a}_2,\dots,\mathfrak{a}_\chi}$. So $P_2$ is an induced subgraph of $K$ and $K$ has degeneracy at most $\chi$ rooted at $V(P_2)$. Thus, $(P_2,f,K,\psi)$ is a $\chi$-compression of $(S',U)$.

For every $i\in[r]$, we have that $\psi(u_{r+1})=\mathfrak{p}_1$, $\psi(w_i)=\mathfrak{a}_{c(w_i)}$, $\psi(u_{r+2})=\mathfrak{p}_3$ and $\psi(u_i)=\mathfrak{p}_2$. By Lemma~\ref{lem:C4-switcher}, there exists a $(\chi+1)$-compressible $(C_4)_F$-switcher. Hence, by Corollary~\ref{cor:model attaching}, we can attach graphs $S_1,\dots,S_r$ to $S'$ such that
\begin{itemize}
\item $S_i$ is a $(u_{r+1},w_i,u_{r+2},u_i)_F$-switcher;
\item $(S,U)$ is $(\chi+1)$-compressible with respect to $(P_2,f)$, where $S:=S'\cup S_1 \cup \dots \cup S_r$.
\end{itemize}
We claim that $S$ is an $(E^+,E^-)_F$-switcher, where $E^+:= \set{u_{r+1}u_i}{i\in[r]}$ and $E^-:= \set{u_{r+2}u_i}{i\in[r]}$. To this end, observe that $S_1\cup \dots \cup S_r$ is an $(E^+\cup D^-,E^-\cup D^+)_F$-switcher. Hence, $$S\cup E^+=((F-v)\cup D^+)\cup(S_1\cup \dots \cup S_r\cup E^+\cup D^-)$$ and $$S\cup E^-=((F-v)\cup D^-)\cup(S_1\cup \dots \cup S_r\cup E^-\cup D^+)$$ are both $F$-decomposable.
\endproof

\section{Divisibility and threshold relations} \label{sec:lower bounds}

In this section, we make some observations regarding the relationship between the `auxiliary thresholds' $\delta_F^e,\delta_F^{vx},\delta_F^{0+},\delta_F^\ast$ and $\delta_F$. For this (and for later use when constructing extremal examples), we first gather some tools that allow us to remove a sparse subgraph of a given graph $G$ in order to make $G$ divisible.

\begin{lemma} \label{lem:make div}
Let $1/n\ll \mu,1/r$. Let $G$ be a graph on $n$ vertices with $\delta(G)\ge (1/2+\mu)n$ and let $\xi\colon V(G) \rightarrow \Set{0,\dots,r-1}$ be any function satisfying $r\mid \sum_{x\in V(G)}\xi(x)$. Then, there exists a subgraph $H$ of $G$ such that $\Delta(H)\le \mu n$ and $d_{G-H}(x)\equiv \xi(x) \mod{r}$ for all $x\in V(G)$.
\end{lemma}

\proof
We will find $H$ as a union of many small subgraphs which shift excess degree from one vertex to the next one. Let $wv\in E(K_{r,r})$ and let $u$ be a new vertex. Let $Q:=(K_{r,r}-\Set{wv})\cup\Set{uw}$. Let $Q'$ be the graph obtained from $K_{r,r}$ by subdividing one edge with a new vertex $z$. Note that $d_Q(u)=1$, $d_Q(v)\equiv -1\mod{r}$, $d_{Q'}(z)= 2$, while all other vertices of $Q$ and $Q'$ have degree $r$. 

\medskip
\noindent\emph{Claim 1:} $(Q,\Set{u,v})$ is $1/2$-embeddable.
\medskip

It is sufficient to show that if $x,x'\in V(G)$ are distinct, then $G$ contains a copy of $Q$ with $x,x'$ playing the roles of $u,v$. So suppose that $x$ and $x'$ are distinct vertices in $G$. Let $y$ be a neighbour of $x$ in $G$. Let $Y$ be a subset of $N_G(y)$ of size $\mu n/2$ disjoint from $\Set{x,x'}$. For every $y'\in Y$, $y'$ has at least $\mu n$ neighbours in $X:=N_G(x')\sm (Y\cup\Set{x,y})$. Let $H:=G[X,Y]$. Then, $e(H)\ge |Y|\mu n \ge \mu^2|H|^2/2$. Hence, $K_{r-1,r-1}$ must be a subgraph of $H$. Together with $x,x',y$, this yields the desired copy of $Q$.

\medskip
\noindent\emph{Claim 2:} $(Q',\Set{z})$ is $1/2$-embeddable.
\medskip

Let $x$ be any vertex in $G$. It is easy to see that $\delta(G[N_G(x)])\ge 2\mu n$.\COMMENT{and hence $e(G[N_G(x)])\ge \mu d_G(x)^2$ since $2e(G[N_G(x)])\ge d_G(x)2\mu n\ge 2\mu d_G(x)^2$} Thus, $K_{r,r}$ must be a subgraph of $G[N_G(x)]$, which together with $x$ contains a copy of $Q'$ with $x$ playing the role of $z$.

\medskip

Now, let $x_1,\dots,x_n$ be an enumeration of $V(G)$, and for each $i\in[n-1]$, let $a_i\in\Set{0,\dots,r-1}$ be such that $a_i\equiv \sum_{j\in[i]}(d_{G}(x_j)-\xi(x_j)) \mod{r}$. Let $a_n\in\Set{0,\dots,r-1}$ be such that $a_n\equiv e(G) \mod{r}$. Moreover, set $a_0:=0$.

For every $i\in [n-1]$, we want to find $a_i$ edge-disjoint copies of $Q$ in $G'$ such that $x_i,x_{i+1}$ play the roles of $u,v$. To this end, for every $i\in[n-1]$ and $j\in[a_i]$, let $\Lambda_{i,j}\colon \Set{u,v}\rightarrow V(G)$ be a $G$-labelling defined as $\Lambda_{i,j}(u):=\Set{x_i}$ and $\Lambda_{i,j}(v):=\Set{x_{i+1}}$. Note that there are at most $rn$ labellings and every vertex of $G$ is the image of a root at most $2r$ times. We can thus apply Lemma~\ref{lem:finding} to find edge-disjoint embeddings $(\phi_{i,j})_{i\in[n-1],j\in[a_i]}$ of $Q$ such that $\phi_{i,j}$ respects $\Lambda_{i,j}$ and $\Delta(H')\le \mu n/2$, where $H':=\bigcup_{i\in[n-1],j\in[a_i]}\phi_{i,j}(Q)$.

Note that $\delta(G-H')\ge (1/2+\mu/2)n$. Thus, we can greedily find $a_n$ edge-disjoint copies of $Q'$ in $G-H'$ with $x_n$ playing the role of $z$. Then, let $H$ be the union of $H'$ and the copies of $Q'$. Clearly, $\Delta(H)\le \mu n$.

Moreover, we have $d_H(x_i)\equiv a_i-a_{i-1} \equiv d_G(x_i)-\xi(x_i) \mod{r}$ for every $i\in [n-1]$. Finally, $$d_H(x_n)\equiv 2a_n-a_{n-1} \equiv 2e(G)-\sum_{j\in [n-1]}(d_G(x_j)-\xi(x_j))\equiv d_G(x_n)-\xi(x_n) \mod{r}.$$ Hence, we have $d_{G-H}(x)\equiv \xi(x) \mod{r}$ for all $x\in V(G)$.
\endproof

Roughly speaking, the above lemma allows us to make a graph $F$-degree-divisible. The following proposition allows us to make a graph $F$-edge-divisible, without destroying degree-divisibility. We will also use it in Section~\ref{sec:bip} to construct extremal examples.

\begin{prop} \label{prop:Walecki}
Let $F$ be a graph and $r:=gcd(F)$. Let $n\in \bN$ and suppose that $G$ is a graph on $n$ vertices with $\delta(G)\ge n/2+2e(F)(r+1)$. Then, for any number $e$ with $r\mid 2e$, there exists an $r$-divisible subgraph $H$ of $G$ such that $e(H) \equiv e \mod{e(F)}$ and $\Delta(H)\le 2e(F)r$.
\end{prop}

\proof
Let $V'\In V(G)$ be such that $gcd\Set{|V'|,e(F)}=1$ and $|V\sm V'|<e(F)$. Let $G':=G[V']$.
Clearly, $r\mid 2e(F)$. Let $a:=r$ if $r$ is odd and $a:=r/2$ if $r$ is even. So $a \mid e$ and $a \mid e(F)$. Let $0\le t < e(F)/a$ be an integer such that $e\equiv ta \mod{e(F)}$. 
Let $\alpha,\beta\in \mathbb{Z}$ be such that $\alpha e(F)+\beta|V'|=t$. We can assume that $0\le \beta < e(F)$.

Observe that $\delta(G')\ge |G'|/2 + 2(\beta a -1)$. Hence, by Dirac's theorem, we can take
$H$ to be the union of $\beta a$ edge-disjoint Hamilton cycles in $G'$. We then have $e(H)=|V'|\beta a= a(t-\alpha e(F)) \equiv e \mod{e(F)}$. Moreover, $H$ is $r$-divisible and $\Delta(H)\le 2e(F)r$.
\endproof

In order to show that $\delta_F^{0+}\le \delta_F^\ast$, we use the following result.

\begin{theorem}[Haxell and R\"odl \cite{HR}] \label{thm:Haxell-Rodl}
Let $F$ be a graph and $\eta>0$. There exists $n_0\in \bN$ such that whenever $G$ is a graph on $n\ge n_0$ vertices that has a fractional $F$-decomposition, then all but at most $\eta n^2$ edges of $G$ can be covered with edge-disjoint copies of $F$.
\end{theorem}

\begin{cor} \label{cor:threshold relations}
Let $F$ be any graph with $\chi(F)\ge 3$. Then $\max\Set{\delta_F^{0+},\delta_F^e}\le \delta_F^\ast\le \delta_F$ and $\delta_F^{vx}\le \delta_F$.
\end{cor}

\proof
Let $r:=gcd(F)$. Clearly, $1/2\le \delta_F^\ast\le \delta_F$.

Firstly, we show $\delta_F^{vx}\le \delta_F$. Let $1/n \ll \mu,1/|F|$ and let $G$ be a graph on $n$ vertices with $\delta(G) \ge (\delta_F + \mu)n$. Suppose that $x^\ast\in V(G)$ with $r\mid d_G(x^\ast)$. Let $G':=G-x^\ast$ and define $\xi(x):=r-1$ for all $x\in N_G(x^\ast)$ and $\xi(x):=0$ for all $x\in V(G')\sm N_G(x^\ast)$. Apply Lemma~\ref{lem:make div} to obtain a subgraph $H$ of $G'$ such that $d_{G'-H}(x)\equiv \xi(x) \mod{r}$ for all $x\in V(G')$ and $\Delta(H) \le \mu n/2$. Let $G'':=G-H$. Hence, $G''$ is $r$-divisible and $\delta(G'')\ge (\delta_F+\mu/2)n$. Apply Proposition~\ref{prop:Walecki} to $G''-x^\ast$ with $e(G'')$ playing the role of $e$ in order to obtain an $r$-divisible subgraph $H'$ of $G''-x^\ast$ such that $e(H') \equiv e(G'') \mod{e(F)}$ and $\Delta(H')\le 2e(F)r$. Let $G''':=G''-H'$. Observe that $G'''$ is $F$-divisible, $\delta(G''') \ge (\delta_F+ \mu/4)n$ and $N_{G'''}(x^\ast)=N_{G}(x^\ast)$. Now, $G'''$ has an $F$-decomposition. In particular, all edges at $x$ are covered.

We continue by showing $\delta_F^e\le \delta_F^\ast$. Let $1/n \ll \mu, 1/|F|$ and let $G$ be a graph on $n$ vertices with $\delta(G) \ge (\delta_F^\ast + \mu)n$. Suppose that $e'=x'y'\in E(G)$. We need to show that $e'$ is contained in a copy of $F$. Using Lemma~\ref{lem:make div} and Proposition~\ref{prop:Walecki} in the same way as above, it is easy to find a spanning subgraph $G'''$ of $G$ such that $G'''$ is $F$-divisible, $\delta(G''') \ge (\delta_F^\ast+ \mu/4)n$ and $e'\in E(G''')$. Now, $G'''$ has a fractional $F$-decomposition, which is only possible if every edge of $G'''$ is contained in a copy of $F$. In particular, $e'$ is contained in a copy of $F$.

Lastly, we show that $\delta_F^{\eta}\le \delta_F^\ast$ for all $\eta>0$, implying that $\delta_F^{0+}\le \delta_F^\ast$. Let $\eta>0$. Let $1/n \ll \mu, 1/|F|$ and let $G$ be a graph on $n$ vertices with $\delta(G) \ge (\delta_F^\ast + \mu)n$. We may assume that $\mu\ll \eta$. Using Lemma~\ref{lem:make div} and Proposition~\ref{prop:Walecki}, it is easy to find a subgraph $H$ of $G$ such that $G':=G-H$ is $F$-divisible, $\delta(G')\ge (\delta_F^\ast+\mu/4)n$ and $e(H) \le \mu n^2$. So $G'$ has a fractional $F$-decomposition. Thus, by Theorem~\ref{thm:Haxell-Rodl}, all but at most $\mu n^2$ edges of $G'$ can be covered by edge-disjoint copies of $F$, giving an $\eta$-approximate $F$-decomposition of $G$.\COMMENT{Apply Lemma~\ref{lem:make div} to obtain a subgraph $H$ of $G$ such that $r\mid d_{G-H}(x)$ for all $x\in V(G)$ and $\Delta(H) \le \mu n/2$. Let $G':=G-H$. Since $G'$ is $r$-divisible, we have $r\mid 2e(G')$. Hence, by Proposition~\ref{prop:Walecki} there exists an $r$-divisible subgraph $H'$ of $G'$ such that $e(H') \equiv e(G') \mod{e(F)}$ and $\Delta(H')\le 2e(F)r$. Let $G'':=G'-H'$. So $G''$ is $F$-divisible and $\delta(G'')\ge (\delta_F^\ast+\mu/4)n$, so $G''$ has a fractional $F$-decomposition. Thus, by Theorem~\ref{thm:Haxell-Rodl}, all but at most $\mu n^2$ edges of $G''$ can be covered by edge-disjoint copies of $F$. Since $e(G)-e(G'')\le \mu n^2$, this is an $\eta$-approximate $F$-decomposition of $G$. Hence, $\delta_F^{\eta}\le \delta_F^\ast$.}
\endproof

\section{Absorbers} \label{sec:absorbers}

The aim of this section is to prove the following lemma. Having done this, we can then bound $\delta_F$ in terms of $\delta^{0+}$, $\delta_F^{vx}$ and $\chi(F)$ (see Theorem~\ref{thm:almost main}).

\begin{lemma} \label{lem:transf2abs}
Let $F$ be any graph and let $\delta\ge 1/2$. If $F$ is $\delta$-transforming, then $F$ is $\delta$-absorbing.
\end{lemma}

Roughly speaking, we obtain the desired absorber by concatenating several suitable transformers. In particular, as intermediate steps, we `transform' a given graph into certain special graphs which we now define.

For a graph $F$, $e\in E(F)$ and $h\in \bN$, let $L(h;F,e)$ be the graph obtained from $h$ vertex-disjoint copies of $F$ by subdividing all copies of $e$ with one new vertex and identifying the new vertices. Note that $\chi(L(h;F,e))\le \max\Set{\chi(F),3}$.\COMMENT{Can colour the copies of $F$ identically and then colour the subdivision vertices with a colour that is not used for one of its two neighbours. In the bipartite case, we need a third colour.}

For a graph $F$, $v\in V(F)$, a graph $H$ and an orientation $\hat{H}$ of $H$, let $H^{att(\hat{H};F,v)}$ be the graph obtained from $H$ by adding $d^+_{\hat{H}}(x)$ copies of $F$ for every $x\in V(H)$ and identifying the copies of $v$ with $x$.\COMMENT{This is not the most precise description, but most readable.}

We need the following result from \cite{BKLO}. (Recall that we write  $H\rightsquigarrow H'$ if there is an edge-bijective homomorphism from $H$ to $H'$.)

\begin{lemma}[see {\cite[Lemma~8.7]{BKLO}}] \label{lem:loop graph}
Let $F$ be a graph, $uv\in E(F)$ and $r:=gcd(F)$. Then for every $F$-degree-divisible graph $H$ and any orientation $\hat{H}$ of $H$, there exists an $r$-regular graph $H_0$ such that $|H_0|\le 4e(H)e(F)$, $H_0\rightsquigarrow H^{att(\hat{H};F,v)}$ and $H_0\rightsquigarrow L(e(H);F,uv)$.
\end{lemma}
The statement here is slightly more general than that in \cite{BKLO}, as we do not require $F$ to be $r$-regular here, but the same proof goes through.\COMMENT{
\proof
Let $H$ be any $F$-degree-divisible graph and $\hat{H}$ any orientation of $H$. 
Let $(F_e)_{e\in E(\hat{H})}$ be a collection of vertex-disjoint copies of $F$, where $u_e,v_e\in F_e$ are the copies of $u,v\in F$. Let $H^{exp}$ be the graph obtained from $H\cup \bigcup_{e\in E(H)} F_e$ by deleting all edges of $H$, all edges $u_ev_e$, and adding the edges $xu_e$ and $yv_e$ for all $e=\vv{xy}\in E(\hat{H})$. Observe that identifying the vertices of $H$ in $H^{exp}$ yields a copy of $L(e(H);F,uv)$. Moreover, by identifying the vertices $x$ and $v_e$ for every edge $e=\vv{xy}\in E(\hat{H})$, one obtains a copy of $H^{att(\hat{H};F,v)}$. So $H^{exp}\rightsquigarrow H^{att(\hat{H};F,v)}$ and $H^{exp}\rightsquigarrow L(e(H);F,uv)$. Clearly, $H^{exp}$ is $F$-degree-divisible, so by splitting every vertex $x$ of $H^{exp}$ into $d_{H^{exp}}(x)/r$ vertices of degree $r$ each, we can obtain $H_0$, where $|H_0|=\sum_{x\in V(H^{exp})}d_{H^{exp}}(x)/r\le 2e(H^{exp}) =2(e(H)+e(H)e(F))$.
\endproof}

\lateproof{Lemma~\ref{lem:transf2abs}}
Let $F$ be a graph, $\chi:=\chi(F)$, $r:=gcd(F)$, $\delta\ge 1/2$ and assume that $F$ is $\delta$-transforming. Note that this implies that $\delta \ge 1-1/(\chi-1)$.
Let $1/n \ll 1/k_0',\eps \ll \alpha, 1/b \ll 1/m, \mu,1/|F|$ and suppose that $G$ is a graph on $n$ vertices with $\delta(G) \ge (\delta + \mu) n$ which has an $(\alpha,\eps,k)$-partition $V_1,\dots,V_k$ for some $k \le k_0'$, and $H$ is any $F$-divisible subgraph of $G$ of order at most $m$. We need to show that $G$ contains an $F$-absorber for $H$ of order at most $b$.
Let $\hat{H}$ be any orientation of $H$. Moreover, fix $uv\in E(F)$.

We first extend $H$ to $H^{att(\hat{H};F,v)}$ in $G$. Let $e_1,\dots,e_t$ be an enumeration of $E(\hat{H})$ and $e_i=\vv{x_iy_i}$. We want to find copies $F_1,\dots,F_t$ of $F$ in $G$ such that
\begin{enumerate}[label=(\roman*)]
\item $V(F_i)\cap V(H) = \Set{x_i}$ for all $i\in[t]$, where $x_i$ plays the role of $v$ in $F_i$;
\item $V(F_i-x_i)\cap V(F_j-x_j) = \emptyset$ for all $1\le i<j \le t$.
\end{enumerate}
Suppose that for some $s\in[t]$, we have already found $F_1,\dots,F_{s-1}$.\COMMENT{such that (i) and (ii) hold with $t$ being replaced by $s-1$} Let $X:=V(H)\cup V(F_1) \cup \dots \cup V(F_{s-1})$. For every $i\in[k]$, let $V_i':=V_i\sm (X\sm \Set{x_s})$, and let $G':=G[V_1'\cup\dots\cup V_k']$. Since $|X|\le m + m^2|F|$, we have $|V_i\sm V_i'| \le \eps^2|V_i|$ for all $i\in[k]$, and thus $V_1',\dots,V_k'$ is an $(\alpha/2,3\eps,k)$-partition of $G'$.\COMMENT{mentioned as fact where partition is defined} We now view $(F,\Set{v})$ as a model. Let $J$ be a graph of order one, say, with vertex $\mathfrak{j_1}$, and let $f$ map $v$ to $\mathfrak{j_1}$. Label $v$ with $\Set{x_s}$ and observe that this is an $(\alpha/2,3\eps,k)$-admissible $G'$-labelling. Since $(F,\Set{v})$ is $(\chi-1)$-compressible\COMMENT{Let $K$ be the complete graph on $\Set{\mathfrak{j_1},\dots,\mathfrak{j_\chi}}$ and $\psi$ a $\chi$-coloring of $F$ with colour set $\Set{\mathfrak{j_1},\dots,\mathfrak{j_\chi}}$ and $\psi(v)=\mathfrak{j_1}$.} and $\delta(G')\ge (1-1/(\chi-1)+\mu/2)|G'|$, Lemma~\ref{lem:rooted embedding} implies that there exists a copy $F_s$ of $F$ in $G'$ with $x_s$ playing the role of $v$. Then (i) and (ii) hold with $t$ replaced by $s$.

Let $H_{att}:=H\cup \bigcup_{j\in [t]} F_j$. So $H_{att}$ is isomorphic to $H^{att(\hat{H};F,v)}$.
Let $L:=L(e(H);F,uv)$. Let $pF$ be the vertex-disjoint union of $p$ copies of $F$, where $p:=e(H)/e(F)$, and let $pF_{att}:=(pF)^{att(p\hat{F};F,v)}$, where $p\hat{F}$ is some orientation of $pF$. By Lemma~\ref{lem:loop graph}, there exist $r$-regular graphs $H_0$ and $pF_0$ of order at most $2m^2e(F)$ such that $H_0\rightsquigarrow H_{att}$, $H_0\rightsquigarrow L$, $pF_0\rightsquigarrow pF_{att}$ and $pF_0 \rightsquigarrow L$.
Recall that $\chi(L)\le \max\Set{\chi,3}$. Since $H_0\rightsquigarrow L$, we have $\chi(H_0)\le \chi(L)$. Similarly, $\chi(pF_0)\le \chi(pF_{att})=\chi$. Since $\delta\ge 1/2$, we have that $\delta \ge 1-1/(\chi(H^\ast)-1)$ for all $H^\ast\in\Set{H_0,L,pF_{att},pF_0}$. Therefore, we can find copies of these graphs in $G$, and we may assume that these copies, which we call again $H_0,L,pF_{att},pF_0$, are vertex-disjoint and vertex-disjoint from $H_{att}$.\COMMENT{\Erd--Stone}

Since $F$ is $\delta$-transforming, $G$ contains subgraphs $T_1,\dots,T_4$ of order at most $b/4$ such that $T_1$ is an $(H_0,H_{att})_F$-transformer, $T_2$ is an $(H_0,L)_F$-transformer, $T_3$ is a $(pF_0,pF_{att})_F$-transformer and $T_4$ is a $(pF_0,L)_F$-transformer, and we may assume that $T_1,\dots,T_4$ are vertex-disjoint and vertex-disjoint from $H_{att},H_0,L,pF_{att},pF_0$ except for the obviously necessary intersections required by the definition of transformers, that is, e.g., $V(T_1)\cap V(T_2)=V(H_0)$ and $V(T_4)\cap V(L)=V(L)$.\COMMENT{So formally, we find the $T_i$'s in steps and for each step define a graph $G'$ by excluding a constant number of vertices and apply the definition of $\delta$-transforming with $0.9n,k_0',3\eps,\alpha/2,b/4,2m^2e(F),\mu/2$.} In particular, $H_{att}-H$, $H_0$, $L$, $pF_{att}$, $pF_0$, $T_1,\dots,T_4$ are edge-disjoint and contain no edge of $G[V(H)]$.
Let $$A:=(H_{att}-H)\cup T_1 \cup H_0 \cup T_2 \cup L \cup T_4 \cup pF_0 \cup T_3 \cup pF_{att}.$$
We claim that $A$ is an $F$-absorber for $H$. Indeed, $A$ has an $F$-decomposition since each of $H_{att}-H$, $T_1 \cup H_0$, $T_2 \cup L$, $T_4 \cup pF_0$, $T_3 \cup pF_{att}$ have $F$-decompositions. Secondly, $A\cup H$ has an $F$-decomposition as $H_{att}\cup T_1$, $H_0 \cup T_2$, $L \cup T_4$, $pF_0 \cup T_3$ and $pF_{att}$ are all $F$-decomposable. Moreover, $|A|\le b$ since $V(A)\In V(T_1\cup \dots \cup T_4)$.
\endproof

\begin{cor} \label{cor:absorbing}
Let $F$ be any graph and $\chi:=\chi(F)$. Then the following are true:
\begin{enumerate}[label=(\roman*)]
\item $F$ is $(1-1/(\chi+1))$-absorbing;
\item if $\delta_F<1-1/(\chi+1)$ and $\chi\ge 4$, then $F$ is $(1-1/\chi)$-absorbing;
\item if $\delta_F<1-1/\chi$ and $\chi \ge 5$, then $F$ is $(1-1/(\chi-1))$-absorbing.
\end{enumerate}
\end{cor}

\proof
Let $r:=gcd(F)$.
(i) By Lemma~\ref{lem:transf2abs}, it is enough to show that $F$ is $(1-1/(\chi+1))$-transforming. By Lemmas~\ref{lem:C4-to-C6 reduction} and~\ref{lem:C4-switcher}, there exists a $(\chi+1)$-compressible $(C_6)_F$-switcher. By Lemma~\ref{lem:K2r-switcher}, there exists a $(\chi+1)$-compressible $(K_{2,r})_F$-switcher. Therefore, Lemma~\ref{lem:switch2transform} implies that $F$ is $(1-1/(\chi+1))$-transforming.

(ii) By Lemma~\ref{lem:discretisation}, there exist a $\chi$-compressible $(C_4)_F$-switcher with augmentation $\Set{\mathfrak{c}_1\mathfrak{c}_3,\mathfrak{c}_2\mathfrak{c}_4}$ and a $\chi$-compressible $(K_{2,r})_F$-switcher with augmentation $\Set{\mathfrak{p}_1\mathfrak{p}_3}$. Since $\chi\ge 4$, Lemmas~\ref{lem:star-switcher reduction}, \ref{lem:C4-to-C6 reduction} and~\ref{lem:C_4-switcher reduction} imply that there exist a $\chi$-compressible $(C_6)_F$-switcher and a $\chi$-compressible $(K_{2,r})_F$-switcher without augmentations. Therefore, Lemma~\ref{lem:switch2transform} implies that $F$ is $(1-1/\chi)$-transforming. Lemma~\ref{lem:transf2abs} finally implies that $F$ is $(1-1/\chi)$-absorbing.

(iii) follows in the same way since $\chi-1\ge 4$.
\endproof

We are now able to deduce the following theorem, which is already close to Theorem~\ref{thm:main}.

\begin{theorem} \label{thm:almost main}
Let $F$ be a graph with $\chi:=\chi(F)$.
\begin{enumerate} [label=(\roman*)]
\item Then $\delta_F\le \max\Set{\delta_F^{0+},\delta_F^{vx},1-1/(\chi+1)}$.
\item If $\chi\ge 5$, then $\delta_F\in\Set{\max\Set{\delta_F^{0+},\delta_F^{vx}},1-1/\chi,1-1/(\chi+1)}$.
\end{enumerate}
\end{theorem}

\proof
Firstly, (i) follows from Theorem~\ref{thm:general dec} and Corollary~\ref{cor:absorbing}(i).

To prove (ii), suppose that $\chi\ge 5$. By Corollary~\ref{cor:threshold relations}, $\delta_F \ge \max\Set{\delta_F^{0+},\delta_F^{vx}}$. Hence, since $\delta_F^{0+}\ge 1-1/(\chi-1)$, Theorem~\ref{thm:general dec} and Corollary~\ref{cor:absorbing} imply that $\delta_F\in \Set{\max\Set{\delta_F^{0+},\delta_F^{vx}},1-1/\chi,1-1/(\chi+1)}$.
\endproof

\section{The decomposition threshold of bipartite graphs} \label{sec:bip}

In this section, we will determine $\delta_F$ for every bipartite graph. We first make some preliminary observations. We then consider the case $\tilde{\tau}(F)=1$ of coprime component sizes, where $\delta_F\in\Set{0,1/2}$. The main part of the section is devoted to the case $\tau(F)=1$. Finally, we consider extremal examples.

\subsection{Preliminary observations}
Whilst not much is known about the value of $\delta_F^{0+}$ in general, we can use the following observation in the bipartite case.

\begin{fact} \label{fact:bip approx}
If $F$ is bipartite, then $\delta_F^{0+}=0$.
\end{fact}

\proof
The \Erd--Stone theorem tells us that the Tur\'an density of $F$ is $0$. Hence, for all $\eta>0$, there exists an $n_0$ such that from every graph $G$ with $n\ge n_0$ vertices, we can greedily remove copies of $F$ until at most $\eta n^2$ edges remain. 
\endproof

In order to achieve upper bounds on $\delta_F$ using Theorem~\ref{thm:general dec}, we will investigate the absorbing behaviour of a given bipartite graph $F$. 
To this end, we recall the definitions of $\tau(F)$ and $\tilde{\tau}(F)$ (see Theorem~\ref{thm:bipartite char}). Let $F$ be a bipartite graph. A set $X\In V(F)$ is called \defn{$C_4$-supporting in $F$} if there exist distinct $a,b\in X$ and $c,d\in V(F)\sm X$ such that $ac,bd,cd\in E(F)$. We defined 
\begin{align*}
\tau(F)&:=gcd\set{e(F[X])}{X\In V(F)\text{ is not }C_4\text{-supporting in }F},\\
\tilde{\tau}(F)&:=gcd\set{e(C)}{C\text{ is a component of }F}.
\end{align*}

\begin{fact} \label{fact:bip par rel}
Let $F$ be bipartite. Then $\tau(F)\mid gcd(F)$ and $gcd(F)\mid \tilde{\tau}(F)$.
\end{fact}

\proof
For the first assertion, note that for every vertex $v\in V(F)$, $X:=N_F(v)\cup\Set{v}$ is not $C_4$-supporting and $e(F[X])=d_F(v)$.
Secondly, for every component $C$ of $F$, we have $gcd(C)\mid e(C)$, since the edges in $C$ can be counted by summing the degrees of the vertices in one colour class of $C$, and clearly $gcd(F)\mid gcd(C)$.
\endproof

Recall that we defined $(C_{2\ell})_F$-switchers and $(K_{2,r})_F$-switchers as special models in Section~\ref{sec:transform}.

\begin{prop} \label{prop:bip-star-switcher}
Let $F$ be a bipartite graph and $r\in \bN$ with $gcd(F)\mid r$. There exists a $0$-compressible $(K_{2,r})_F$-switcher.
\end{prop}

\proof
By Lemma~\ref{lem:degree2gcd-switcher reduction}, we can assume that $r=d_F(v)$ for some $v\in V(F)$. Let $c$ be a $\Set{1,2}$-colouring of $F$ with $c(v)=1$. Let $S$ be the graph obtained from $F$ by deleting all edges at $v$ and adding a new vertex $v'$. Since $F$ is bipartite, $U:=\Set{v,v'}\cup N_F(v)$ is independent in $S$. Let $E^+:=\set{vu}{u\in N_F(v)}$ and $E^-:=\set{v'u}{u\in N_F(v)}$. Clearly, $S$ is an $(E^+,E^-)_F$-switcher. Define $f(u):=\mathfrak{p}_2$ for all $u\in N_F(v)$, $f(v):=\mathfrak{p}_1$ and $f(v'):=\mathfrak{p}_3$. Moreover, define $\psi(x):=\mathfrak{p}_{c(x)}$ for all $x\in V(F)$ and $\psi(v'):=\mathfrak{p}_3$. Then, $(P_2,f,P_2,\psi)$ is a $0$-compression of $(S,U)$.
\endproof

Let $F$ be a graph and $d\ge 0$. Let $(T,\Set{u_1,\dots,u_4})$ be a model such that $T$ is a $(\Set{u_1u_2},\Set{u_3u_4})_F$-switcher. Note that a necessary condition for this to exist is that $gcd(F)=1$. 

If $(T,\Set{u_1,\dots,u_4})$ is $d$-compressible with respect to the root-compression $(P_1,f)$, where $f(u_1)=f(u_3)=\mathfrak{p}_1$ and $f(u_2)=f(u_4)=\mathfrak{p}_2$, then we call $(T,\Set{u_1,\dots,u_4})$ a \defn{$d$-compressible internal $(P_1)_F$-teleporter}.

Let $2P_1$ be the graph with $V(2P_1)=\Set{\mathfrak{p}_1,\mathfrak{p}_2,\mathfrak{p}_1',\mathfrak{p}_2'}$ and $E(2P_1)=\Set{\mathfrak{p}_1\mathfrak{p}_2,\mathfrak{p}_1'\mathfrak{p}_2'}$.
If $(T,\Set{u_1,\dots,u_4})$ is $d$-compressible with respect to the root-compression $(2P_1,f')$, where $f'(u_1)=\mathfrak{p}_1$, $f'(u_2)=\mathfrak{p}_2$, $f'(u_3)=\mathfrak{p}_1'$ and $f'(u_4)=\mathfrak{p}_2'$, then we call $(T,\Set{u_1,\dots,u_4})$ a \defn{$d$-compressible external $(P_1)_F$-teleporter}.

Note that by Fact~\ref{fact:model shrink}, every $d$-compressible external $(P_1)_F$-teleporter is also a $d$-compressible internal $(P_1)_F$-teleporter. Loosely speaking, when considering an $(\alpha,\eps,k)$-partition of a graph $G$, then an internal $(P_1)_F$-teleporter allows us to switch between two edges lying in the same regular pair of clusters, whereas an external $(P_1)_F$-teleporter would allow us to switch between two edges that may belong to different regular pairs.

The following proposition gives an easy way of constructing an internal $(P_1)$-teleporter. We will use it as a tool in the proof of Lemmas~\ref{lem:bip-ex teleporter} and~\ref{lem:bipC6switcher}.

\begin{prop} \label{prop:bip-int teleporter}
Let $F$ be a bipartite graph with $gcd(F)=1$. Then there exists a $0$-compressible internal $(P_1)_F$-teleporter $(T,\Set{u_1,\dots,u_4})$.
\end{prop}

Note that in this case, $T$ itself has a homomorphism onto $P_1$ and is thus bipartite.

\proof
Let $T'$ be the graph with vertex set $\Set{u_1,\dots,u_4,w}$ and edge set $\Set{u_2w,wu_4}$. Consider the model $(T',\Set{u_1,\dots,u_4})$ with compression $(P_1,f,P_1,\psi)$, where $\psi(u_i):=f(u_i):=\mathfrak{p}_1$ for $i\in\Set{1,3}$, $\psi(u_i):=f(u_i):=\mathfrak{p}_2$ for $i\in\Set{2,4}$, and $\psi(w):=\mathfrak{p}_1$. By Proposition~\ref{prop:bip-star-switcher}, there exists a $0$-compressible $(K_{2,1})_F$-switcher. So by Corollary~\ref{cor:model attaching}, we can attach graphs $S_1,S_2,S_3$ to $T'$ such that
\begin{itemize}
\item $S_1$ is a $(\Set{u_1u_2},\Set{wu_2})_F$-switcher;
\item $S_2$ is a $(\Set{u_2w},\Set{u_4w})_F$-switcher;
\item $S_3$ is a $(\Set{wu_4},\Set{u_3u_4})_F$-switcher;
\item $(T,\Set{u_1,\dots,u_4})$ is $0$-compressible with respect to $(P_1,f)$, where $T:=T'\cup S_1 \cup S_2 \cup S_3$.
\end{itemize}\COMMENT{Here, the homomorphisms $\beta_i$ used to attach the $S_i$'s are: $\beta_i(\mathfrak{p}_1)=\beta_i(\mathfrak{p}_3):=\mathfrak{p}_1$ and $\beta_i(\mathfrak{p}_2):=\mathfrak{p}_2$ for $i\in\Set{1,3}$, while $\beta_2(\mathfrak{p}_1)=\beta_2(\mathfrak{p}_3):=\mathfrak{p}_2$ and $\beta_2(\mathfrak{p}_2):=\mathfrak{p}_1$.}
Then, $T$ is clearly a $(\Set{u_1u_2},\Set{u_3u_4})_F$-switcher.
\endproof

\subsection{Coprime component sizes}

We first analyse the case when $\tilde{\tau}(F)=1$. 
In the proof of the following lemma, we construct an external $(P_1)_F$-teleporter, which we will use in the proof of Lemma~\ref{lem:bip abs 2} to show that $F$ is $0$-absorbing.

\begin{lemma} \label{lem:bip-ex teleporter}
Let $F$ be a bipartite graph with $\tilde{\tau}(F)=1$. Then there exists a $0$-compressible external $(P_1)_F$-teleporter.
\end{lemma}

\proof
Let $M_1,M_2$ be disjoint multisets containing components of $F$ such that $1+\sum_{C\in M_1}e(C)=\sum_{C\in M_2}e(C)$. Let $U=\Set{u_1,\dots,u_4}$ and define $f'\colon U\rightarrow V(2P_1)$ as $f'(u_1):=\mathfrak{p}_1$, $f'(u_2):=\mathfrak{p}_2$, $f'(u_3):=\mathfrak{p}_1'$ and $f'(u_4):=\mathfrak{p}_2'$.

We want to construct a $0$-compressible model $(T,U)$ with respect to $(2P_1,f')$ such that $T$ is a $(\Set{u_1u_2},\Set{u_3u_4})_F$-switcher. 

Fix some component $C^\ast\in M_2$ and let $vw\in E(C^\ast)$. Let $c$ be a $\Set{1,2}$-colouring of $F$ such that $c(v)=1$ and $c(w)=2$.

Let $F^+_{C^\ast}$ be a copy of $C^\ast-vw$ such that $u_1,u_2$ play the roles of $v,w$. Likewise, let $F^-_{C^\ast}$ be a copy of $C^\ast-vw$ such that $u_3,u_4$ play the roles of $v,w$. Let $F^0_{C^\ast}$ be a copy of $F-V(C^\ast)$. Moreover, for each component $C\in M_1\cup(M_2\sm\Set{C^\ast})$, let $F^+_C$ and $F^-_C$ be two copies of $C$ and let $F^0_C$ be a copy of $F-V(C)$. We may assume that all these copies are vertex-disjoint. Let $T':=\bigcup_{C\in M_1\cup M_2}F_C^+\cup F_C^- \cup F_C^0$.

For all $x\in V(T')$ that belong to some $F^0_C$ or $F^+_C$, define $\psi(x):=\mathfrak{p}_{c(z)}$, where $z$ is the vertex of $F$ whose role $x$ is playing. Likewise, for all $x\in V(T')$ that belong to some $F^-_C$, define $\psi(x):=\mathfrak{p}'_{c(z)}$, where $z$ is the vertex of $F$ whose role $x$ is playing. Observe that $\psi\colon T' \rightarrow 2P_1$ is a homomorphism such that $\psi{\restriction_U}=f'$. Hence, $(2P_1,f',2P_1,\psi)$ is a $0$-compression of $(T',U)$.

Let $E^\odot_i:=E(\bigcup_{C\in M_i}F^\odot_C)$ for $i\in \Set{1,2}$ and $\odot\in\Set{+,-}$. By the definitions of $M_1,M_2$ and $C^\ast$, we have $|E^+_1|=|E^+_2|$ and $|E^-_1|=|E^-_2|$.
Let $\phi^+\colon E^+_1\rightarrow E^+_2$ and $\phi^-\colon E^-_1\rightarrow E^-_2$ be arbitrary bijections. For every edge $e\in E^+_1$, we have $e=xy$ and $\phi^+(e)=x'y'$ for suitable distinct $x,y,x',y'$ such that $\psi(x)=\psi(x')=\mathfrak{p}_1$ and $\psi(y)=\psi(y')=\mathfrak{p}_2$. Similarly, for every edge $e\in E^-_1$, we have $e=xy$ and $\phi^-(e)=x'y'$ for suitable distinct $x,y,x',y'$ such that $\psi(x)=\psi(x')=\mathfrak{p}_1'$ and $\psi(y)=\psi(y')=\mathfrak{p}_2'$. By Fact~\ref{fact:bip par rel} and Proposition~\ref{prop:bip-int teleporter}, there exists a $0$-compressible internal $(P_1)_F$-teleporter. We can therefore use Corollary~\ref{cor:model attaching} to attach graphs $(T_e)_{e\in E^+_1\cup E^-_1}$ to $T'$ such that
\begin{itemize}
\item $T_e$ is an $(\Set{e},\Set{\phi^+(e)})_F$-switcher for every $e\in E^+_1$;
\item $T_e$ is an $(\Set{e},\Set{\phi^-(e)})_F$-switcher for every $e\in E^-_1$;
\item $(T,U)$ is $0$-compressible with respect to $(2P_1,f')$, where $T:=T'\cup S^+\cup S^-$, $S^+:=\bigcup_{e\in E^+_1}T_e$ and $S^-:=\bigcup_{e\in E^-_1}T_e$.
\end{itemize}
It remains to show that $T$ is a $(\Set{u_1u_2},\Set{u_3u_4})_F$-switcher. Note that $S^+$ is an $(E^+_1,E^+_2)_F$-switcher and $S^-$ is an $(E^-_1,E^-_2)_F$-switcher. Thus, 
\begin{align*}
T\cup \Set{u_1u_2} &= (F^0_{C^\ast}\cup F^+_{C^\ast}\cup \Set{u_1u_2})  \cup (S^+ \cup E^+_1) \cup (S^- \cup E^-_2) \\
                   & \cup \bigcup_{C\in M_1}(F^0_C \cup F^-_C) \cup \bigcup_{C\in M_2\sm\Set{C^\ast}}(F^0_C \cup F^+_C), \\
T\cup \Set{u_3u_4} &= (F^0_{C^\ast}\cup F^-_{C^\ast}\cup \Set{u_3u_4})  \cup (S^+ \cup E^+_2) \cup (S^- \cup E^-_1) \\
                   & \cup \bigcup_{C\in M_1}(F^0_C \cup F^+_C) \cup \bigcup_{C\in M_2\sm\Set{C^\ast}}(F^0_C \cup F^-_C),
\end{align*}
are both $F$-decomposable.
\endproof

\begin{lemma} \label{lem:bip abs 2}
Let $F$ be bipartite and $\tilde{\tau}(F)=1$. Then $F$ is $0$-absorbing.
\end{lemma}

Using Lemma~\ref{lem:bip-ex teleporter}, we will be able to transform any given leftover into a union of copies of $F$ resulting in the desired absorber.

\proof
By Lemma~\ref{lem:bip-ex teleporter}, there exists a $0$-compressible external $(P_1)_F$-teleporter $(T,\Set{u_1,\dots,u_4})$, say, that is, $T$ is a $(\Set{u_1u_2},\Set{u_3u_4})_F$-switcher and $(T,\Set{u_1,\dots,u_4})$ is $0$-compressible with respect to $(2P_2,f')$, where $f'(u_1)=\mathfrak{p}_1$, $f'(u_2)=\mathfrak{p}_2$, $f'(u_3)=\mathfrak{p}_1'$ and $f'(u_4)=\mathfrak{p}_2'$.

Let $t:=|T|$ and assume that $1/n \ll 1/k_0',\eps \ll \alpha, 1/b \ll 1/m,\mu,1/|F|$. Since $t$ only depends on $F$ this implies that $\alpha,1/b\ll 1/t$.\COMMENT{Should we say choose $T$ such that $|T|$ is minimised?} Suppose that $G$ is a graph on $n$ vertices with $\delta(G) \ge \mu n$ that has an $(\alpha,\eps,k)$-partition $V_1,\dots,V_k$ for some $k \le k_0'$. Suppose also that $H$ is any $F$-divisible subgraph of $G$ of order at most $m$. We are to show that $G$ contains an $F$-absorber for $H$ of order at most $b$.
Let $p:=e(H)/e(F)$ and let $pF$ be the vertex-disjoint union of $p$ copies of $F$. Clearly, we can find a copy $H'$ of $pF$ as a subgraph in $G$\COMMENT{\Erd-Stone} such that $H'$ is vertex-disjoint from $H$. Let $e_1,\dots,e_h$ be an enumeration of the edges of $H$ and let $e_1',\dots,e_h'$ be an enumeration of the edges of $H'$. We now want to find an $(\Set{e_i},\Set{e_i'})_F$-switcher for all $i\in[h]$. More precisely, we want to find edge-disjoint copies $T_1,\dots,T_h$ of $T$ in $G$ such that
\begin{enumerate}[label=(\roman*)]
\item $T_i$ is an $(\Set{e_i},\Set{e_i'})_F$-switcher;
\item $T_i[V(H)]$ and $T_i[V(H')]$ are empty.
\end{enumerate}
Once again, we find them one by one using Lemma~\ref{lem:rooted embedding}. Suppose that for some $s\in[h]$, we have already found $T_1,\dots,T_{s-1}$. Write $e_s=xy\in E(H)$ and $e_s'=x'y'\in E(H')$.\COMMENT{all distinct since $H'$ was chosen vertex-disjoint from $H$} Let $X:=V(H)\cup V(H')\cup V(T_1) \cup \dots \cup V(T_{s-1})$. For every $i\in[k]$, let $V_i':=V_i\sm(X\sm\Set{x,y,x',y'})$. Then, $V_1',\dots,V_k'$ is an $(\alpha/2,3\eps,k)$-partition of $G':=G[V_1'\cup \dots \cup V_k']$ and $\delta(G')\ge \mu|G'|/2$. Let $\Lambda\colon \Set{u_1,\dots,u_4} \rightarrow V(G')$ be defined as $\Lambda(u_1):=\Set{x}$, $\Lambda(u_2):=\Set{y}$, $\Lambda(u_3):=\Set{x'}$ and $\Lambda(u_4):=\Set{y'}$. We claim that $\Lambda$ is $(\alpha/2,3\eps,k)$-admissible. Clearly, $\Lambda$ respects $(2P_2,f')$. Let $R$ be the reduced graph of $V_1',\dots,V_k'$ (with respect to $G'$) and $\sigma\colon V(G') \rightarrow R$ the cluster function. Define $j:V(2P_1)\rightarrow V(R)$ by $j(\mathfrak{p}_1):=\sigma(x)$, $j(\mathfrak{p}_2):=\sigma(y)$, $j(\mathfrak{p}_1'):=\sigma(x')$ and $j(\mathfrak{p}_2'):=\sigma(y')$. Since $xy,x'y'\in E(G')$, we have $\sigma(x)\sigma(y),\sigma(x')\sigma(y')\in E(R)$, so \ref{adm:hom} holds. Observe that $W_{\mathfrak{p}_1}=\Set{y}$.
Since $d_{G'}(y,V_{\sigma(x)})>0$, we have $d_{G'}(W_{\mathfrak{p}_1},V_{j(\mathfrak{p}_1)}) \ge \alpha|V_{j(\mathfrak{p}_1)}|/2$. The same applies to $\mathfrak{p}_2,\mathfrak{p}_1',\mathfrak{p}_2'$, so  \ref{adm:conform} holds as well. \ref{adm:label size} and \ref{adm:label loc} hold trivially.\COMMENT{as `$J_2=\emptyset$'} Thus, Lemma~\ref{lem:rooted embedding} implies that there exists an embedding $\rho$ of $(T,\Set{u_1,\dots,u_4})$ into $G'$ respecting $\Lambda$. Take $T_s:=\rho(T)$.

Let $A:=H'\cup T_1 \cup \dots \cup T_h$. Then $|A|\le e(H)|F| + e(H)t \le b$ and $A[V(H)]$ is empty. Since $A=\bigcup_{i\in[h]}(T_i\cup \Set{e_i'})$ is $F$-decomposable and $A\cup H=H'\cup \bigcup_{i\in[h]}(T_i\cup \Set{e_i})$ is $F$-decomposable too, $A$ is an $F$-absorber for $H$ in $G$ of order at most $b$.
\endproof

\subsection{Coprime non-\texorpdfstring{$C_4$}{C4}-supporting sets}

Here we show that $\tau(F)=1$ implies that $F$ is $1/2$-absorbing (see Corollary~\ref{cor:bip abs 1}). The remaining step towards this goal is to show that there exists a $2$-compressible $(C_6)_F$-switcher. In fact, we will construct a $0$-compressible $(C_6)_F$-switcher. For this we need to construct a model $(S,\Set{u_1,\dots,u_6})$ and a homomorphism $\psi\colon S\to C_6$ such that $\psi(u_i)=\mathfrak{c}_i$ for all $i\in[6]$ and such that both $S\cup \Set{u_1u_2,u_3u_4,u_5u_6}$ and $S\cup \Set{u_2u_3,u_4u_5,u_6u_1}$ are $F$-decomposable.

\begin{fact} \label{fact:support connected}
Let $F$ be bipartite. Then $$\tau(F)=gcd\set{e(F[X])}{X\In V(F)\text{ is not }C_4\text{-supporting and }F[X]\text{ is connected}}.$$
\end{fact}
\proof It is sufficient to show that for any $X\In V(F)$ that is not $C_4$-supporting, we have that $gcd\set{e(F[Y])}{Y\In V(F)\text{ is not }C_4\text{-supporting and }F[Y]\text{ is connected}}$ divides $e(F[X])$. Therefore, suppose that $X\In V(F)$ is not $C_4$-supporting. There is a partition $X=Y_1\cup \dots \cup Y_t$ such that $F[Y_1],\dots,F[Y_t]$ are the components of $F[X]$. Then, $e(F[X])=e(F[Y_1])+\dots+e(F[Y_t])$ and every $Y_i$ itself is not $C_4$-supporting.
\endproof

For graphs $H$ and $J$, a homomorphism $\psi\colon H\rightarrow J$, and vertices $\mathfrak{j}_1,\dots,\mathfrak{j}_k\in V(J)$, we slightly abuse notation for the sake of readability and write $H[\psi^{-1}(\mathfrak{j}_1\dots\mathfrak{j}_k)]$ instead of $H[\psi^{-1}(\Set{\mathfrak{j}_1,\dots,\mathfrak{j}_k})]$.

The following lemma is the key building block for the construction of the desired $C_6$-switcher.

\begin{lemma}\label{lem:step1}
Let $F$ be bipartite and $\tau(F)=1$. Then there exist $F$-decomposable graphs $G_+,G_-$, an edge $e^0\in E(G_+)$ and a homomorphism $\rho\colon G_+\cup G_-\to C_6$ such that $\rho(e^0)= \mathfrak{c}_1\mathfrak{c}_2$ and $G_-[\rho^{-1}(\mathfrak{c}_1\mathfrak{c}_2)] = G_+[\rho^{-1}(\mathfrak{c}_1\mathfrak{c}_2)]-\Set{e^0}$.
\end{lemma}

Note that the condition $\tau(F)=1$ is crucial here. Indeed, if $\rho'\colon F\to C_6$ is a homomorphism, then the number of edges mapped to $\mathfrak{c}_1\mathfrak{c}_2$ is divisible by $\tau(F)$ as $\rho'^{-1}(\Set{\mathfrak{c}_1,\mathfrak{c}_2})$ is not $C_4$-supporting. Hence, $\tau(F)$ must divide both $e(G_+[\rho^{-1}(\mathfrak{c}_1\mathfrak{c}_2)])$ and $e(G_-[\rho^{-1}(\mathfrak{c}_1\mathfrak{c}_2)])$.

\proof
Let $$\cC:=\set{F[X]}{X\In V(F)\text{ is not }C_4\text{-supporting and }F[X]\text{ is connected}}.$$ A graph $G$ is called \defn{$\cC$-decomposable} if $G$ can be decomposed into copies of elements of $\cC$.
It is sufficient to show that there exist a graph $G^0$, a homomorphism $\rho^0\colon G^0\to C_6[\Set{\mathfrak{c}_1,\mathfrak{c}_2}]\simeq P_1$ and an edge $e^0\in E(G^0)$ such that both $G^0$ and $G^0-\Set{e_0}$ are $\cC$-decomposable. Indeed, if $F[X]\in \cC$ appears in the decomposition of $G^0$ (or $G^0-\Set{e_0}$), then one can extend $F[X]$ to a copy of $F$ and extend $\rho^0$ appropriately without mapping new edges to $\mathfrak{c}_1\mathfrak{c}_2$. Clearly, all these extensions can be carried out edge-disjointly.

We now construct $G_0$. By assumption and Fact~\ref{fact:support connected}, $gcd(\set{e(C)}{C\in \cC})=1$. Thus, there exist disjoint sets $M_+,M_-$ containing copies of elements of $\cC$ such that $1+\sum_{C\in M_-}e(C)=\sum_{C\in M_+}e(C)$. We may assume that all the elements of $M_+\cup M_-$ are vertex-disjoint. Let $G':=\bigcup M_+ \cup \bigcup M_-$. Clearly, $G'$ is bipartite, i.e.~there exists a homomorphism $\rho'\colon G'\rightarrow P_1$.
Let $E^\odot:=E(\bigcup M_\odot)$ for $\odot\in\Set{+,-}$. Hence, $|E^+|=|E^-|+1$. Let $e^0$ be any edge in $E^+$. Let $\tilde{F}$ be the disjoint union of all elements of $\cC$. We want to construct an $(E_+\sm\Set{e^0},E_-)_{\tilde{F}}$-switcher $S$. To this end, let $\phi$ be an arbitrary bijection from $E_-$ to $E_+\sm\Set{e^0}$, and for every $e\in E_-$, let $V_e$ be the set of vertices incident to $e$ or $\phi(e)$.

Clearly, $\tilde{\tau}(\tilde{F})=1$. By Fact~\ref{fact:bip par rel}, it follows that $gcd(\tilde{F})=1$. Thus, Proposition~\ref{prop:bip-int teleporter} implies that there exists a $0$-compressible internal $(P_1)_{\tilde{F}}$-teleporter.
Therefore, there exist graphs $(T_e)_{e\in E_-}$ such that
\begin{itemize}
\item $T_e$ is an $(\Set{e},\Set{\phi(e)})_{\tilde{F}}$-switcher for each $e\in E_-$;
\item $V(G')\cap V(T_e)=V_e$ for each $e\in E_-$;
\item $V(T_e)\cap V(T_{e'})=V_e \cap V_{e'}$ for all distinct $e,e'\in E_-$;
\item there exists a homomorphism $\rho_e\colon T_e \rightarrow P_1$ with $\rho_e(x)=\rho'(x)$ for all $x\in V_e$.
\end{itemize}
Let $S:=\bigcup_{e\in E_-}T_e$ and $G^0:=G'\cup S$. Clearly, $S$ is an $(E_+\sm\Set{e^0},E_-)_{\tilde{F}}$-switcher and $\rho^0:=\rho'\cup \bigcup_{e\in E_-}\rho_e$ is a homomorphism $G^0\rightarrow P_1$.

In particular, both $S\cup (E^+\sm\Set{e^0})$ and $S\cup E^-$ are $\cC$-decomposable.
Thus, $G^0=(S\cup E^-)\cup \bigcup M_+$ is $\cC$-decomposable. Similarly, $G^0-\Set{e^0}=(S\cup(E^+\sm\Set{e^0})) \cup \bigcup M_-$ is $\cC$-decomposable.
\endproof

\begin{figure}[htbp]
\centering
\begin{tikzpicture}[scale=0.6]
		
				\draw [green, pattern = crosshatch, pattern color=green] (-6,3.75) rectangle (6,3.2);
				\draw [green, pattern = crosshatch, pattern color=green] (-6.75,3) rectangle (-6.2,0);
				\draw [green, pattern = crosshatch, pattern color=green] (6.75,3) rectangle (6.2,0);
				
				\draw [blue, pattern = crosshatch dots, pattern color= blue] (-6,-0.75) rectangle (6,-0.2);
				\draw [blue, pattern = crosshatch dots, pattern color= blue]  (-5.8,3) rectangle (-5.25,0);
				\draw [blue, pattern = crosshatch dots, pattern color= blue]  (5.8,3) rectangle (5.25,0);

				\draw [red, pattern =north west lines
, pattern color= red]  (-6,0.75) rectangle (6,0.2);

				\draw [fill = white] (6,3) circle (1);
				\draw [fill = white] (-6,3) circle (1);
				\draw [fill = white] (6,0) circle (1);
				\draw [fill = white] (-6,0) circle (1);
				\draw [fill = white] (2,0) circle (1);
				\draw [fill = white] (-2,0) circle (1);
				
				\filldraw[fill=black] (6,3) circle (2pt);
				\filldraw[fill=black] (-6,3) circle (2pt);
				\filldraw[fill=black] (6,0) circle (2pt);
				\filldraw[fill=black] (-6,0) circle (2pt);
				
		\begin{scope}[red,line width = 1pt]
		\draw (6,3) -- (-6,3);
		\end{scope}
		
		\begin{scope}[blue,line width = 1pt]
		\draw (6,3) -- (6,0);
		\draw (-6,3) -- (-6,0);
		\end{scope}
			
		\node  at (-6.4,3)  {$u_1$};
		\node  at (6.4,3)  {$u_2$};
		\node  at (6.4,0)  {$u_3$};
		\node  at (-6.4,0)  {$u_6$};
		
		\node  at (-2,-1.4)  {$\mathfrak{a}_1$};
		\node  at (2,-1.4)  {$\mathfrak{a}_2$};
		\node  at (-7.4,3)  {$\mathfrak{c}_1$};
		\node  at (7.4,3)  {$\mathfrak{c}_2$};
		\node  at (7.4,0)  {$\mathfrak{c}_3$};
		\node  at (-7.4,0)  {$\mathfrak{c}_6$};	
		
		\draw [green, pattern = crosshatch, pattern color=green] (-5,-2) rectangle (-6,-2.5);
		\node  at (-4,-2.25)  {$\tilde{H}_0$};
		
		\draw [blue, pattern = crosshatch dots, pattern color= blue] (1,-2) rectangle (2,-2.5);
		\node  at (4.5,-2.25)  {$\tilde{H}_-\cup \tilde{H}_1 \cup \tilde{H}_2$};
		
			\draw [red, pattern =north west lines
, pattern color= red] (-2.5,-2) rectangle (-1.5,-2.5);
		\node  at (-0.5,-2.25)  {$\tilde{H}_+$};
		
\end{tikzpicture}

\caption{The structure to be constructed in Step 1.}
\label{fig:bipswitchStep1}
\end{figure}
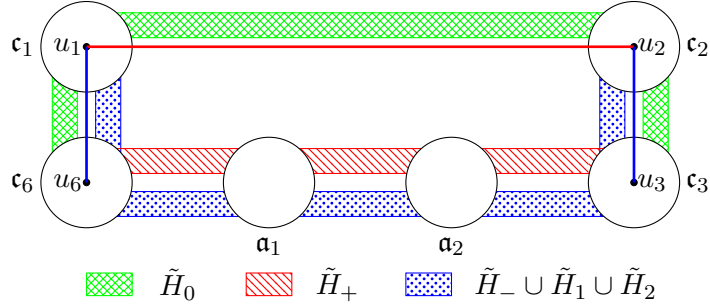

\begin{figure}[htbp]
\centering
\begin{tikzpicture}[scale=0.6]

				\draw [blue, pattern = crosshatch dots, pattern color= blue]  (-5.8,3) rectangle (-5.25,0);
				\draw [blue, pattern = crosshatch dots, pattern color= blue]  (5.8,3) rectangle (5.25,0);
				
				\draw [red, pattern =north west lines
, pattern color= red] (-5.8,-3) rectangle (-5.25,0);
				\draw [red, pattern =north west lines
, pattern color= red] (5.8,-3) rectangle (5.25,0);

				\draw [fill = white] (6,3) circle (1);
				\draw [fill = white] (-6,3) circle (1);
				\draw [fill = white] (6,-3) circle (1);
				\draw [fill = white] (-6,-3) circle (1);
				\draw [fill = white] (6,0) circle (1);
				\draw [fill = white] (-6,0) circle (1);
				\draw [fill = white] (2,0) circle (1);
				\draw [fill = white] (-2,0) circle (1);
				
				\filldraw[fill=black] (6,3) circle (2pt);
				\filldraw[fill=black] (-6,3) circle (2pt);
				\filldraw[fill=black] (6,-3) circle (2pt);
				\filldraw[fill=black] (-6,-3) circle (2pt);
				\filldraw[fill=black] (6,0) circle (2pt);
				\filldraw[fill=black] (-6,0) circle (2pt);
				
		\begin{scope}[red,line width = 1pt]
		\draw (6,3) -- (-6,3);
		\draw (6,-3) -- (6,0);
		\draw (-6,-3) -- (-6,0);
		\end{scope}
		
		\begin{scope}[blue,line width = 1pt]
		\draw (6,-3) -- (-6,-3);
		\draw (6,3) -- (6,0);
		\draw (-6,3) -- (-6,0);
		\end{scope}
			
		\node  at (-6.4,3)  {$u_1$};
		\node  at (6.4,3)  {$u_2$};
		\node  at (6.4,0)  {$u_3$};
		\node  at (6.4,-3)  {$u_4$};
		\node  at (-6.4,-3)  {$u_5$};
		\node  at (-6.4,0)  {$u_6$};
		
		\node  at (-2,-1.4)  {$\mathfrak{a}_1$};
		\node  at (2,-1.4)  {$\mathfrak{a}_2$};
		\node  at (-7.4,3)  {$\mathfrak{c}_1$};
		\node  at (7.4,3)  {$\mathfrak{c}_2$};
		\node  at (7.4,0)  {$\mathfrak{c}_3$};
		\node  at (7.4,-3)  {$\mathfrak{c}_4$};
		\node  at (-7.4,-3)  {$\mathfrak{c}_5$};
		\node  at (-7.4,0)  {$\mathfrak{c}_6$};
		\node  at (-4.4,1.5)  {$\tilde{H}_1$};
		\node  at (4.4,1.5)  {$\tilde{H}_2$};
		\node  at (-4.4,-1.5)  {$\tilde{H}_1'$};
		\node  at (4.4,-1.5)  {$\tilde{H}_2'$};
	\end{tikzpicture}
\caption{By mirroring the structure obtained in Step 1 (see Figure~\ref{fig:bipswitchStep1}), we obtain the desired $(u_1,\dots,u_6)_F$-switching, plus some unwanted (but symmetric) switchings. The subgraphs $\tilde{H}_+$ and $\tilde{H}_-$ are contained in both $F$-decompositions.}
\label{fig:bipswitchStep2}
\end{figure}
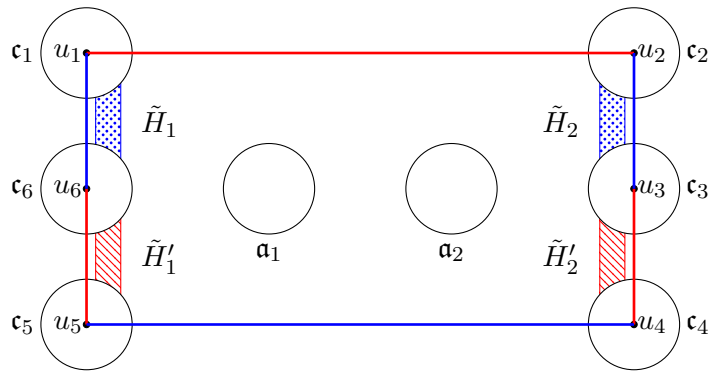

We will now construct the desired $C_6$-switcher. We will first translate the structure obtained from Lemma~\ref{lem:step1} into a `pseudo'-$(\Set{u_1u_2},\Set{u_1u_6,u_2u_3})_F$-switcher in the sense that there are some additional unwanted switchings (see Figure~\ref{fig:bipswitchStep1}). Here $u_1u_2$ plays the role of the edge $e^0$ from Lemma~\ref{lem:step1}. We will then mirror this structure so that the mirror image of $u_1u_2$ is $u_5u_4$ and $u_3,u_6$ are fixed points. The original structure together with its mirror image form a `pseudo'-$(u_1,\dots,u_6)_F$-switcher. In fact, some of the unwanted switchings will cancel out (see Figure~\ref{fig:bipswitchStep2}). By the inherent symmetry, we can pair up the unwanted switchings and decompose them in a number of `double-stars', which we can incorporate by adding suitable `double-star-switchers'.

\begin{lemma} \label{lem:bipC6switcher}
Let $F$ be bipartite and $\tau(F)=1$. Then there exists a $0$-compressible $(C_6)_F$-switcher.
\end{lemma}

\proof
Let $J$ be a copy of $C_6$ with vertices $\mathfrak{c}_1,\mathfrak{c}_2,\mathfrak{c}_3,\mathfrak{a}_2,\mathfrak{a}_1,\mathfrak{c}_6$ appearing in this order on the cycle.

\medskip
\noindent\emph{Step 1}
\medskip

Let $G_+,G_-$ be $F$-decomposable graphs, $u_1u_2\in E(G_+)$ and let $\rho\colon G_+\cup G_-\to J$ be a homomorphism such that $\rho(u_i)=\mathfrak{c}_i$ for $i\in\Set{1,2}$ and $G_-[\rho^{-1}(\mathfrak{c}_1\mathfrak{c}_2)]=G_+[\rho^{-1}(\mathfrak{c}_1\mathfrak{c}_2)]-\Set{u_1u_2}$, which exist by Lemma~\ref{lem:step1}. We may assume that $G_+$ and $G_-$ are otherwise edge-disjoint.

Let $H:=(G_+\cup G_-)-\Set{u_1u_2}$. For $\odot\in\Set{+,-}$, define $H_\odot:=G_\odot[\rho^{-1}(\mathfrak{c}_6\mathfrak{a}_1\mathfrak{a}_2\mathfrak{c}_3)]$, $H_{1,\odot}:=G_\odot[\rho^{-1}(\mathfrak{c}_1\mathfrak{c}_6)]$ and $H_{2,\odot}:=G_\odot[\rho^{-1}(\mathfrak{c}_2\mathfrak{c}_3)]$. Thereby, we have the following:

\begin{enumerate}[label=(H\arabic*)]
\item $\Set{u_1,u_2}\In V(H)$ is independent in $H$ and $\rho(u_i)=\mathfrak{c}_i$ for $i\in\Set{1,2}$; \label{1st step roots}
\item $H=H_0\cupdot H_+\cupdot H_- \cupdot H_{1,+} \cupdot H_{1,-} \cupdot H_{2,+} \cupdot H_{2,-}$, where
\begin{enumerate}[label=$\bullet$]
\item $H_0:=H[\rho^{-1}(\mathfrak{c}_1\mathfrak{c}_2)]$,
\item $H_+\cupdot H_-=H[\rho^{-1}(\mathfrak{c}_6\mathfrak{a}_1\mathfrak{a}_2\mathfrak{c}_3)]$,
\item $H_{1,+}\cupdot H_{1,-}=H[\rho^{-1}(\mathfrak{c}_1\mathfrak{c}_6)]$, and
\item $H_{2,+}\cupdot H_{2,-}=H[\rho^{-1}(\mathfrak{c}_2\mathfrak{c}_3)]$;
\end{enumerate} \label{1st step partition}
\item $H_0\cup \Set{u_1u_2} \cup H_+ \cup H_{1,+} \cup H_{2,+}$ and $H_0\cup H_- \cup H_{1,-} \cup H_{2,-}$ are both $F$-decomposable. \label{1st step decomposable}
\end{enumerate}

We now want to extend $H$ to $\tilde{H}$ and $\rho$ to $\tilde{\rho}$ such that the following hold (see Figure~\ref{fig:bipswitchStep1}):

\begin{enumerate}[label=(H\arabic*$'$)]
\item $\Set{u_1,u_2,u_3,u_6}\In V(\tilde{H})$ is independent in $\tilde{H}$ and $\tilde{\rho}(u_i)=\mathfrak{c}_i$ for $i\in\Set{1,2,3,6}$; \label{2nd step roots}
\item $\tilde{H}=\tilde{H}_0\cupdot \tilde{H}_+\cupdot \tilde{H}_- \cupdot \tilde{H}_{1} \cupdot \tilde{H}_{2}$, where
\begin{enumerate}[label=$\bullet$]
\item $\tilde{H}_+\cupdot \tilde{H}_-=\tilde{H}[\tilde{\rho}^{-1}(\mathfrak{c}_6\mathfrak{a}_1\mathfrak{a}_2\mathfrak{c}_3)]$,
\item $\tilde{H}_{1,0}\cupdot \tilde{H}_{1}=\tilde{H}[\tilde{\rho}^{-1}(\mathfrak{c}_1\mathfrak{c}_6)]$,
\item $\tilde{H}_{2,0}\cupdot \tilde{H}_{2}=\tilde{H}[\tilde{\rho}^{-1}(\mathfrak{c}_2\mathfrak{c}_3)]$, and
\item $\tilde{H}_0:=\tilde{H}[\tilde{\rho}^{-1}(\mathfrak{c}_1\mathfrak{c}_2)] \cup \tilde{H}_{1,0} \cup \tilde{H}_{2,0}$;
\end{enumerate} \label{2nd step partition}
\item $\tilde{H}_0\cup \Set{u_1u_2} \cup \tilde{H}_+$ and $\tilde{H}_0\cup \Set{u_1u_6,u_2u_3} \cup \tilde{H}_- \cup \tilde{H}_{1} \cup \tilde{H}_{2}$ are both $F$-decomposable. \label{2nd step decomposable}
\end{enumerate}
To this end, let $u_3,u_6$ be new vertices (i.e.~vertices not in $H$). Pick $vw\in E(F)$ and let $F_1$ be a copy of $F-vw$ such that $u_1,u_6$ play the roles of $v,w$, and let $F_2$ be a copy of $F-vw$ such that $u_2,u_3$ play the roles of $v,w$, and all other vertices are new vertices. Clearly, there exists a homomorphism $\rho_1\colon F_1 \rightarrow J[\Set{\mathfrak{c}_1,\mathfrak{c}_6}]$ such that $\rho_1(u_1)=\mathfrak{c}_1$ and $\rho_1(u_6)=\mathfrak{c}_6$. Similarly, let $\rho_2\colon F_2 \rightarrow J[\Set{\mathfrak{c}_2,\mathfrak{c}_3}]$ be a homomorphism such that $\rho_2(u_2)=\mathfrak{c}_2$ and $\rho_2(u_3)=\mathfrak{c}_3$.

For every $e\in E(H_{1,+} \cup H_{2,+})$, let $F_e$ be a copy of $F$ that contains $e$ and consists of new vertices apart from the endpoints of $e$. For $e\in E(H_{1,+})$, let $\rho_e\colon F_e \rightarrow J[\Set{\mathfrak{c}_1,\mathfrak{c}_6}]$ be a homomorphism such that $\rho_e(x)=\rho(x)$ for $x\in V(e)$. For $e\in E(H_{2,+})$, let $\rho_e\colon F_e \rightarrow J[\Set{\mathfrak{c}_2,\mathfrak{c}_3}]$ be a homomorphism such that $\rho_e(x)=\rho(x)$ for $x\in V(e)$.

Let $$\tilde{H}:=H\cup F_1\cup F_2 \cup \bigcup_{e\in E(H_{1,+} \cup H_{2,+})}F_e$$ and define $\tilde{\rho}:=\rho \cup \rho_1 \cup \rho_2 \cup \bigcup_{e\in E(H_{1,+} \cup H_{2,+})}\rho_e$. Then, \ref{2nd step roots} holds. Let $\tilde{H}_\odot:=H_\odot$ for $\odot\in\Set{+,-}$ and $\tilde{H}_{i,0}:=H_{i,+}$ for $i\in\Set{1,2}$. Note that $\tilde{H}[\tilde{\rho}^{-1}(\mathfrak{c}_1\mathfrak{c}_2)]=H[\rho^{-1}(\mathfrak{c}_1\mathfrak{c}_2)]=H_0$. Therefore, defining $\tilde{H}_0:=\tilde{H}[\tilde{\rho}^{-1}(\mathfrak{c}_1\mathfrak{c}_2)] \cup \tilde{H}_{1,0} \cup \tilde{H}_{2,0}$ yields $\tilde{H}_0=H_0\cup H_{1,+} \cup H_{2,+}$.
Finally, for $i\in\Set{1,2}$, let $$\tilde{H}_{i}:=H_{i,-}\cup F_i \cup \bigcup_{e\in E(H_{i,+})}(F_e-e).$$ This way, \ref{2nd step partition} holds. Now, $\tilde{H}_0\cup \Set{u_1u_2} \cup \tilde{H}_+$ is $F$-decomposable by \ref{1st step decomposable}. Moreover, $\tilde{H}_0\cup \Set{u_1u_6,u_2u_3} \cup \tilde{H}_- \cup \tilde{H}_{1} \cup \tilde{H}_{2}$ can be decomposed into $H_0\cup H_- \cup H_{1,-} \cup H_{2,-}$, $F_1 \cup \Set{u_1u_6}$, $F_2\cup \Set{u_2u_3}$ and $H_{1,+} \cup H_{2,+} \cup \bigcup_{e\in E(H_{1,+} \cup H_{2,+})}(F_e-e)$. The former is $F$-decomposable by \ref{1st step decomposable} and the others are all trivially $F$-decomposable, so \ref{2nd step decomposable} holds.

\medskip
\noindent\emph{Step 2}
\medskip

Let $J^\ast$ be the graph obtained from $J$ by mirroring $J$ with fixed points $\Set{\mathfrak{c}_6,\mathfrak{a}_1,\mathfrak{a}_2,\mathfrak{c}_3}$, that is, add new vertices $\mathfrak{c}_5,\mathfrak{c}_4$ to $J$ together with the edges $\mathfrak{c}_6\mathfrak{c}_5$, $\mathfrak{c}_5\mathfrak{c}_4$ and $\mathfrak{c}_4\mathfrak{c}_3$. 

Alternatively, $J^\ast$ can be viewed as the graph obtained from $C_6$ (with vertices $\mathfrak{c}_1,\dots,\mathfrak{c}_6$ in the usual order) by adding two new vertices $\mathfrak{a}_1,\mathfrak{a}_2$ and the edges $\mathfrak{c}_6\mathfrak{a}_1,\mathfrak{a}_1\mathfrak{a}_2,\mathfrak{a}_2\mathfrak{c}_3$. 

For every vertex $x\in \tilde{\rho}^{-1}(\mathfrak{c}_1\mathfrak{c}_2)$, let $x'$ be a new vertex. Let $\tilde{H}'$ be the copy of $\tilde{H}$ obtained by replacing every $x\in \tilde{\rho}^{-1}(\mathfrak{c}_1\mathfrak{c}_2)$ with $x'$, and let $\tilde{H}_0',\tilde{H}_{1}',\tilde{H}_{2}'$ be the subgraphs of $\tilde{H}'$ corresponding to $\tilde{H}_0,\tilde{H}_{1},\tilde{H}_{2}$. So $\tilde{H}'=\tilde{H}_0'\cupdot \tilde{H}_+\cupdot \tilde{H}_- \cupdot \tilde{H}_{1}' \cupdot \tilde{H}_{2}'$.
Let $$S^\ast:=\tilde{H}\cup \tilde{H}'.$$

We can extend $\tilde{\rho}$ to a homomorphism $\psi\colon S^\ast\rightarrow J^\ast$ by defining $\psi(x'):=\mathfrak{c}_5$ for every $x\in \tilde{\rho}^{-1}(\mathfrak{c}_1)$ and $\psi(x'):=\mathfrak{c}_4$ for every $x\in \tilde{\rho}^{-1}(\mathfrak{c}_2)$. 

We let $u_5:=u_1'$ and $u_4:=u_2'$, thus $u_5u_4$ is the mirror image of $u_1u_2$. Note that $U:=\Set{u_1,\dots,u_6}$ is independent in $S^\ast$. Let $f$ be defined as $f(u_i):=\mathfrak{c}_i$ for all $i\in[6]$. Therefore, $(J^\ast,f,J^\ast,\psi)$ is a $0$-compression of $(S^\ast,U)$.

Note that by symmetry and \ref{2nd step decomposable},
\begin{enumerate}[label=(H4$'$)]
\item $\tilde{H}_0'\cup \Set{u_4u_5} \cup \tilde{H}_+$ and $\tilde{H}_0'\cup \Set{u_5u_6,u_3u_4} \cup \tilde{H}_- \cup \tilde{H}_{1}' \cup \tilde{H}_{2}'$ are both $F$-decomposable.\label{strange}
\end{enumerate}
Let 
\begin{align}
S^{\ast\ast}:=\tilde{H}_0\cup \tilde{H}_+\cup \tilde{H}_- \cup \tilde{H}_0'=S^\ast-(\tilde{H}'_{1} \cup \tilde{H}'_{2}\cup \tilde{H}_{1} \cup \tilde{H}_{2}).\label{mirror}
\end{align}
Combining \ref{2nd step decomposable} and \ref{strange}, we conclude that
\begin{enumerate}[label=(H5$'$)]
\item $S^{\ast\ast}\cup \Set{u_1u_2,u_3u_4,u_5u_6} \cup \tilde{H}'_{1} \cup \tilde{H}'_{2}$ and $S^{\ast\ast}\cup\Set{u_2u_3,u_4u_5,u_6u_1} \cup \tilde{H}_{1} \cup \tilde{H}_{2}$ are both $F$-decomposable (see Figure~\ref{fig:bipswitchStep2}).\label{combined}
\end{enumerate}
Consider $x\in \tilde{\rho}^{-1}(\mathfrak{c}_1)$. Note that by \ref{2nd step partition}, $x$ is not incident with any edge from $\tilde{H}_+,\tilde{H}_-,\tilde{H}_{2},\Set{u_2u_3}$. Let $r:=gcd(F)$. By \ref{2nd step decomposable}, we have $r\mid d_{\tilde{H}_0\cup \Set{u_1u_2}}(x)$ and $r\mid d_{\tilde{H}_0\cup \Set{u_1u_6}\cup \tilde{H}_{1}}(x)$, implying that $r\mid d_{\tilde{H}_{1}}(x)$. Moreover, $N_{\tilde{H}_{1}}(x)=N_{\tilde{H}_{1}'}(x')$ since $N_{\tilde{H}_{1}}(x)\In \tilde{\rho}^{-1}(\mathfrak{c}_6)$. The same applies to $x\in \tilde{\rho}^{-1}(\mathfrak{c}_2)$ with $\tilde{H}_{2},\tilde{H}_{2}'$. Hence, by Corollary~\ref{cor:model attaching} and Proposition~\ref{prop:bip-star-switcher}, we can attach graphs $(S_x)_{x\in \tilde{\rho}^{-1}(\mathfrak{c}_1\mathfrak{c}_2)}$ to $S^\ast$ such that
\begin{itemize}
\item $S_x$ is an $(\tilde{H}_{i}[x,N_{\tilde{H}_{i}}(x)],\tilde{H}_{i}'[x',N_{\tilde{H}_{i}'}(x')])_F$-switcher for all $i\in\Set{1,2}$, $x\in \tilde{\rho}^{-1}(\mathfrak{c}_i)$;
\item $(S,U)$ is $0$-compressible with respect to $(J^\ast,f)$, where $S:=S^\ast\cup S_1 \cup S_2$ and $S_i:=\bigcup_{x\in \tilde{\rho}^{-1}(\mathfrak{c}_i)}S_x$ for $i\in\Set{1,2}$.
\end{itemize}
Clearly, $S_1$ is an $(\tilde{H}_{1},\tilde{H}_{1}')_F$-switcher and $S_2$ is an $(\tilde{H}_{2},\tilde{H}_{2}')_F$-switcher. Thus, together with \eqref{mirror} and \ref{combined}, it follows that $S$ is a $(u_1,\dots,u_6)_F$-switcher.

By identifying $\mathfrak{a}_1$ with $\mathfrak{c}_1$ and $\mathfrak{a}_2$ with $\mathfrak{c}_2$,
we can see that $(S,U)$ is $0$-compressible with respect to $(C_6,f)$ by Fact~\ref{fact:model shrink}.
\endproof

\begin{cor} \label{cor:bip abs 1}
Let $F$ be bipartite and $\tau(F)=1$. Then $F$ is $1/2$-absorbing.
\end{cor}

\proof
This follows from Lemma~\ref{lem:bipC6switcher}, Proposition~\ref{prop:bip-star-switcher}, Lemma~\ref{lem:switch2transform} and Lemma~\ref{lem:transf2abs}.
\endproof

\subsection{Lower bounds}

The remaining steps towards the proof of Theorem~\ref{thm:bipartite char} are extremal examples giving lower bounds on $\delta_F$.

\begin{prop} \label{prop:exex bip 1}
Let $F$ be bipartite. If $\tau(F)>1$, then $\delta_F\ge 2/3$.
\end{prop}

\proof
Let $r:=gcd(F)$. By Fact~\ref{fact:bip par rel}, $\tau(F)\mid r$. We show that there are $F$-divisible graphs $G$ of arbitrarily large order which are not $F$-decomposable, while $\delta(G) \geq \lfloor 2 |G| /3 \rfloor  -2 r (e(F)+1)$. 
Let $m$ be sufficiently large. 

\medskip
\noindent \textbf{Case 1:} \emph{$\tau(F)$ is odd.}
\medskip

Let $V_1,V_2, V_3$ be disjoint sets with $|V_1| = 2r m + \tau(F)-1$, $|V_2| = 2r m - \tau(F) + 2$, $|V_3| = 2r m - \tau(F)+1$.
Let $G$ be the graph on vertex set $V_1 \cup V_2 \cup  V_3$ consisting of two cliques on~$V_1$ and~$V_3$ and a complete bipartite graph with vertex classes $V_1 \cup V_3$ and $V_2$. 
Hence, $G$ has $6r m -\tau(F)+2$ vertices. Moreover, $d_G(v)=4rm$ for all $v\in V_1\cup V_2$ and $d_G(v)=4rm-2\tau(F)+2$ for all $v\in V_3$.

Let $G'$ be the graph obtained from $G$ by removing the edges of $r+1-\tau(F)$ edge-disjoint Hamilton cycles in $G[V_3]$. Hence, $G'$ is $r$-divisible and thus $r \mid 2e(G')$. Apply Proposition~\ref{prop:Walecki} to $G'[V_3]$ with $e(G')$ playing the role of $e$ to obtain an $r$-divisible subgraph $H$ of $G'[V_3]$ such that $e(G') \equiv e(H) \mod{e(F)}$ and $\Delta(H)\le 2e(F)r$. Let $G'':=G'-H$. Hence, $G''$ is $F$-divisible and $$\delta(G'') \ge 4rm -2\tau(F)+2-2(r+1-\tau(F))-2e(F)r \ge \lfloor 2|G''|/3 \rfloor - 2r(e(F)+1).$$

We will now see that $G''$ is not $F$-decomposable. Let $F'$ be any copy of $F$ in $G''$. Note that $V(F')\cap V_1$ is not $C_4$-supporting in $F'$, hence $\tau(F)\mid e(F'[V_1])$. So a necessary condition for $G''$ having an $F$-decomposition is that $e(G''[V_1])$ is divisible by $\tau(F)$. However, $e(G''[V_1]) = (r m + (\tau(F)-1)/2)(2 r m+\tau(F)-2)$. Recall that $\tau(F)\mid r$. Hence, the first factor is not divisible by $\tau(F)$ since $\tau(F)>1$ and the second factor is coprime to $\tau(F)$ as $\tau(F)$ is odd.

\medskip
\noindent \textbf{Case 2:} \emph{$\tau(F)$ is even.}
\medskip

Define $G$ as before, but this time the sizes of $V_1,V_2,V_3$ satisfy $|V_1| = 2r m + \tau(F)$, $|V_2| = 2r m-\tau(F)+1$, $|V_3| = 2r m - \tau(F)$.
Hence, $G$ has $6r m -\tau(F)+1$ vertices. Moreover, $d_G(v)=4rm$ for all $v\in V_1\cup V_2$ and $d_G(v)=4rm-2\tau(F)$ for all $v\in V_3$.

Let $G'$ be the graph obtained from $G$ by removing the edges of $r-\tau(F)$ edge-disjoint Hamilton cycles in $G[V_3]$. Hence, $G'$ is $r$-divisible and thus $r \mid 2e(G')$. Apply Proposition~\ref{prop:Walecki} to $G'[V_3]$ with $e(G')$ playing the role of $e$ to obtain an $r$-divisible subgraph $H$ of $G'[V_3]$ such that $e(G') \equiv e(H) \mod{e(F)}$ and $\Delta(H)\le 2e(F)r$. Let $G'':=G'-H$. Hence, $G''$ is $F$-divisible and $$\delta(G'') \ge 4rm -2\tau(F)-2(r-\tau(F))-2e(F)r \ge \lfloor 2|G''|/3 \rfloor - 2r(e(F)+1).$$

As before, a necessary condition for $G''$ having an $F$-decomposition is that $e(G''[V_1])$ is divisible by $\tau(F)$. However, $e(G''[V_1]) = (r m + \tau(F)/2)(2 r m+ \tau(F)-1)$, where the first factor is not divisible by $\tau(F)$ and the second factor is coprime to $\tau(F)$, so $e(G''[V_1])$ is not divisible by $\tau(F)$.
\endproof

\begin{prop} \label{prop:exex bip 2}
Let $F$ be bipartite. If $\tilde{\tau}(F)>1$ or every edge of $F$ is contained in a cycle, then $\delta_F\ge 1/2$.
\end{prop}

\proof
Let $r:=gcd(F)$. Suppose that $\tilde{\tau}(F)=1$. By Fact~\ref{fact:bip par rel}, $r=1$. Moreover, our assumption implies that every edge of $F$ is contained in a cycle. For any $m\in \bN$, the graph $G$ obtained from two disjoint cliques of order $me(F)$ each by deleting one edge and adding a bridge between the two components is $F$-divisible and satisfies $\delta(G)\ge |G|/2-2$, but is not $F$-decomposable.

We can therefore assume that $\tilde{\tau}(F)>1$.

\medskip
\noindent \textbf{Case 1:} \emph{$r$ is even or $\tilde{\tau}(F)>r$.}
\medskip

Let $a:=r$ if $r$ is odd and $a:=r/2$ if $r$ is even. Note that $a<\tilde{\tau}(F)$ since $r \mid \tilde{\tau}(F)$. For any $m\in \bN$, let $V_1,V_2$ be disjoint sets of size $2me(F)\tilde{\tau}(F)+1$ each and let $G$ be the graph consisting of two cliques on~$V_1$ and~$V_2$. Clearly, $G$ is $F$-divisible. Let $G'$ be the graph obtained from $G$ by removing the edges of $a$ edge-disjoint Hamilton cycles from $G[V_1]$ and $(e(F)-1)a$ edge-disjoint Hamilton cycles from $G[V_2]$. Observe that $G'$ is still $F$-divisible and $\delta(G')\ge 2m e(F)\tilde{\tau}(F) - 2e(F)r = |G'|/2-1-2e(F)r$.

However, $e(G'[V_1])\equiv -|V_1|a \mod{\tilde{\tau}(F)}$. Since $a<\tilde{\tau}(F)$ and $gcd\Set{|V_1|,\tilde{\tau}(F)}=1$, we deduce that $e(G'[V_1])$ is not divisible by $\tilde{\tau}(F)$, implying that $G'$ cannot be $F$-decomposable.

\medskip
\noindent \textbf{Case 2:} \emph{$r$ is odd and $\tilde{\tau}(F)=r$.}
\medskip

Since $\tilde{\tau}(F)>1$, we have $r>1$. We first claim that every edge of $F$ is contained in a cycle. Suppose that $xy$ is not contained in a cycle. Then there exists a partition of $V(F)$ into sets $A_1,A_2,B_1,B_2$ such that $x\in A_1$, $y\in B_2$ and $E(F)=E(F[A_1,B_1])\cup E(F[A_2,B_2])\cup\Set{xy}$. Hence, $$e(F)=e(F[A_1,B_1])+e(F[A_2,B_2])+1=\sum_{v\in B_1}d_F(v)+\sum_{v\in A_2}d_F(v)+1\equiv 1 \mod{r},$$ which contradicts $r \mid e(F)$ and $r>1$.

Let $Q$ be a graph with one vertex $q$ of degree $1$ and all other vertices of degree $r$. (To construct such a $Q$, start with $K_{r,r}$ and remove the edges of a matching of size $(r-1)/2$. Add a new vertex $q'$ and join $q'$ to all vertices that were incident with an edge from the matching. Add $q$ and join $q$ to $q'$.)

Now, for any sufficiently large $m\in \bN$, let $V_1,V_2$ be disjoint sets of size $rm+1$ each and let $G$ be the graph consisting of two cliques on~$V_1$ and~$V_2$. Clearly, $G$ is $r$-divisible. Fix $v_1\in V_1$ and $v_2\in V_2$. For $i\in\Set{1,2}$, let $Q_i$ be a copy of $Q$ in $G[V_i]$ such that $v_i$ plays the role of $q$. Let $G':=(G\cup\Set{v_1v_2})-Q_1-Q_2$. Clearly, $G'$ is $r$-divisible. In particular, $r\mid 2e(G')$. Apply Proposition~\ref{prop:Walecki} to $G'[V_1]$ with $e(G')$ playing the role of $e$ in order to obtain an $r$-divisible subgraph $H$ of $G'[V_1]$ such that $e(H)\equiv e(G') \mod{e(F)}$ and $\Delta(H)\le 2e(F)r$. Let $G'':=G'-H$. Thus, $G''$ is $F$-divisible and $\delta(G'')\ge rm-r-2e(F)r = |G''|/2-1-r-2e(F)r$. However, $G''$ is not $F$-decomposable because $v_1v_2$ cannot be covered.
\endproof

Let $F$ be a bipartite graph. In Section~\ref{sec:deltafa} we will see that $\delta_F^{vx}=0$ if $F$ contains a bridge and $\delta_F^{vx}=1/2$ otherwise (see Corollary~\ref{cor:deltafa}(ii)). Using this, we can now prove Theorem~\ref{thm:bipartite char}.

\lateproof{Theorem~\ref{thm:bipartite char}}
Note that $\delta_F^{vx}\le 1/2$ by Corollary~\ref{cor:deltafa}(ii) and $\delta_F^{0+}=0$ by Fact~\ref{fact:bip approx}. By Corollary~\ref{cor:absorbing}, $F$ is $2/3$-absorbing. Hence, by Theorem~\ref{thm:general dec}, $\delta_F\le 2/3$. Now, if $\tau(F)>1$, then Proposition~\ref{prop:exex bip 1} implies that $\delta_F=2/3$. On the other hand, if $\tau(F)=1$, then we can deduce from Corollary~\ref{cor:bip abs 1} and Theorem~\ref{thm:general dec} that $\delta_F\le 1/2$. If $\tilde{\tau}(F)>1$ or every edge is contained in a cycle, we deduce $\delta_F=1/2$ with Proposition~\ref{prop:exex bip 2}. So assume that $\tilde{\tau}(F)=1$ and that $F$ contains a bridge. Then $\delta_F^{vx}= 0$ by Corollary~\ref{cor:deltafa}(ii) and $F$ is $0$-absorbing by Lemma~\ref{lem:bip abs 2}. Hence, Theorem~\ref{thm:general dec} implies that $\delta_F=0$.
\endproof

\section{Covering the edges at a vertex} \label{sec:deltafa}

In this section, we investigate $\delta_F^{vx}$, i.e.~the threshold at which we can cover all edges at a vertex. In particular, we will determine $\delta_F^{vx}$ for all bipartite graphs. In the general case, we determine $\delta_F^{vx}$ as a function of $\delta_F^e$, that is, we reduce the problem of covering all edges at one vertex to the problem of covering one edge.
We will use an iterative absorbing approach which has many parallels to the main proof.

In Section~\ref{subsec:approx}, we will show how to obtain an approximate cover at some vertex $x$. In Section~\ref{subsec:almost}, we will show how to turn an approximate cover into a near-optimal cover. Roughly speaking, the neighbourhood of $x$ will be partitioned into sets $N_1,\dots,N_\ell$ of successively smaller size, where $N_\ell$ has constant size. Using the result of Section~\ref{subsec:approx}, we can cover all but a small fraction of the edges from $x$ to $N_1$. By assuming that the minimum degree is above $\delta_F^e$, we can cover the leftover edges one by one by using some edges from $x$ to $N_2$. We then cover all but a small fraction of the remaining edges from $x$ to $N_2$ and continue as above until all edges at $x$ are covered except some from $x$ to $N_\ell$. In Section~\ref{subsec:absorb}, we will see how these remaining edges can be absorbed.

\subsection{Approximate cover} \label{subsec:approx}
In order to determine the threshold which guarantees an approximate cover with copies of $F$ at a vertex, we will use a result of Koml\'os~\cite{Ko}. He showed that the minimum degree threshold that guarantees the existence of vertex-disjoint copies of a given graph $F$ covering almost all vertices of the host graph is governed by the so-called critical chromatic number of $F$. We will apply his result to a reduced graph.

For a graph $F$, let $Col(F)$ denote the set of all $[\chi(F)]$-colourings of $F$ and let $\sigma(F):=\min_{c\in Col(F)}|c^{-1}(1)|$. The \emph{critical chromatic number} of $F$ is defined as $$\chi_{cr}(F):=(\chi(F)-1)\frac{|F|}{|F|-\sigma(F)},$$ where we set $\chi_{cr}(F):=0$ if $\chi(F)=1$.

\begin{theorem}[Koml\'os~\cite{Ko}] \label{thm:Komlos}
For every graph $H$ and $\mu >0$ there exists an $n_0 = n_0(\mu,H)$ such that every graph $G$ on $n \geq n_0$ vertices with $\delta(G) \ge (1-1/\chi_{cr}(H)) n $ contains vertex-disjoint copies of $H$ covering all but at most $\mu n$ vertices of $G$.
\end{theorem}
(Theorem~\ref{thm:Komlos} was further improved by Shokoufandeh and Zhao~\cite{SZ} who replaced the $\mu n$ term with a constant depending only on $H$. However, the result of Koml\'os is sufficient for our purposes.) We will apply Theorem~\ref{thm:Komlos} to an appropriate subgraph of the reduced graph.

We now define the graph parameter $\chi^{vx}$ that governs the existence of an approximate cover at one vertex, which turns out to be closely related to $\chi_{cr}$.

For a graph $F$ with $\chi:=\chi(F)$ and a vertex $v\in V(F)$, let $Col(F,v):=\set{c\in Col(F)}{c(v)=\chi}$ and define $$\sigma(F,v):=\min_{c\in Col(F,v)}|N_F(v)\cap c^{-1}(1)|.$$
So if $F$ can be coloured in a way such that $N_F(v)$ requires fewer than $\chi-1$ colours, then $\sigma(F,v)=0$. Note that if $F$ is bipartite, then $\sigma(F,v)=d_F(v)$ for all $v\in V(F)$. But if $\chi\ge 3$, then $\sigma(F,v)<d_F(v)$ for all $v\in V(F)$.
Thus, if $\chi\ge 3$, then $$\chi^{vx}(F):=(\chi-2)\min_{v\in V(F)}\frac{d_F(v)}{d_F(v)-\sigma(F,v)}$$ is well-defined. Note that for all $v\in V(F)$,
\begin{align}
\sigma(F,v)\ge d_F(v)(1-(\chi-2)/\chi^{vx}(F)).\label{eq:chi_vx degrees}
\end{align}
Moreover, we set $\chi^{vx}(F):=0$ if $F$ is bipartite. Clearly,
\begin{align}
\chi-2 \le \chi^{vx}(F) \le \chi-1.\label{eq:chi_vx relation}
\end{align} 
\COMMENT{We don't have $<$ for the left inequality, even not if $\chi>2$.}

\begin{prop} \label{prop:space extremal}
For all graphs $F$, $\min\Set{\delta_F^{vx},\delta_F^{0+}} \ge 1-1/(\chi^{vx}(F)+1)$.
\end{prop}

We will see that $\delta_F^{vx}$ can be strictly larger than $1-1/(\chi^{vx}(F)+1)$. Roughly speaking, $1-1/(\chi^{vx}(F)+1)$ represents only a `space' barrier, whereas $\delta_F^{vx}$ is also subject to other kinds of barriers. The reason why we also show that $\delta_F^{0+}\ge 1-1/(\chi^{vx}(F)+1)$ is that because of this inequality, we can omit the term $1-1/(\chi^{vx}(F)+1)$ in the discretisation given in Theorem~\ref{thm:main}.

\proof
Let $\chi:=\chi(F)$, $\chi^{vx}:=\chi^{vx}(F)$ and $\delta:=\min\Set{\delta_F^{vx},\delta_F^{0+}}$. Clearly, $\delta \ge 1-1/(\chi-1)$. We may therefore assume that $\chi^{vx}>\chi-2$. In particular, $\chi\ge 3$. Let $\mu:=1-1/(\chi^{vx}+1) - \delta$ and suppose, for a contradiction, that $\mu>0$. We will construct graphs $G$ of arbitrarily large order with $\delta(G)\ge(\delta+\mu/2)|G|$ such that any $F$-collection in $G$ fails to cover $\Omega(|G|^2)$ edges and that there is a vertex $x\in V(G)$ with $gcd(F)\mid d_G(x)$ such that the edges at $x$ cannot be covered with edge-disjoint copies of $F$, contradicting the definition of $\delta$.

More precisely, choose $\nu>0$ small enough such that $1-1/(\chi^{vx}+1-\nu)\ge 1-1/(\chi^{vx}+1)-\mu/2$ and $\alpha:=\chi^{vx}-(\chi-2)-\nu>0$. Let $m\in \bN$ be sufficiently large. We may assume that $(\chi-2)m+\alpha m$ is an integer and divisible by $gcd(F)$. Let $G$ be the complete $\chi$-partite graph with vertex classes $V_1,\dots,V_\chi$ such that $|V_i|=m$ for $i\in[\chi-1]$ and $|V_\chi|=\alpha m$. Thus, $|G|=(\chi-1+\alpha)m$ and $d_G(x)\ge (\chi-2+\alpha)m$ for all $x\in V(G)$ since $\alpha< 1$. By our choice of $\nu$, it follows that $\delta(G)\ge (\delta+\mu/2)|G|$.\COMMENT{$$\frac{\chi-2+\alpha}{\chi-1+\alpha}=1-\frac{1}{\chi-1+\alpha}=1-\frac{1}{\chi^{vx}+1-\nu}\ge 1-1/(\chi^{vx}+1)-\mu/2$$}

Observe that $$\eta:=1-\frac{\alpha}{\alpha+\nu\frac{\chi-2}{\chi^{vx}}}>0.$$

Let $x$ be any vertex in $V(G)\sm V_{\chi}$. By our assumption, $gcd(F)\mid d_G(x)$.
We will now see that any $F$-collection in $G$ fails to cover at least $\eta d_G(x)$ edges at $x$. Let $F_1,\dots,F_t$ be edge-disjoint copies of $F$ in $G$, all containing $x$, and let $\tilde{d}:=\sum_{j=1}^td_{F_j}(x)$ denote the number of edges covered at $x$.
For $j\in[t]$, let $v_j$ be the vertex of $F$ whose role $x$ is playing in $F_j$. Thus, for all $j\in[t]$, $d_{F_j}(x,V_\chi)\ge \sigma(F,v_j)$. Hence, $\sigma(F,v_1)+\dots+\sigma(F,v_t) \le \alpha m$, and \eqref{eq:chi_vx degrees} thus implies that $\tilde{d}\cdot(1-(\chi-2)/\chi^{vx}) \le \alpha m$.
Since $d_G(x)=(\chi-2+\alpha )m = (\chi^{vx}-\nu)m$, we deduce that $$\frac{\tilde{d}}{d_G(x)}\le \frac{\alpha m}{(1-\frac{\chi-2}{\chi^{vx}})(\chi^{vx}-\nu)m}=\frac{\alpha}{\alpha+\nu\frac{\chi-2}{\chi^{vx}}}=1-\eta ,$$ proving the claim.

In particular, there exists a vertex $x\in V(G)$ with $gcd(F)\mid d_G(x)$ such that the edges at $x$ cannot be covered with edge-disjoint copies of $F$, implying that $\delta^{vx}>\delta$.

Moreover, let $\cF$ be any $F$-collection in $G$ and let $H$ be the subgraph of $G$ consisting of the edges that are not covered by $\cF$. By the above, we have $d_H(x)\ge \eta d_G(x)$ for all $x\in V(G)\sm V_\chi$. Thus, $2e(H)\ge \eta(\chi-2+\alpha)m\cdot (\chi-1)m \ge \eta (\chi-2)^2m^2$. Since $|G|\le \chi m$, we deduce that $e(H)\ge \eta(1-2/\chi)^2|G|^2/2 \ge \eta|G|^2/18$. Thus, $\delta_F^{0+}\ge \delta_F^{\eta/18}>\delta$, contradicting the definition of $\delta$.
\endproof

We will now show how to obtain an approximate cover of the edges at a specified vertex~$x$. 
The following lemma is an easy application of the key lemma. It describes a structure within which we can cover almost all edges at $x$ with copies of $F$.

\begin{lemma} \label{lem:key lemma application}
Let $F$ be a $\chi$-chromatic graph and $u\in V(F)$. Let $c\in Col(F,u)$ be such that $a_1,\dots,a_s\ge 1$ and $a_{s+1},\dots,a_{\chi-1}=0$ for some $s\in[\chi-1]$, where $a_i:=|N_F(u)\cap c^{-1}(i)|$ for $i\in[\chi-1]$. 

Let $1/n\ll \eps \ll \alpha,1/|F|$ and suppose that $G$ is a graph with
\begin{enumerate}[label=(\roman*)]
\item $V(G)=\Set{x}\cupdot V_1 \cupdot \cdots \cupdot V_{\chi}$;
\item $N_G(x)=V_1\cup \dots \cup V_s$;
\item $G[V_i,V_j]$ is $\eps$-regular with density at least $\alpha$ for all $1\le i<j \le \chi$;
\item $|V_i|=a_in$ for all $i\in[s]$;
\item $|V_i|=n$ for all $s<i\le \chi$.
\end{enumerate}
Then there exists an $F$-collection covering all but at most $\sqrt{\eps} d_G(x)$ edges at $x$.
\end{lemma}

\proof
Let $t:=\lceil(1-\sqrt{\eps})d_G(x)/d_F(u)\rceil$. Note that $d_G(x)=(a_1+\dots+a_s)n=d_F(u)n$. We will greedily find injective homomorphisms $\phi_1,\dots,\phi_t$ from $F$ into $G$ such that $\phi_1(F),\dots,\phi_t(F)$ are edge-disjoint, $\phi_j(u)=x$ for all $j\in[t]$ and $\phi_j(v)\in V_{c(v)}$ for all $j\in[t]$ and $v\in V(F)\sm\Set{u}$. Suppose that for some $j\in [t]$, we have already found $\phi_1,\dots,\phi_{j-1}$. We now want to find $\phi_j$. Let $H:=G-\phi_1(F)-\dots-\phi_{j-1}(F)$. Note that for every $i\in[s]$, we have $$d_H(x,V_i)=|V_i|-(j-1)a_i\ge |V_i|-(1-\sqrt{\eps})a_id_G(x)/d_F(u)=\sqrt{\eps}|V_i|\ge \sqrt{\eps} n.$$ Hence, for every $i\in[s]$, we can pick a set $V_i'\In N_H(x,V_i)$ of size $\sqrt{\eps} n$. Furthermore, for every $s<i\le \chi$, we can pick a set $V_i'\In V_i$ of size $\sqrt{\eps} n$. By Fact~\ref{fact:slicing} and since $e(G-H)\le ne(F)$, we know that $H[V_i',V_j']$ is ${\eps}^{1/3}$-regular with density at least $\alpha/2$ for all $1\le i<j\le \chi$. Therefore, the key lemma (Lemma~\ref{lem:key lemma}) implies that there exists an embedding $\phi_j'\colon (F-u)\rightarrow H$ such that $\phi_j'(v)\in V'_{c(v)}$ for all $v\in V(F)\sm \Set{u}$. Defining $\phi_j(u):=x$ and $\phi_j(v):=\phi'_j(v)$ for all $v\in V(F)\sm \Set{u}$ yields the desired embedding~$\phi_j$.
\endproof

\begin{lemma} \label{lem:approximate cover}
Let $F$ be a graph and $1/n \ll \mu,1/|F|$. Let $G$ be a graph on $n$ vertices with $\delta(G)\ge (1-1/(\chi^{vx}(F)+1)+\mu)n$ and $x\in V(G)$. Then there exists an $F$-collection covering all but at most $\mu n$ edges at $x$.
\end{lemma}

To prove Lemma~\ref{lem:approximate cover}, we will apply the regularity lemma and consider the corresponding reduced graph to find the structures described in Lemma~\ref{lem:key lemma application}. Applying Lemma~\ref{lem:key lemma application} will then give the desired $F$-collection.

\proof
Let $u$ be a vertex of $F$ such that $\chi^{vx}:=\chi^{vx}(F)=(\chi-2)d_F(u)/(d_F(u)-\sigma(F,u))$. By definition, there exists $c\in Col(F,u)$ such that $\sigma(F,u)=a_{\chi-1}$ and $a_1,\dots,a_s>0$ and $a_{s+1},\dots,a_{\chi-1}=0$ for some $s\in[\chi-1]$, where $a_i:=|N_F(u)\cap c^{-1}(i)|$ for $i\in[\chi-1]$. 

Let $H$ be the complete $s$-partite graph with class sizes $a_1,\dots,a_s$ and observe that
\begin{align}
\chi_{cr}(H)\le \chi^{vx}.\label{critical and vx}
\end{align}\COMMENT{If $s\le \chi-2$, then $\chi_{cr}(H)\le \chi-2\le \chi^{vx}$. If $s=\chi-1$, then $\chi(H)-1= \chi-2$, $|H|=d_F(u)$, $\sigma(H)= \sigma(F,u)$ and hence $\chi_{cr}(H)=\chi^{vx}$.} Let $A:=\prod_{i\in[s]}a_i$.

Choose new constants $k_0,k_0'\in \bN$ and $\eps,\alpha>0$ such that $1/n\ll 1/k_0' \ll \eps \ll \alpha , 1/k_0 \ll \mu,1/|F|$. (Then also $\alpha,1/k_0\ll 1/A$.) Suppose that $G,n,x$ are as in the hypothesis. Apply the regularity lemma (Lemma~\ref{lem:regularity}) to obtain a partition $V_0,V_1,\dots,V_k$ of $V(G)-x$ and a spanning subgraph $G'$ of $G-x$ satisfying the following:
\begin{enumerate}[label=(R\arabic*)]
\item $k_0 \le k \le k_0'$;
\item $|V_0| \le \eps n$;
\item $|V_1|=\dots=|V_k|=:L$;
\item $d_{G'}(z) \ge d_G(z)-(\alpha+\eps)n$ for all $z \in V(G)\sm\Set{x}$;
\item $G'[V_y]$ is empty for all $y\in [k]$;
\item for all $1\le y < y' \le k$, $G'[V_y,V_{y'}]$ is either $\eps$-regular with density at least $\alpha$ or empty;
\item for all $y\in [k]$, $V_y\In N_G(x)$ or $V_y\cap N_G(x)=\emptyset$.
\end{enumerate}
By adding some vertices of $V(G')\sm V_0$ to $V_0$ if necessary, we can assume that $L$ is divisible by $A$. Let $X\In [k]$ be the set of indices $y$ for which $V_y\In N_G(x)$. Let $R$ be the reduced graph of $V_1,\dots,V_k$ with respect to $G'$, and let $R_x:=R[X]$. By Proposition~\ref{prop:reduced graph}, $\delta(R)\ge (1-1/(\chi^{vx}+1)+\mu/2)k$. Firstly, this implies $\delta(R)\ge (1-1/(\chi-1)+\mu/2)k$. Secondly, since $|X|L\ge d_G(x)-\eps n$ and thus $|X|\ge (1-1/(\chi^{vx}+1))k$, we have that $\delta(R_x)\ge (1-1/\chi^{vx})|X|\ge (1-1/\chi_{cr}(H))|R_x|$ by \eqref{critical and vx}. Hence, by Theorem~\ref{thm:Komlos}, there exist vertex-disjoint copies $H_1,\dots,H_t$ of $H$ in $R_x$ covering all but at most $\mu |R_x|/2$ vertices.

Ideally we would now like to extend each $H_i$ to a copy of $F-u$ in $R$ so that these copies are edge-disjoint, and then try to apply the key lemma to the corresponding subgraphs of $G$. However, the number of vertices in $H$ is too large for this to work in general. Instead, we will extend every $s$-clique in $H_i$ into a $\chi$-clique and then construct (in $G$) structures appropriate for Lemma~\ref{lem:key lemma application}.

For each $j\in[t]$, let $H_{j,1},\dots,H_{j,s}$ be the vertex classes of $H_j$ with $|H_{j,i}|=a_i$ for all $i\in[s]$. Let $$\mathcal{A}:=\set{\xi=(y_1,\dots,y_s)}{\exists j\in[t] \mbox{ so that } y_i\in H_{j,i}\mbox{ for all }i\in[s]}.$$
If $\xi=(y_1,\dots,y_s)$, we write $\xi^\ast$ for $\Set{y_1,\dots,y_s}$. So for every $\xi\in\mathcal{A}$, $R[\xi^\ast]$ is a clique. We want to extend each such clique to a clique on $\chi$ vertices by attaching additional vertices. Moreover, we want that all attachments are edge-disjoint. More precisely, for every $\xi\in \mathcal{A}$, we want to find a set $att(\xi)$ of $\chi-s$ vertices in $R$ such that
\begin{enumerate}[label=(A\arabic*)]
\item $R[\xi^\ast\cup att(\xi)]$ is a clique on $\chi$ vertices;\label{attach clique}
\item the graphs $R[att(\xi)]\cup R[att(\xi),\xi^\ast]$ are all pairwise edge-disjoint and edge-disjoint from $\bigcup_{j\in [t]}H_j$.\label{attachments edge-disjoint}
\end{enumerate}
We will achieve this using Lemma~\ref{lem:finding}. First, let $R^\ast:=R-\bigcup_{j\in [t]}H_j$ and note that $\delta(R^\ast)\ge(1-1/(\chi-1)+\mu/4)k$. Moreover, let $K^\ast:=K_{\chi}-K_{\chi}[S]$, where $S\In V(K_\chi)$ is of size $s$. Clearly, $(K^\ast,S)$ is a model and $K^\ast$ has degeneracy at most $\chi-1$ rooted at $S$. For every $\xi \in \mathcal{A}$, let $\Lambda_\xi$ be any bijection from $S$ to $\set{\Set{y}}{y\in\xi^\ast}$. Hence, $\Lambda_\xi$ is an $R^\ast$-labelling of $S$. Now, since $|\mathcal{A}|\le At\le Ak$ and for every $y\in V(R)$, $|\set{\xi\in\mathcal{A}}{y\in \xi^\ast}|\le A$, we can apply Lemma~\ref{lem:finding} in order to obtain edge-disjoint embeddings $(\phi_\xi)_{\xi\in\mathcal{A}}$ of $K^\ast$ into $R^\ast$ such that $\phi_\xi$ respects $\Lambda_\xi$. Thus, we can take $att(\xi):=\phi_\xi(V(K_\chi)\sm S)$. This satisfies \ref{attach clique} and \ref{attachments edge-disjoint}.

We note that we will only need some of the sets $att(\xi)$. For reasons that will become clear later, we partition $V_y$ for every $y\in V(R)$ arbitrarily into $A$ equal-sized parts $V_{y,1},\dots,V_{y,A}$.

Fix $j\in[t]$. For each $i\in[s]$ and $y\in H_{j,i}$, we can cut $V_y$ into $A/a_i$ equal-sized parts in order to obtain a partition of $\bigcup_{y\in H_{j,i}}V_y$ into $A$ parts $W_{j,i,1},\dots,W_{j,i,A}$ of size $La_i/A$ each. Fix $\ell\in[A]$. For $i\in[s]$, let $y_{j,i,\ell}$ be the vertex of $H_{j,i}$ such that $W_{j,i,\ell}$ is contained in $V_{y_{j,i,\ell}}$. We have $\xi_{j,\ell}:=(y_{j,1,\ell},\dots,y_{j,s,\ell})\in \mathcal{A}$. Let $y_{j,s+1,\ell},\dots,y_{j,\chi,\ell}$ be the elements of $att(\xi_{j,\ell})$ and then define $W_{j,i,\ell}:=V_{y,\ell}$ with $y=y_{j,i,\ell}$ for each $s<i\le \chi$.
Observe that by \ref{attachments edge-disjoint}, for all $\ell,\ell'\in[A]$,
\begin{itemize}
\item[(A3)] $R[\xi^\ast_{j,\ell}\cup att(\xi_{j,\ell})]$ and $R[\xi^\ast_{j',\ell'}\cup att(\xi_{j',\ell'})]$ are edge-disjoint if $j\neq j'$.
\end{itemize}
(But it may happen that $\xi^\ast_{j,\ell}\cup att(\xi_{j,\ell})=\xi^\ast_{j,\ell'}\cup att(\xi_{j,\ell'})$ even if $\ell\neq\ell'$.)

Now, for every $(j,\ell)\in[t]\times [A]$, let $G_{j,\ell}$ be the graph obtained from $G'[\bigcup_{i\in[\chi]}W_{j,i,\ell}]$ by adding $x$ and all edges from $x$ to $\bigcup_{i\in[s]}W_{j,i,\ell}$. As we shall see, the graph $G_{j,\ell}$ is the desired structure to which we can apply Lemma~\ref{lem:key lemma application}.

We claim the following:
\begin{enumerate}[label=(\roman*)]
\item $(G_{j,\ell})_{j\in[t],\ell\in[A]}$ is a family of edge-disjoint subgraphs of $G$;\label{Gjl:edge-disjoint}
\item the graphs $G_{j,\ell}$ cover all but at most $(\eps+\mu/2)n$ edges at $x$ in $G$;\label{Gjl:xcover}
\item $G_{j,\ell}[W_{j,i,\ell},W_{j,i',\ell}]$ is $2A\eps$-regular with density at least $\alpha/2$ for all $(j,\ell)\in[t]\times [A]$ and $1\le i<i'\le \chi$;\label{Gjl:regular}
\item $|W_{j,i,\ell}|=La_i/A$ for all $(j,i,\ell)\in[t]\times [s]\times [A]$;\label{Gjl:weighted sizes}
\item $|W_{j,i,\ell}|=L/A$ for all $(j,\ell)\in[t]\times [A]$ and $s<i\le \chi$.\label{Gjl:additional sizes}
\end{enumerate}

Firstly, whenever $(j,i,\ell),(j',i',\ell')\in[t]\times [s]\times [A]$ are distinct, then $W_{j,i,\ell}\cap W_{j',i',\ell'}=\emptyset$. Hence, every edge at $x$ is contained in at most one of the $G_{j,\ell}$. Moreover, $xz$ is covered whenever $z\in V_{y}$ with $y\in \bigcup_{j\in[t]}V(H_j)$, implying \ref{Gjl:xcover}.

We continue with proving \ref{Gjl:edge-disjoint}. If $j\neq j'$, then $G_{j,\ell}$ and $G_{j',\ell'}$ are edge-disjoint by (A3). It remains to check that $G_{j,\ell}$ and $G_{j,\ell'}$ are edge-disjoint for fixed $j\in[t]$ and distinct $\ell,\ell'\in[A]$. In fact, they are vertex-disjoint (except for $x$ of course). Clearly, $W_{j,i,\ell}\cap W_{j,i',\ell'}=\emptyset$ whenever $i,i'\in [s]$. Moreover, $y_{j,i,\ell}\notin V(H_j)$ for $i>s$ since $R[att(\xi_{j,\ell}),\xi^\ast_{j,\ell}]$ is complete bipartite but also edge-disjoint from $H_j$. Hence, $W_{j,i,\ell}\cap W_{j,i',\ell'}=\emptyset$ whenever $i\notin [s]$ and $i'\in[s]$ (and vice versa). Finally, for $i,i'\notin[s]$, $W_{j,i,\ell}\cap W_{j,i',\ell'}=\emptyset$ by our partition of each $V_y$ into $A$ parts.

The size conditions \ref{Gjl:weighted sizes} and \ref{Gjl:additional sizes} follow directly from the definitions. Moreover, \ref{Gjl:regular} holds by \ref{attach clique} and by Fact~\ref{fact:slicing}.

Hence, we can apply Lemma~\ref{lem:key lemma application} to each $G_{j,\ell}$ in order to find edge-disjoint copies of $F$ covering all but at most $\sqrt{2A\eps}d_{G_{j,\ell}}(x)$ edges at $x$ in $G_{j,\ell}$. Thus, by~\ref{Gjl:xcover}, all but at most $(\eps+\mu/2+\sqrt{2A\eps})n\le \mu n$ edges at $x$ are covered, as desired.
\endproof

\subsection{Covering all edges} \label{subsec:almost}

The following lemma is an analogue to Lemma~\ref{lem:near optimal}. It guarantees a `near-optimal' cover of the edges at $x$.

\begin{lemma} \label{lem:almost cover}
Let $F$ be a graph and $\delta:=\max\Set{1-1/(\chi^{vx}(F)+1),\delta_F^e}$. Assume that $1/m \ll \mu,1/|F|$. Let $G$ be a graph with $\delta(G) \ge (\delta + 2\mu)|G|$ and let $U_1\In V(G)$ be such that $|U_1|=\lfloor \mu |G|\rfloor$ and $d_G(y,U_1)\ge (\delta+3\mu)|U_1|$ for all $y\in V(G)$. Suppose that $U_1 \supseteq U_2 \supseteq \dots \supseteq U_\ell$ is a $(\delta+4\mu,\mu,m)$-vortex in $G[U_1]$ and $x\in U_\ell$. Then there exist edge-disjoint copies of $F$ covering all edges at $x$ except possibly some edges from $x$ to $U_\ell$.
\end{lemma}

\proof
Let $\gamma>0$ be such that $1/m\ll \gamma\ll \mu,1/|F|$. We proceed by induction on $\ell$. For $\ell=0$, there is nothing to prove. So assume that $\ell>0$ and that the statement is true for $\ell-1$. 

Let $R:=U_1\sm\Set{x}$, $L:=(V(G)\sm U_1) \cup \Set{x}$ and $G':=G[L]$. Note that $\delta(G')\ge (\delta+\gamma)|G'|$. By Lemma~\ref{lem:approximate cover}, there exists an $F$-collection $\cF_1$ in $G'$ covering all but at most $\gamma |G'|$ edges at $x$. 

Let $H$ be the subgraph of $G$ consisting of all those edges from $x$ to $V(G)\sm U_1$ which are not covered by $\cF_1$. Let $G'':=G-G[U_2]$ if $\ell\ge 2$ and $G'':=G$ otherwise. So $\Delta(H)=d_H(x)\le \gamma |G'| \le \gamma |G''|$. Since $d_{G''}(y,R) \ge (\delta+\mu/2)|R|$ for every $y\in V(G'')$, we can apply Proposition~\ref{prop:cover-sparse-graph} in order to obtain a subgraph $A$ of $G''$ such that $A[L]$ is empty, $H\cup A$ has an $F$-decomposition $\cF_2$ and $\Delta(A[R])\le \mu^2|R|/4$. By deleting copies of $F$ from $\cF_2$ which do not contain any edge of $H$, we can assume that $d_A(x)\le |F|\gamma |G|$. Hence, $\Delta(A[U_1])\le \max\Set{|F|\gamma |G|,\mu^2|R|/4+1} \le \mu^2 |U_1|/2$.

Note that $\bigcup \cF_1$ and $\bigcup \cF_2$ are edge-disjoint and together cover all edges from $x$ to $V(G)\sm U_1$. If $\ell=1$, this completes the proof.

If $\ell\ge 2$, then let $G''':=G[U_1] - \bigcup \cF_1 - \bigcup \cF_2= G[U_1]-A[U_1]$. Clearly, $G'''[U_2]=G[U_2]$, and so $U_2 \supseteq U_3 \supseteq \dots \supseteq U_\ell$ is a $(\delta+4\mu,\mu,m)$-vortex in $G'''[U_2]$. Let $y\in U_1$. Then $d_{G'''}(y)\ge d_G(y,U_1)-\mu^2 |U_1|/2 \ge (\delta+2\mu)|U_1|$ and $d_{G'''}(y,U_2)\ge d_G(y,U_2)-\mu^2 |U_1|/2 \ge (\delta+3\mu)|U_2|$. 

By induction, there exists an $F$-collection $\cF_3$ covering all edges at $x$ in $G'''$ except possibly some edges from $x$ to $U_\ell$. Finally, $\cF_1\cup \cF_2 \cup \cF_3$ covers all edges at $x$ in $G$ except possibly some edges from $x$ to $U_\ell$.
\endproof

We will now prove an analogue of Theorem~\ref{thm:general dec}, that is, assuming that we are able to absorb a small number of leftover edges, we can cover all edges at one vertex. We use the concept of absorption in the following form.

Given a graph $F$, a vertex $x$ and a vertex set $W$, an \defn{$F$-neighbourhood-absorber for $(x,W)$} is a graph $A$ such that
\begin{enumerate}[label=$\bullet$]
\item $W\cup \Set{x}$ is an independent set in $A$;
\item $A$ contains an $F$-collection covering all edges at $x$;
\item $A+xW$ contains an $F$-collection covering all edges at $x$, where $A+xW$ is obtained from $A$ by adding all edges between $x$ and $W$.
\end{enumerate}

Call $F$ \defn{$\delta$-neighbourhood-absorbing} if the following is true:
\begin{itemize}
\item[] Let $1/n \ll 1/b \ll \mu,1/|F|$ and suppose that $G$ is a graph on $n$ vertices with $\delta(G) \ge (\delta + \mu) n$. Let $x\in V(G)$ and $W\In V(G)\sm\Set{x}$ with $|W|=gcd(F)$. Then $G$ contains an $F$-neighbourhood-absorber for $(x,W)$ of order at most~$b$. 
\end{itemize}

The following result states that if the minimum degree is sufficiently large to ensure an approximate cover of the edges at $x$, a copy of $F$ covering any edge, and the existence of an $F$-neighbourhood-absorber, then we can cover \emph{all} edges at $x$ by edge-disjoint copies of~$F$.

\begin{lemma} \label{lem:general neighbourhood}
Let $F$ be a $\delta$-neighbourhood-absorbing graph and suppose that $\delta\ge\max\Set{1-1/(\chi^{vx}(F)+1),\delta_F^e}$. Then $\delta^{vx}\le \delta$.
\end{lemma}

\proof
Let $r:=gcd(F)$ and let $1/n \ll 1/b \ll 1/m' \ll \mu,1/|F|$. Suppose that $G$ is a graph on $n$ vertices with $\delta(G)\ge (\delta+5\mu)n$ and $x\in V(G)$ with $r\mid d_G(x)$. We have to show that there exists an $F$-collection covering all edges at $x$.
By Lemma~\ref{lem:get vortex}, there exists a $(\delta+4\mu,\mu,m)$-vortex $U_0 \supseteq U_1 \supseteq \dots \supseteq U_\ell$ in $G$ such that $x\in U_\ell$ and $\lfloor \mu m' \rfloor \le m \le m'$.

Let $W_1,\dots,W_s$ be an enumeration of all $r$-subsets of $U_\ell\sm\Set{x}$. We aim to find an $F$-neighbourhood-absorber for each $(x,W_i)$.

For this, let $G':=G[(U_0\sm U_1)\cup U_\ell]-G[U_\ell]$ and observe that $\delta(G')\ge (\delta+ 4\mu)|G'|$. We want to find edge-disjoint subgraphs $A_1,\dots,A_s$ in $G'$ such that $A_i$ is an $F$-neighbourhood-absorber for $(x,W_i)$ of order at most $b$. Suppose that for some $j\in[s]$, we have already found $A_1,\dots,A_{j-1}$. Let $G_j:=G'-(A_1\cup \dots \cup A_{j-1})$. Clearly, $\delta(G_j) \ge (\delta+ 3\mu)|G_j|$. Since $F$ is $\delta$-neighbourhood-absorbing, $G_j$ contains an $F$-neighbourhood-absorber for $(x,W_j)$ of order at most~$b$.

Let $G_{app}:=G-(A_1\cup \dots \cup A_s)$. Hence, $\delta(G_{app}) \ge (\delta+ 2\mu)n$. Moreover, $G_{app}[U_1]=G[U_1]$ and so $U_1 \supseteq U_2 \supseteq \dots \supseteq U_\ell$ is a $(\delta+4\mu,\mu,m)$-vortex in $G_{app}[U_1]$. Finally, since $\Delta(A_1\cup \dots \cup A_s)\le \mu |U_1|$, we have $d_{G_{app}}(y,U_1)\ge (\delta+3\mu)|U_1|$ for all $y\in V(G_{app})$. Thus, by Lemma~\ref{lem:almost cover}, there exists an $F$-collection $\cF$ covering all edges $x$ in $G_{app}$ except possibly some going to $U_\ell$. Let $W\In U_\ell$ be the set of neighbours of $x$ in $G_{app}-\bigcup \cF$. Since $r\mid d_G(x)$ and $r\mid d_{A_i}(x)$ for all $i\in[s]$, we have $r\mid |W|$. Hence, there exists a set $I\In[s]$ such that $\set{W_i}{i\in I}$ is a partition of $W$. For all $i\in I$, let $\cF_i$ be an $F$-collection covering all edges at $x$ in $A_i+xW_i$. For all $i\in [s]\sm I$, let $\cF_i$ be an $F$-collection covering all edges at $x$ in $A_i$. Then, $\cF\cup \cF_1 \cup \dots \cup \cF_s$ is an $F$-collection in $G$ covering all edges at $x$.
\endproof

\subsection{Absorbing} \label{subsec:absorb}

It remains to investigate the absorbing properties of a given graph $F$. In addition to the space barrier for $\delta_F^{vx}$ given by $1-1/(\chi^{vx}(F)+1)$, the following definition gives rise to a divisibility-type barrier. 

For a graph $F$ with $\chi(F)\ge 3$, we define
\begin{align}
\theta(F):=gcd\set{d_F(v, c^{-1}(1))-d_F(v, c^{-1}(2))}{v\in V(F),c\in Col(F,v)},\label{theta def}
\end{align}
where we set $gcd\Set{0}:=2$ for technical reasons only.

\begin{prop} \label{prop:theta extremal}
If $\theta(F)>1$ and $\chi:=\chi(F)\ge 4$, then $\delta_F^{vx}\ge 1-1/\chi$.
\end{prop}

\proof
Let $r:=gcd(F)$, $m\in\bN$ and let $G$ be the complete $\chi$-partite graph with vertex classes $V_1,\dots,V_\chi$ such that $|V_1|=rm+1$, $|V_2|=rm-1$ and $|V_i|=rm$ for all $3\le i\le \chi$. Let $x$ be a vertex in $V_\chi$. Hence, $r\mid d_G(x)$. Clearly, $\delta(G) = (1-1/\chi)|G|-1$.

Suppose that $F_1,\dots,F_t$ are edge-disjoint copies of $F$ covering all edges at $x$. Let $a_j:=d_{F_j}(x,V_1)$ and let $b_j:=d_{F_j}(x,V_3)$. So $\sum_{j\in[t]}a_j=|V_1|$ and $\sum_{j\in[t]}b_j=|V_3|$ and therefore $\sum_{j\in[t]}(a_j-b_j)=1$. However, $\theta(F)\mid (a_j-b_j)$ for all $j\in[t]$, which gives a contradiction.
\endproof

In Corollary~\ref{cor:deltafa}, we will see that $\delta_F^{vx}\le 1-1/\chi(F)$ for all graphs $F$. Hence, if $\theta(F)>1$, then this settles the problem of determining $\delta_F^{vx}$ for all graphs $F$ that are at least $4$-chromatic. The next proposition will exploit the structural information of graphs $F$ for which $\theta(F)=1$.

Let $s\in \bN$ and let $F$ be a graph. Define $CN_s(F)$ to be the set of all $(s-1)$-tuples $(a_1,\dots,a_{s-1})$ such that there exists an $[s]$-colouring $c$ of $F$ and a vertex $v\in V(F)$ such that $c(v)=s$ and $d_F(v,c^{-1}(i))=a_i$ for all $i\in[s-1]$.

Suppose we are given some graph $F$ with $v\in V(F)$ and an $[s]$-colouring $c$ of $F$ with $c(v)=s$. We say that \defn{$F',v',c'$ are obtained from $F$ by rotating $c$ around $v$} if $F'$ is obtained from $s-1$ vertex-disjoint copies $F_1,\dots,F_{s-1}$ of $F$ by identifying the copies of $v$ into a new vertex $v'$, and $c'$ is defined as follows: Let $c'(v'):=s$. For every $w'\in V(F')\sm \Set{v'}$, there is a unique $i\in[s-1]$ with $w'\in V(F_i)$. Let $w\in V(F)$ be the vertex whose role $w'$ is playing in $F_i$. Define $c'(w'):=((1,2,\dots,s-1)^i \circ c)(w)$. In other words, we permute the colours of the colour classes $c^{-1}(1),\dots,c^{-1}(s-1)$ cyclically amongst the $F_i$'s such that ultimately every $w\in V(F)\sm c^{-1}(s)$ has exactly one copy in each of the colours $1,\dots,s-1$.

Clearly, $c'$ is an $[s]$-colouring of $F'$. Moreover, for each $i\in[s-1]$, we have
\begin{align}
d_{F'}(v',c'^{-1}(i))=d_{F}(v,c^{-1}(1))+\dots+d_{F}(v,c^{-1}(s-1))=d_F(v).\label{rotated degrees}
\end{align}

\begin{prop} \label{prop:rotaters}
Let $F$ be a graph and let $\chi:=\chi(F)$. Then there exists an $F$-decomposable graph $F'$ and $m\in \bN$ such that $(m-1,m+1,m,\dots,m),(m,\dots,m)\in CN_s(F')$, with $s=\chi+1$. Moreover, if $\chi\ge 3$ and $\theta(F)=1$, then we can assume that $s=\chi$.
\end{prop}

\proof
First, suppose that $s=\chi+1$. Fix some $v\in V(F)$ and let $c$ be an $[s]$-colouring of $F$ with $c(v)=s$, with $c^{-1}(1)\cap N_F(v)\neq\emptyset$ and $c^{-1}(2)=\emptyset$. Let $m=d_F(v)$ and let $F',v',c_1$ be obtained from $F$ by rotating $c$ around $v$. Hence, $d_{F'}(v',c_1^{-1}(i))=m$ for all $i\in [s-1]$ by \eqref{rotated degrees}. Note that there is a component of $F'-v'$ in which colour $2$ does not appear and at least one neighbour of $v'$ is coloured $1$. Changing the colour of one of those neighbours to $2$ thus gives an $[s]$-colouring $c_2$ such that $d_{F'}(v',c_2^{-1}(1))=m-1$, $d_{F'}(v',c_2^{-1}(2))=m+1$ and $d_{F'}(v',c_2^{-1}(i))=m$ for all $i\in\Set{3,\dots,s-1}$. 

Now, assume that $s=\chi\ge 3$ and $\theta(F)=1$. By the definition of $\theta(F)$, there exist (not necessarily distinct) $v_1,\dots,v_t\in V(F)$ and $c_i\in Col(F,v_i)$ for $i\in[t]$ such that $$\sum_{i=1}^t (d_F(v_i,c_i^{-1}(1))-d_F(v_i,c_i^{-1}(2)))=1.$$
This implies that there are vertex-disjoint copies $F_1,\dots,F_t$ of $F$, $v_i\in V(F_i)$ and $c_i\in Col(F_i,v_i)$ for $i\in[t]$ such that $$\sum_{i=1}^t d_{F_i}(v_i,c_i^{-1}(1))=1+\sum_{i=1}^t d_{F_i}(v_i,c_i^{-1}(2)).$$ Let $F''$ be obtained by identifying $v_1,\dots,v_t$ into a new vertex $v''$. Clearly, $F''$ is $F$-decomposable. Moreover, $c_1,\dots,c_t$ induce a colouring $c''\in Col(F'',v'')$ with $d_1:=d_{F''}(v'',c''^{-1}(1))=d_{F''}(v'',c''^{-1}(2))+1$. Let $m:=d_{F''}(v'')$ and let $F',v',c_1$ be obtained from $F''$ by rotating $c''$ around $v''$. Hence, $d_{F'}(v',c_1^{-1}(i))=m$ for all $i\in [s-1]$ by~\eqref{rotated degrees}. Note that there is a collection of components of $F'-v'$ such that $d_1$ neighbours of $v'$ are coloured $1$ in the union of these components and $d_1-1$ neighbours of $v'$ are coloured $2$ in the union of these components. Hence, exchanging the colours $1$ and $2$ among all the vertices in those components gives an $[s]$-colouring $c_2$ such that $d_{F'}(v',c_2^{-1}(1))=m-1$, $d_{F'}(v',c_2^{-1}(2))=m+1$ and $d_{F'}(v',c_2^{-1}(i))=m$ for all $i\in\Set{3,\dots,s-1}$. 
\endproof

The following proposition gives a construction of a neighbourhood-absorber for $(x,W)$ where $W$ is a subset of a single class in an $s$-partite graph.

\begin{prop} \label{prop:partite abs}
Let $F$ be a graph with $r:=gcd(F)$ and assume that there exists an $F$-decomposable graph $F'$ and $m,s\in \bN$ with $s\ge 3$ such that $(m-1,m+1,m,\dots,m),(m,\dots,m)\in CN_s(F')$. Then, for every $b\in \bN$, there exists a graph $T$ and an $[s]$-colouring $c$ of $T$ such that $T$ is an $F$-neighbourhood-absorber for $(x,W)$, where $W\In c^{-1}(1)$ is of size $br$ and $x\in c^{-1}(s)$.
\end{prop}

\proof 
Let $F_1,\dots,F_t$ be vertex-disjoint copies of $F$ and $v_i\in V(F_i)$ such that $\sum_{i=1}^{t}d_{F_i}(v_i)\equiv br \mod{(s-1)m}$. Let $F''$ be obtained by identifying $v_1,\dots,v_t$ into a new vertex $x$.
Let $c''$ be an $[s]$-colouring of $F''$ such that $c''(x)=s$. We may assume that $d_{F''}(x,c^{\prime\prime-1}(1))\ge br$. Let $W\In N_{F''}(x,c^{\prime\prime-1}(1))$ be of size $br$ and let $F'''$ be obtained from $F''$ by deleting the edges from $x$ to $W$. So $F''=F'''+xW$.

For $i\in[s-1]$, let $a_i:=d_{F'''}(x,c^{\prime\prime-1}(i))$. Hence, $(s-1)m \mid \sum_{i=1}^{s-1}a_i$.
Let $\bar{a}:=(\sum_{i=1}^{s-1}a_i)/(s-1)$, let $I^+:=\set{i\in[s-1]}{a_i>\bar{a}}$ and $I^-:=\set{j\in[s-1]}{a_j<\bar{a}}$.
Define $$p:=\sum_{i\in I^+}(a_i-\bar{a})=\sum_{j\in I^-}(\bar{a}-a_j).$$  
It is easy to see that there exist non-negative integers $(b_{i,j})_{i\in I^+,j\in I^-}$ such that 
\begin{enumerate}[label=(\roman*)]
\item for each $i\in I^+$, $\sum_{j\in I^-}b_{i,j}=a_i-\bar{a}$;
\item for each $j\in I^-$, $\sum_{i\in I^+}b_{i,j}=\bar{a}-a_j$.\COMMENT{Start with $b_{i,j}=0$ for all $i,j$ and update the $b_{i,j}$'s while maintaining that $\sum_{j\in I^-}b_{i,j}\le a_i-\bar{a}$ for all $i\in I^+$ and $\sum_{i\in I^+}b_{i,j}\le \bar{a}-a_j$ for all $j\in I^-$. Now, if there exists some $i\in I^+$ such that $\sum_{j\in I^-}b_{i,j}< a_i-\bar{a}$, then $\sum_{i'\in I^+,j\in I^-}b_{i',j}<p$ and hence there also exists some $j\in I^-$ such that $\sum_{i'\in I^+}b_{i',j}< \bar{a}-a_j$. We can thus increase $b_{i,j}$ by one. After a finite number of steps, this must terminate such that $\sum_{j\in I^-}b_{i,j}= a_i-\bar{a}$ for all $i\in I^+$. Then, we also have $\sum_{i\in I^+}b_{i,j}=\bar{a}-a_j$ for each $j\in I^-$.}
\end{enumerate}
Note that
\begin{align}
\sum_{(i,j)\in I^+\times I^-}b_{i,j}=p.\label{average differences}
\end{align}

Since $(m-1,m+1,m,\dots,m) \in CN_s(F')$, there exists an $[s]$-colouring $c^\ast$ of $F'$ and $u\in V(F')$ such that $c^\ast(u)=s$, $d_{F'}(u,c^{\ast-1}(1))=m-1$, $d_{F'}(u,c^{\ast-1}(2))=m+1$ and $d_{F'}(u,c^{\ast-1}(i))=m$ for all $i\in\Set{3,\dots,s-1}$. For every $(i,j)\in I^+\times I^-$, let $F^{i,j}_1,\dots,F^{i,j}_{b_{i,j}}$ be disjoint copies of $F'$ with new vertices.  Let $T'$ be the graph obtained from $F'''$ and $(F^{i,j}_k)_{i\in I^+,j\in I^-,k\in[b_{i,j}]}$ by identifying the copy of $u$ in each $F^{i,j}_k$ with $x\in V(F''')$.

We now define a colouring $c'$ of $T'$. For every $v\in V(F''')$, let $c'(v):=c''(v)$. For every $v\in V(T')\sm V(F''')$, there are unique $i\in I^+,j\in I^-,k\in[b_{i,j}]$ such that $v\in V(F^{i,j}_k)$. Let $w$ be the vertex of $F'$ whose role $v$ is playing in $F^{i,j}_k$. Let $\phi$ be a permutation on $[s]$ such that $\phi(1)=i$, $\phi(2)=j$ and $\phi(s)=s$. Define $c'(v):=(\phi\circ c^\ast)(w)$. So $c'$ colours $V(F^{i,j}_k)$ such that $m-1$ neighbours of $x$ are coloured $i$, $m+1$ neighbours of $x$ are coloured $j$, and $m$ neighbours of $x$ are coloured $\ell$ for each $\ell\in[s-1]\sm\Set{i,j}$.
Clearly, $c'$ is an $[s]$-colouring of $T'$ with $c'(x)=s$. Moreover, for every $i\in I^+$, we have $$d_{T'}(x,c^{\prime-1}(i))\overset{\eqref{average differences}}{=}a_i+pm+\sum_{j\in I^-}-b_{i,j}\overset{(i)}{=}pm+\bar{a}$$ by the choice of $b_{i,j}$. Similarly, for all $j\in I^-$, $$d_{T'}(x,c^{\prime-1}(j))\overset{\eqref{average differences}}{=}a_j+pm+\sum_{i\in I^+}b_{i,j}\overset{(ii)}{=}pm+\bar{a}.$$ Finally, if $i\notin I^+\cup I^-$, then $a_i=\bar{a}$ and thus $d_{T'}(x,c^{\prime-1}(i))=pm+\bar{a}$ as well.

Note that $m\mid \bar{a}$ and let $q:=p+\bar{a}/m$. For every $i\in[s-1]$, let $V_{i,1},\dots,V_{i,q}$ be a partition of $N_{T'}(x,c^{\prime-1}(i))$ into parts of size $m$. For $j\in[q]$, let $U_j:=V_{1,j}\cup \dots \cup V_{s-1,j}$ and let $\hat{F}_j$ be a copy of $F'$ and $c_j$ an $[s]$-colouring of $\hat{F}_j$ such that $V(\hat{F}_j)\cap V(T')=\Set{x}\cup U_j$, such that $N_{\hat{F_j}}(x)=U_j$ and such that $c_j$ agrees with $c'$ on $\Set{x}\cup U_j$. This is possible since $(m,\dots,m)\in CN_s(F')$. We may assume that the $\hat{F}_j$'s only intersect in $x$. 

Let $T:=T'\cup \hat{F}_1 \cup \dots \cup \hat{F}_q$. Hence, $c:=c'\cup c_1 \cup \dots \cup c_q$ is an $[s]$-colouring of $T$. Moreover, $\hat{F}_1,\dots,\hat{F}_q$ cover all edges at $x$ in $T$, and $F'',(F^{i,j}_k)_{i\in I^+,j\in I^-,k\in[b_{i,j}]}$ cover all edges at $x$ in $T+xW$.\COMMENT{technically, the $F^{i,j}_k$ do not contain $x$, but the correct vertex gets identified with $x$ later.}
\endproof

\begin{lemma} \label{lem:neighbourhood abs}
Let $F$ be a graph and assume that there exists an $F$-decomposable graph $F'$ and $m,s\in \bN$ such that $(m-1,m+1,m,\dots,m),(m,\dots,m)\in CN_s(F')$, where $s\ge 4$. Then $F$ is $\max\Set{\delta_F^e,1-1/(s-1)}$-neighbourhood-absorbing.
\end{lemma}
The idea to prove Lemma~\ref{lem:neighbourhood abs} is as follows. If $\delta(G)\ge (1-1/(s-1)+\mu)|G|$, then we can find a complete $s$-partite graph $H^B$ with vertex classes $B_1,\dots,B_s$ such that $x\in B_s$. Now, if we had $W\In B_1$, then Proposition~\ref{prop:partite abs} would allow us to find the desired neighbourhood-absorber as a subgraph of $H^B$. Suppose now that $W$ is not contained in $H^B$. So our aim is to `move' $W$ to $B_1$. To achieve this, we will use the colouring properties of $F'$. (We have seen in Proposition~\ref{prop:rotaters} that such graphs $F'$ exist.) Suppose for example that $A,B,D_3,\dots,D_s$ are disjoint sets of vertices and $w\in W$, $b_w\in B$  such that $G[A\cup \Set{w},B,D_3,\dots,D_{s-1},D_s]$ is a complete $s$-partite graph and $x$ has $m$ neighbours in each of $A,B,D_3,\dots,D_{s-1}$. We can then `move' $w$ to $b_w$ as follows: If we do not need to cover $xw$, we can embed a copy of $F'$ such that all edges from $x$ to $A,B,D_3,\dots,D_{s-1}$ are covered, thereby covering $b_w$. If we intend to cover $xw$, then we can embed a copy $F'_w$ of $F'$ containing $x$ such that $d_{F'_w}(x,A\cup\Set{w})=m+1$ and $d_{F'_w}(x,B)=m-1$, leaving the edge $xb_w$ uncovered. In order to find these `movers', we will use the regularity lemma. However, for this to work, we would need that $w$ is a `typical' vertex, which we cannot assume. We will therefore use the definition of $\delta_F^e$ to `move' each $w\in W$ to some `typical' vertex first.

\lateproof{Lemma~\ref{lem:neighbourhood abs}}
Let $\delta:=\max\Set{\delta_F^e,1-1/(s-1)}$ and $r:=gcd(F)$. By Proposition~\ref{prop:partite abs}, there exist graphs $T_1,\dots,T_{|F|}$ such that $T_t$ admits an $[s]$-colouring $c_t$ and is an $F$-neighbourhood-absorber for $(z_t,W_t)$, where $W_t\In c_t^{-1}(1)$ is of size $tr$ and $z_t\in c_t^{-1}(s)$. Let $M:=r|F'|^2+\max\Set{|T_1|,\dots,|T_{|F|}|}$.

Let $1/n\ll 1/k_0' \ll \eps \ll \alpha , 1/k_0 , 1/b \ll \mu,1/|F|$. Since $M$ only depends on $F$ we may also assume that $\alpha , 1/k_0 , 1/b \ll 1/M$. Suppose that $G$ is a graph on $n$ vertices with $\delta(G) \ge (\delta + \mu) n$. Let $x\in V(G)$ and $W=\Set{w_1,\dots,w_r}\In V(G)\sm\Set{x}$. We will find an $F$-neighbourhood-absorber for $(x,W)$ in $G$ of order at most~$b$. 

Apply the regularity lemma (Lemma~\ref{lem:regularity}) to obtain a partition $V_0,V_1,\dots,V_k$ of $V(G)$ and a spanning subgraph $G'$ of $G$ satisfying the following:
\begin{enumerate}[label=(R\arabic*)]
\item $k_0 \le k \le k_0'$;
\item $|V_0| \le \eps n$;
\item $|V_1|=\dots=|V_k|=:L$;
\item $d_{G'}(z) \ge d_G(z)-(\alpha+\eps)n$ for all $z \in V(G)$;
\item $G'[V_y]$ is empty for all $y\in [k]$
\item for all $1\le y < y' \le k$, $G'[V_y,V_{y'}]$ is either $\eps$-regular with density at least $\alpha$ or empty.
\item for all $y\in [k]$, $V_y\In N_G(x)$ or $V_y\cap N_G(x)=\emptyset$.
\end{enumerate}
Let $R$ be the reduced graph of $V_1,\dots,V_k$ and let $\sigma$ be the corresponding cluster function. By Proposition~\ref{prop:reduced graph}, $\delta(R)\ge (\delta+\mu/4)k \ge (1-1/(s-1)+\mu/4)k$. 

Let $X\In [k]$ be the set of indices $y$ for which $V_y\In N_G(x)$ and let $R_x:=R[X]$. Since $|X|L\ge d_G(x)-\eps n$ and thus $|X|\ge (1-1/(s-1)) k$, we have that $\delta(R_x)\ge (1-1/(s-2)+\mu/4)|R_x|$. Therefore, there exist $y_1^\ast,\dots,y_{s-1}^\ast\in X$ and $y_s^\ast\in [k]$ such that $Y^\ast:=R[\Set{y_1^\ast,\dots,y_s^\ast}]$ is a copy of $K_s$.

Let $G''$ be obtained from $G'$ by removing, for every $i\in[r]$, all edges from $w_i$ to all clusters $V_j$ with $d_{G'}(w_i,V_j)< \alpha L$. Then we still have $\delta(G'')\ge (\delta+\mu/2)n$. Using the definition of $\delta_F^e$, it is straightforward to find copies $F_1,\dots,F_r$ of $F$ with the following properties:

\begin{enumerate}[label=(F\arabic*)]
\item $V(F_i)\sm\Set{x,w_i}\In V(G)\sm V_0$ and $xw_i\in E(F_i)$ and $E(F_i)\sm\Set{xw_i}\In E(G'')$;
\item $F_1,\dots,F_r$ intersect only in $x$.
\end{enumerate}

Let $U:=(V(F_1)\cup \dots \cup V(F_r))\sm (\Set{x}\cup W)$ and $N_x:=(N_{F_1}(x)\cup \dots \cup N_{F_r}(x))\sm W$. Let $t:=|N_x|/r$ and note that $t\le |F|$. Clearly, $\sigma(u)$ is defined for every $u\in U$, and if $u\in N_x$, then $\sigma(u)\in X$. Consider any $u\in N_x$. In $R_x$, every $s-2$ vertices have a common neighbour. Since $s\ge 4$, there exists a common neighbour $u_C$ of $\sigma(u)$ and $y_1^\ast$ in $R_x$. Moreover, there exist $u_3^+,u_3^-,\dots,u_{s-1}^+,u_{s-1}^-\in X$ and $u_s^+,u_s^-\in [k]$ such that $Y_u^+:=R[\Set{\sigma(u),u_C,u_3^+,\dots,u_s^+}]$ and $Y_u^-:=R[\Set{u_C,y^\ast_1,u_3^-,\dots,u_s^-}]$ are $s$-cliques. By the definition of $G''$, if $u\in N_{F_i}(w_i)\sm\Set{x}$, then
\begin{align}
d_{G}(w_i,V_{\sigma(u)})\ge d_{G''}(w_i,V_{\sigma(u)})\ge \alpha L/2.\label{large neighbourhood}
\end{align}\COMMENT{$L/2$ because we might have $w_iw_{i'}\in E(G')$, $d_{G'}(w_i,V_{\sigma(w_{i'})})< \alpha L$ but $d_{G'}(w_{i'},V_{\sigma(w_i)})= \alpha L$.}

Hence, by applying the key lemma (Lemma~\ref{lem:key lemma}) with suitable candidate sets, we can find a subgraph $H$ in $G''$ whose vertex set can be partitioned into sets $(A_u)_{u\in U}$, $B_1,\dots,B_s$, $(C_u)_{u\in N_x}$, $(D^+_{u,j})_{u\in N_x,j\in\Set{3,\dots,s}}$ and $(D^-_{u,j})_{u\in N_x,j\in\Set{3,\dots,s}}$ which satisfy the following (see Figure~\ref{fig:neighbourhoodabsorber}): 
\begin{enumerate}[label=(H\arabic*)]
\item for all $u\in U$, $A_u\In V_{\sigma(u)}$, for all $j\in[s]$, $B_j\In V_{y_j^\ast}$, for all $u\in N_x$ and $j\in\Set{3,\dots,s}$, $C_u\In V_{u_C}$, $D^+_{u,j}\In V_{u_j^+}$, $D^-_{u,j}\In V_{u_j^-}$;\label{set locations}
\item all those sets are independent in $H$ and $V(H)\cap (\Set{x}\cup W)=\emptyset$;
\item $H[A_u,A_{u'}]$ is complete bipartite whenever $u,u'\in U$ with $uu'\in E(F_i)$ for some $i\in[r]$;
\item $H^+_{u}:=H[A_u,C_u,D^+_{u,3},\dots,D^+_{u,s}]$ and $H^-_{u}:=H[B_1,C_u,D^-_{u,3},\dots,D^-_{u,s}]$ are complete $s$-partite for every $u\in N_x$;
\item $H^B:=H[B_1,\dots,B_s]$ is complete $s$-partite;
\item $A_u\In N_{G}(w_i)$ if $u\in N_{F_i}(w_i)\sm\Set{x}$;\label{w-neighbourhoods}
\item $|A_u|=1$ for all $u\in U\sm N_x$;
\item all sets $(A_u)_{u\in N_x}$, $C_u$, $D^+_{u,j}$, $D^-_{u,j}$ have cardinality $|F'|$;
\item $|B_j|=M$ for $j\in[s]$.
\end{enumerate}
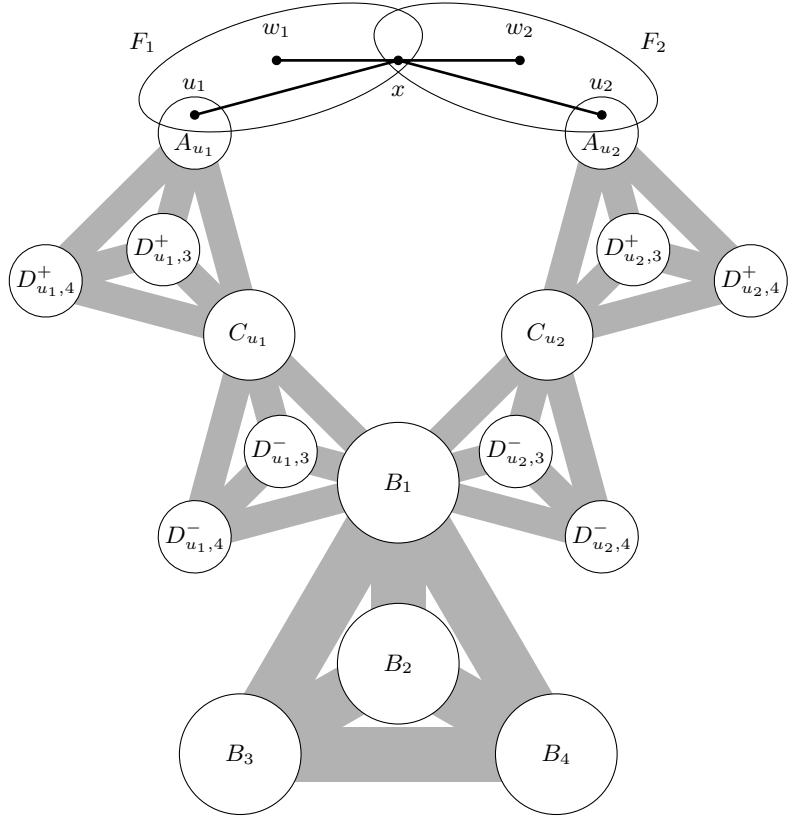
\begin{figure}[htbp]
\begin{center}
\begin{tikzpicture}[scale=0.8]
{\footnotesize
\begin{scope}

		\coordinate (x) at (0,7);
		\filldraw[fill=black]  (x) circle (2pt);
		\node at ($(x)+(0,-0.5)$) {$x$};

	\begin{scope}[rotate =-15, shift={(-2,0)} ]

		\foreach \i in {0,120,240}
		{
				\fill [fill=black!30,rotate=\i] (0,0.25) rectangle (2,-0.25);
				\fill [fill=black!30,rotate = \i -60] ($ (-60:2 ) + (0.25,0)$) rectangle ($(+60:2 ) + (-0.25,0)$);
		}
		\foreach \i in {0,120,240}
		{
				\draw [fill = white] (\i:2) circle (0.6);
		}
		\draw [fill = white] (0,0) circle (0.6);
		
		\node  at ( $(0,0) $)  {$D^-_{u_1,3}$};
		\node  at ( $(-120:2)  $)  {$D^-_{u_1,4}$};
	\end{scope}		
	
	\begin{scope}[rotate =-15, shift={($(-2,0)+(120:2)$)}, rotate =-30, shift={((-2,0))}]

		\foreach \i in {0,120,240}
		{
				\fill [fill=black!30,rotate=\i] (0,0.25) rectangle (2,-0.25);
				\fill [fill=black!30,rotate = \i -60] ($ (-60:2 ) + (0.25,0)$) rectangle ($(+60:2 ) + (-0.25,0)$);
		}
		\foreach \i in {0,120,240}
		{
				\draw [fill = white] (\i:2) circle (0.6);
		}
		\draw [fill = white] (0,0) circle (0.6);
		
		\node  at ( $(0,0) $)  {$D^+_{u_1,3}$};
		\node  at ( $(-120:2)  $)  {$D^+_{u_1,4}$};
		
		\draw [fill = white] (0:2) circle (0.75);
		\node  at ( $(0:2)  $)  {$C_{u_1}$};
		\coordinate (a1) at (120:2);
	\end{scope}

		\node  at ( $(a1) -(0,0.2)$)  {$A_{u_1}$};
		\filldraw[fill=black]  ( $(a1)+(0,0.3)$) circle (2pt);
		\node  at ( $(a1)+(0,0.8)$)  {${u_1}$};

		\coordinate (b1) at (-2,7);
		\filldraw[fill=black]  (b1) circle (2pt);
		\node  at ( $(b1)+(0,0.5)$)   {${w_1}$};
		\draw [line width = 1pt]  ( $(a1)+(0,0.3)$) -- (x) -- (b1);
	
		\draw[shift={($(x)+(0,0.4)$)}, rotate=15] (-2,0) ellipse (2.4 and 0.9);
		\node  at ( -4.2,7.3)  {$F_1$};
	
\end{scope}

\begin{scope}[xscale=-1]

	\begin{scope}[rotate =-15, shift={(-2,0)} ]

		\foreach \i in {0,120,240}
		{
				\fill [fill=black!30,rotate=\i] (0,0.25) rectangle (2,-0.25);
				\fill [fill=black!30,rotate = \i -60] ($ (-60:2 ) + (0.25,0)$) rectangle ($(+60:2 ) + (-0.25,0)$);
		}
		\foreach \i in {0,120,240}
		{
				\draw [fill = white] (\i:2) circle (0.6);
		}
		\draw [fill = white] (0,0) circle (0.6);
		
		\node  at ( $(0,0) $)  {$D^-_{u_2,3}$};
		\node  at ( $(-120:2)  $)  {$D^-_{u_2,4}$};
	\end{scope}		
	
	\begin{scope}[rotate =-15, shift={($(-2,0)+(120:2)$)}, rotate =-30, shift={((-2,0))}]

		\foreach \i in {0,120,240}
		{
				\fill [fill=black!30,rotate=\i] (0,0.25) rectangle (2,-0.25);
				\fill [fill=black!30,rotate = \i -60] ($ (-60:2 ) + (0.25,0)$) rectangle ($(+60:2 ) + (-0.25,0)$);
		}
		\foreach \i in {0,120,240}
		{
				\draw [fill = white] (\i:2) circle (0.6);
		}
		\draw [fill = white] (0,0) circle (0.6);
		
		\node  at ( $(0,0) $)  {$D^+_{u_2,3}$};
		\node  at ( $(-120:2)  $)  {$D^+_{u_2,4}$};
		
		\draw [fill = white] (0:2) circle (0.75);
		\node  at ( $(0:2)  $)  {$C_{u_2}$};
		\coordinate (a2) at (120:2);
	\end{scope}

		\node  at ( $(a2) -(0,0.2)$)  {$A_{u_2}$};
		\filldraw[fill=black]  ( $(a2)+(0,0.3)$) circle (2pt);
		\node  at ( $(a2)+(0,0.8)$)  {${u_2}$};

		\coordinate (b2) at (-2,7);
		\filldraw[fill=black]  (b2) circle (2pt);
		\node  at ( $(b2)+(0,0.5)$)   {${w_2}$};
		\draw [line width = 1pt]  ( $(a2)+(0,0.3)$) -- (x) -- (b2);
	
		\draw[shift={($(x)+(0,0.4)$)}, rotate=15] (-2,0) ellipse (2.4 and 0.9);
		\node  at ( -4.2,7.3)  {$F_2$};
	
\end{scope}

\begin{scope}[shift={(0,-3)}]

		\foreach \i in {90,210,330}
		{
				\fill [fill=black!30,rotate=\i] (0,0.45) rectangle (3,-0.45);
				\fill [fill=black!30,rotate = \i -60] ($ (-60:3) + (0.45,0)$) rectangle ($(+60:3 ) + (-0.45,0)$);
		}
		\foreach \i in {90,210,330}
		{
				\draw [fill = white] (\i:3) circle (1);
		}
		\draw [fill = white] (0,0) circle (1);

		\node  at (90:3) {$B_1$};
		\node  at (0,0)  {$B_2$};
		\node  at (210:3) {$B_3$};
		\node  at (330:3)  {$B_4$};

\end{scope}
}
\end{tikzpicture}
\end{center}
\caption{An example illustrating the construction of $H$ for $F=C_5$, $r=2$, $s=4$. Note that $N_x=\Set{u_1,u_2}$. The four singleton sets $A_u$ for $u\in U\sm N_x$ are not shown in the figure.}
\label{fig:neighbourhoodabsorber}
\end{figure}
So $H_u^+$ arises from the clique $Y_u^+$, $H_u^-$ arises from $Y_u^-$, and $H^B$ arises from $Y^\ast$. Note that to ensure \ref{w-neighbourhoods} we use \eqref{large neighbourhood}. Also, we do not require that $u\in A_u$.

We now describe how to construct a neighbourhood-absorber $A$ from $H$ by attaching $x$ and $W$ in a suitable way. Note that \ref{set locations} implies that
\begin{enumerate}[label=(H10)]
\item $(B_{j})_{j\in[s-1]}$, $(A_u)_{u\in N_x}$, $(C_u)_{u\in N_x}$, $(D^+_{u,j})_{u\in N_x,j\in\Set{3,\dots,s-1}}$, $(D^-_{u,j})_{u\in N_x,j\in\Set{3,\dots,s-1}}$ are contained in $N_G(x)$.\label{locations}
\end{enumerate}
For every $u\in N_x$ and $j\in\Set{3,\dots,s-1}$, let $A_u'\In A_u$, $C_u'\In C_u$, $D^{+\prime}_{u,j}\In D^+_{u,j}$, $D^{-\prime}_{u,j}\In D^-_{u,j}$ be such that $|A_u'|=|D^{+\prime}_{u,j}|=|D^{-\prime}_{u,j}|=m$ and $|C_u'|=2m$. Moreover, let $a_u$ be some element of $A_u'$. For every $u\in U\sm N_x$, let $a_u$ be the unique element of $A_u$. Let $(B'_u)_{u\in N_x}$ be disjoint subsets of $B_1$ of size $m+1$ and fix some $b_u\in B'_u$ for each $u\in N_x$. Define $W':=\set{b_u}{u\in N_x}$. So $|W'|=|N_x|=tr$. Let $T''$ be a subgraph of $H^B$ and $z\in B_s$ such that $T''$ is an $F$-neighbourhood-absorber for $(z,W')$ and $N_{T''}(z)\cap B'_u=\emptyset$ for all $u\in N_x$. (So we might choose $T''$ to be a copy of the graph $T_t$ defined at the beginning of the proof.)

Let $A$ be the graph obtained from $H$ by adding $x$ and $W$ and the following edges: For every $i\in[r]$, add the edge from $w_i$ to $a_u$ for all $u\in N_{F_i}(w_i)\sm\Set{x}$. For every $u\in N_x$ and $j\in\Set{3,\dots,s-1}$, add all edges from $x$ to $A_u',C_u',D^{+\prime}_{u,j},D^{-\prime}_{u,j},B'_u$. Moreover, add all edges from $x$ to $N_{T''}(z)$. By \ref{locations} and \ref{w-neighbourhoods}, $A$ is a subgraph of $G$. Moreover, $A$ has order at most $b$ and $\Set{x}\cup W$ is independent in $A$. Let $T'$ be obtained from $T''$ by replacing $z$ with $x$. So $T'$ is an $F$-neighbourhood-absorber for $(x,W')$.

We claim that $A$ is an $F$-neighbourhood-absorber for $(x,W)$. Let $C'_{u,1},C'_{u,2}$ be a partition of $C_u'$ into two sets of size $m$ and let $C'_{u,+},C'_{u,-}$ be a partition of $C_u'$ into two sets of sizes $m+1,m-1$.

The edges at $x$ in $A$ can be covered by edge-disjoint copies of $F$ as follows: Let $u\in N_x$. Since $(m,\dots,m)\in CN_s(F')$, there exists a copy $F'_{u,+}$ of $F'$ in $H^+_{u}$ such that $N_{F'_{u,+}}(v)=A_u'\cup C_{u,1}' \cup D^{+\prime}_{u,3} \cup \dots \cup D^{+\prime}_{u,s-1}$ for some $v\in D^{+\prime}_{u,s}$. Exchanging $v$ with $x$ yields a copy of $F'$ that covers all edges from $x$ to $N_{F'_{u,+}}(v)$ and otherwise uses only edges inside $H^+_{u}$. Similarly, there exists a copy $F'_{u,-}$ of $F'$ that covers all edges from $x$ to $(B_u'\sm\Set{b_u})\cup C_{u,2}' \cup D^{-\prime}_{u,3} \cup \dots \cup D^{-\prime}_{u,s-1}$ and otherwise uses only edges inside $H^-_{u}$. This can be done for all $u\in N_x$ without interference. This way, all edges at $x$ in $A$ are covered except the ones that have an endpoint in $N_{T'}(x)\cup\bigcup_{u\in N_x}\Set{b_u}=N_{T'}(x)\cup W'$. Finally, these edges can be covered by edge-disjoint copies of $F$ since $T'$ is an $F$-neighbourhood-absorber for $(x,W')$.

The edges at $x$ in $A+xW$ can be covered by edge-disjoint copies of $F$ as follows: Let $U_i:=V(F_i)\sm\Set{x,w_i}$. Then $F_i':=A[\Set{x,w_i}\cup\bigcup_{u\in U_i}\Set{a_u}]+xw_i$ is a copy of $F$. Moreover, $F_1',\dots,F_r'$ are edge-disjoint subgraphs of $A+xW$ and cover the edges from $x$ to $W\cup \bigcup_{u\in N_x}\Set{a_u}$.
Let $u\in N_x$. Since $(m-1,m+1,m,\dots,m)\in CN_s(F')$, there exists a copy of $F'$ that covers all edges from $x$ to $(A_u'\sm\Set{a_u})\cup C_{u,+}' \cup D^{+\prime}_{u,3} \cup \dots \cup D^{+\prime}_{u,s-1}$ and otherwise uses only edges inside $H^+_{u}$, and a copy of $F'$ that covers all edges from $x$ to $C_{u,-}' \cup B_u' \cup D^{-\prime}_{u,3} \cup \dots \cup D^{-\prime}_{u,s-1}$ and otherwise uses only edges inside $H^-_{u}$. Finally, by definition of $T'$, the edges at $x$ in $T'$ can be covered by edge-disjoint copies of $F$.
\endproof

The next lemma analyses the bipartite case. We want to establish that $\delta_F^{vx}=0$ if $F$ contains a bridge, and $\delta_F^{vx}=1/2$ otherwise. The results leading to the near-optimal cover include the bipartite case.\COMMENT{In fact the proofs become trivial at many stages.} It only remains to determine the absorbing properties of a given bipartite graph $F$.

\begin{lemma} \label{lem:bipartite neighbourhood}
Let $F$ be bipartite and $r:=gcd(F)$. Then $F$ is $1/2$-neighbourhood-absorbing. Moreover, if $F$ contains a bridge, then $F$ is $0$-neighbourhood-absorbing.
\end{lemma}

\proof
Let $uz\in E(F)$ and assume that this edge is a bridge if one exists. Let $d^\ast:=d_F(u)$. There exist (not necessarily distinct) vertices $v_1,\dots,v_t\in V(F)$ such that $D:=\sum_{i=1}^t d_F(v_i)\equiv -r \mod{d^\ast}$. We may assume that $D\ge r(d^\ast-1)$. Moreover, let $\delta:=1/2$ if $F$ contains no bridge and $\delta:=0$ otherwise.

Let $1/n\ll 1/b \ll \mu,1/|F|$. Since $D$ only depends on $F$, this means that $1/b\ll 1/D$. Suppose that $G$ is a graph on $n$ vertices with $\delta(G) \ge (\delta + \mu) n$. Let $x\in V(G)$ and $W=\Set{w_1,\dots,w_r}\In V(G)\sm\Set{x}$. We will find an $F$-neighbourhood-absorber for $(x,W)$ in $G$ of order at most~$b$.

Suppose that $\delta=1/2$. We can use the regularity lemma to find disjoint sets $B_1,\dots,B_{r},B_x$ of size $(D+r)|F|$ each in $V(G)\sm({\Set{x}\cup W})$ such that $B_x\In N_G(x)$ and $G[B_i,B_x]$ is complete bipartite and $B_i\In N_G(w_i)$ for all $i\in[r]$.

Let $A$ be the graph obtained from $\bigcup_{i=1}^{r}G[B_i,B_x]$ by adding $\Set{x}\cup W$, all edges from $w_i$ to $B_i$ for all $i\in[r]$, and exactly $D$ edges from $x$ to $B_x$. Clearly, $\Set{x}\cup W$ is independent in $A$ and $A$ has order at most $r+1+(r+1)(D+r)|F|\le b$.

Moreover, the edges at $x$ in $A$ can be covered by edge-disjoint copies of $F$. Indeed, since $d_{A}(x)=D=\sum_{i=1}^t d_F(v_i)$, we can let $x$ play the roles of $v_1,\dots,v_t$ and use edges of $A[B_x,B_1]$ otherwise.

In order to see that all edges at $x$ in $A+xW$ can be covered, partition $N_A(x)$ into sets $N_1,\dots,N_{r+1}$, where $|N_i|=d^\ast-1$ for all $i\in[r]$. It follows that $d^\ast\mid |N_{r+1}|$.
Now, for every $i\in[r]$, we can cover the edge $xw_i$ with a copy $F_i$ of $F$ such that $x,w_i$ play the roles of $u,z$ and $N_{F_i}(x)=N_i\cup\Set{w_i}$. The edges from $x$ to $N_{r+1}$ can also be covered by letting $x$ play the role of $u$ in every copy, and using edges of $A[B_x,B_1]$ otherwise.

The case $\delta=0$ is very similar. Using the \Erd-Stone theorem, we can find disjoint sets $B_1,\dots,B_r,B_x$ and $B_1',\dots,B_{r}',B_x'$ of size $(D+r)|F|$ each in $V(G)\sm({\Set{x}\cup W})$ such that $B_i\In N_G(w_i)$ and $G[B_i,B_i']$ is complete bipartite for each $i\in[r]$, and $B_x\In N_G(x)$ and $G[B_x,B_x']$ is complete bipartite. Let $A$ be the graph obtained from $G[B_x,B_x']\cup\bigcup_{i=1}^{r}G[B_i,B_i']$ by adding $\Set{x}\cup W$, all edges from $w_i$ to $B_i$ for all $i\in[r]$, and exactly $D$ edges from $x$ to $B_x$. That $A$ is the desired neighbourhood-absorber follows analogously to the above case.
\endproof

We can now combine our previous results. We make use of the following simple bounds on~$\delta_F^e$.

\begin{fact}\label{fact:delta-trivial}
Let $F$ be a graph. Then
\begin{enumerate}[label=(\roman*)]
\item $\delta_F^e\le 1-1/\chi(F)$;\COMMENT{Not sure how obvious that really is.}
\item $\delta_F^e=0$ if $F$ is bipartite and contains a bridge.
\end{enumerate}
\end{fact}

Indeed (i) and~(ii) follow easily using the regularity lemma and the key lemma. In general, it seems very difficult to give an explicit formula for $\delta_F^e$. This also seems an interesting problem in its own right. Recall that $\theta(F)$ was defined in \eqref{theta def}.

\begin{cor} \label{cor:deltafa}
Let $F$ be a graph with $\chi:=\chi(F)$ and $\chi^{vx}:=\chi^{vx}(F)$.
\begin{enumerate}[label=(\roman*)]
\item $\delta_F^{vx}\le 1-1/\chi$;
\item if $F$ is bipartite, then $\delta_F^{vx}=\begin{cases} 0 &\mbox{if }F \mbox{ contains a bridge}; \\
1/2 & \mbox{otherwise;} \end{cases}$
\item if $\chi\ge 4$, then $\delta_F^{vx} = \begin{cases} 1-1/\chi &\mbox{if } \theta(F)>1; \\
\max\Set{1-1/(\chi^{vx}+1),\delta_F^e} &\mbox{if } \theta(F)=1. \end{cases}$
\end{enumerate}
\end{cor}

\proof
By Proposition~\ref{prop:rotaters} and Lemmas~\ref{lem:neighbourhood abs} and~\ref{lem:bipartite neighbourhood}, $F$ is $\max\Set{\delta_F^e,1-1/\chi}$-neighbourhood-absorbing. Then, Lemma~\ref{lem:general neighbourhood} implies that $\delta_F^{vx}\le \max\Set{1-1/\chi,1-1/(\chi^{vx}+1),\delta_F^e}$.
By Fact~\ref{fact:delta-trivial}(i) and \eqref{eq:chi_vx relation}, the general upper bound for $\delta_F^{vx}$ stated in (i) follows.

Moreover, if every edge of $F$ is contained in a cycle, it is easy to see that $\delta_F^{vx}\ge\delta_F^e\ge 1/2$, where an extremal example consists of two disjoint cliques with one edge joining the cliques. Hence, if $F$ is bipartite, then $\delta_F^{vx}=1/2$ if $F$ contains no bridge. Otherwise, by Lemma~\ref{lem:bipartite neighbourhood}, Lemma~\ref{lem:general neighbourhood} and Fact~\ref{fact:delta-trivial}(ii), $\delta_F^{vx}=0$.

Now, if $\chi\ge 4$ and $\theta(F)>1$, then $\delta_F^{vx}=1-1/\chi$ by Proposition~\ref{prop:theta extremal} and (i), so suppose that $\theta(F)=1$. Then, by Proposition~\ref{prop:rotaters} and Lemma~\ref{lem:neighbourhood abs}, $F$ is $\max\Set{\delta_F^e,1-1/(\chi-1)}$-neighbourhood-absorbing. Lemma~\ref{lem:general neighbourhood} thus implies that $$\delta_F^{vx}\le \max \Set{\delta_F^e,1-1/(\chi-1),1-1/(\chi^{vx}+1)}\overset{\eqref{eq:chi_vx relation}}{=} \max\Set{\delta_F^e,1-1/(\chi^{vx}+1)}.$$ Hence, by Proposition~\ref{prop:space extremal}, we have $\delta_F^{vx}=\max\Set{\delta_F^e,1-1/(\chi^{vx}+1)}$.
\endproof

We are now ready to prove Theorem~\ref{thm:main}.

\lateproof{Theorem~\ref{thm:main}}
Firstly, note that (i) follows from Theorem~\ref{thm:almost main} and Corollary~\ref{cor:deltafa}(i).

To prove (ii), suppose that $\chi\ge 5$. By Theorem~\ref{thm:almost main}, we have $\delta_F\in \Set{\max\Set{\delta_F^{0+},\delta_F^{vx}},1-1/\chi,1-1/(\chi+1)}$.
By Corollary~\ref{cor:deltafa}(iii), we have $\delta_F^{vx}\in \Set{1-1/(\chi^{vx}(F)+1),\delta_F^e,1-1/\chi}$. Since $\delta_F^{0+}\ge 1-1/(\chi^{vx}(F)+1)$ by Proposition~\ref{prop:space extremal} and since $\delta_F^{vx}\ge \delta_F^e$, it follows that $\max\Set{\delta_F^{0+},\delta_F^{vx}}\in \Set{\max\Set{\delta_F^{0+},\delta_F^{e}},1-1/\chi}$, implying (ii).

Finally, (iii) follows from (ii) and Corollary~\ref{cor:threshold relations}.
\endproof

\section{Concluding remarks}\label{sec:conclusion}

We conclude this paper with some final remarks. In the light of Theorem~\ref{thm:main}, for all graphs $F$ with $\chi(F)\ge 3$ the limiting factor in giving good explicit bounds on $\delta_F$ are now the bounds available in the literature for $\delta^\ast_F$.
The original aim of this project was to be able to determine from the value of $\delta^\ast_F$ the value of $\delta_F$. This we come close to achieving when $\chi(F)\ge 5$, showing that $\delta_F$ is either $\delta^\ast_F$ or one of only two other values. We note that, in order to determine which of these values $\delta_F$ takes, it is left only to determine the minimum $d$ for which there exists an augmented $d$-compressible $(C_4)_F$-switcher and an augmented $d$-compressible $(K_{2,gcd(F)})_F$-switcher (see Section~\ref{sec:switchers}). Furthermore, in the light of the proof of Lemma~\ref{lem:discretisation}, it is sufficient to ask what the minimum $d$ is such that the following holds:
There exists some $n_0$ such that any $F$-divisible balanced $d$-partite graph with at least $dn_0$ vertices and which is missing at most $e(F)^2$ edges (say) between vertex classes is $F$-decomposable.
In other words, if $\delta_F$ is not equal to $\delta^\ast_F$, then there will exist extremal graphs which are extremely close to large complete $\chi(F)$-partite or $(\chi(F)+1)$-partite graphs.

Finally, we briefly consider the case when $\chi:=\chi(F)\in\Set{3,4}$. When $\chi\ge 5$, we reduced finding an $F$-decomposition to constructing two different augmented switchers, before showing that the smallest minimum degree ratio above which we can construct these switchers takes one of three values: $1-1/(\chi-1)$, $1-1/\chi$ or $1-1/(\chi+1)$ (see Section~\ref{sec:switchers}). This result concerning augmented switchers also holds when $\chi\in\Set{3,4}$, but $\delta_F$ may not be large enough to perform the reductions from augmented switchers to unaugmented switchers. We have no version of the `discretisation lemma' (Lemma~\ref{lem:discretisation}) for unaugmented switchers, and therefore get no `discretisation result' for $\delta_F$ here, only an upper bound (Theorem~\ref{thm:main}(i)). Furthermore, it seems likely that the smallest minimum degree ratio above which these unaugmented switchers appear can take many different values outside of $\{1-1/(\chi-1), 1-1/\chi, 1-1/(\chi+1)\}$, and therefore we do not expect a simple `discretisation lemma' to hold in this case. Due to this, we suspect that there exist graphs $F$ with $\chi(F)=3$ and graphs $F$ with $\chi(F)=4$ for which Theorem~\ref{thm:main}(ii)
and~(iii) do not hold, but to show this there remains much work to do in giving good bounds on~$\delta_F^\ast$.

\vspace{2cm}

{\footnotesize \obeylines \parindent=0pt

Stefan Glock, Daniela K\"{u}hn, Allan Lo, Richard Montgomery, Deryk Osthus
\vspace{0.3cm}
School of Mathematics
University of Birmingham
Edgbaston
Birmingham
B15 2TT
UK

\vspace{0.3cm}
\begin{flushleft}
{\it{E-mail addresses}:}
\tt{[s.glock,d.kuhn,s.a.lo,r.h.montgomery,d.osthus]@bham.ac.uk}
\end{flushleft}

\end{document}